\documentclass[a4paper,12pt,reqno]{amsart}

\usepackage[all]{xy}
\xyoption{web}
\usepackage{amssymb}

\numberwithin{equation}{section}

\theoremstyle{plain}
\newtheorem{thm}{Theorem}[section]
\newtheorem*{thm*}{Theorem}

\newtheorem{lem}[thm]{Lemma}
\newtheorem{prop}[thm]{Proposition}
\newtheorem{conj}[thm]{Conjecture}

\theoremstyle{remark}
\newtheorem{rem}[thm]{Remark}

\theoremstyle{definition}

\def\U{{\mathcal U_q(\sl_2)}}
\def\Uhat{{\mathcal U_q(\widehat\sl_2)}}

\def\C{\mathbb C}

\def\Z{\mathbb Z}

\def\o{\otimes}

\def\dual{\star}

\def\adm#1{\overline{#1}}

\def\ZZ{\mathcal Z}

\def\ex{\noindent{\bf Example. }}

\def\End{{\operatorname{End}}}
\def\Res{{\operatorname{Res}}}
\def\Adm{{\operatorname{Adm}}}
\def\Id{{\operatorname{Id}}}
\def\Tr{{\operatorname{Tr}}}

\def\Opp{{\operatorname{opp}}}

\def\Groth{\mathfrak{Gr}}
\def\Path{\operatorname{Path}}

\def\Fun{{\mathfrak {Fun}}}
\def\FunBA{{\mathfrak{\Fun}}}
\def\Harm{{\mathfrak{Harm}}}

\def\Ver{{\operatorname{Ver}}}
\def\Span{{\operatorname{Span}}}
\def\Conf{{\operatorname{Conf}}}
\def\cb{\vartheta}

\def\Re{\operatorname{Re}}
\def\Im{\operatorname{Im}}

\def\1{{\mathbf 1}}

\def\a{\alpha}
\def\z{\zeta}
\def\eps{\varepsilon}

\def\F{{\mathcal F}}

\def\sl{\mathfrak {sl}}

\def\bh{{\mathbf h}}

\def\h{{\mathfrak h}}

\def\m{{\mathfrak m}}
\def\<{\langle}
\def\>{\rangle}
\def\({\left(}
\def\){\right)}
\def\[{\left[}
\def\]{\right]}

\def\qqR{\boldsymbol{\mathcal R}}
\def\qqRch{\boldsymbol{\check{\mathcal R}}}
\def\qqKZ{\boldsymbol{\mathcal K}}
\def\qqH{\boldsymbol{\mathcal H}}

\def\qR{\boldsymbol{\rm R}}
\def\qRch{\boldsymbol{\check{\rm R}}}
\def\qKZ{\boldsymbol{\rm K}}

\def\dR{\boldsymbol{\mathbb R}}
\def\dRch{\boldsymbol{\check{\mathbb R}}}
\def\qKZB{\boldsymbol{\mathbb K}}
\def\dH{\boldsymbol{\mathbb H}}

\def\qbinom#1#2{{\begin{bmatrix} #1 \\ #2 \end{bmatrix}}}
\def\hyper{{{}_{{}_3} \boldsymbol{\digamma}_{{}_2}}}
\def\hyp{{{}_{{}_2} \boldsymbol{\digamma}_{{}_1}}}

\begin{document}

\title[Resonance relations]
{\strut\\[-3.3\baselineskip]
Resonance relations, holomorphic trace functions \\
and hypergeometric solutions to qKZB and Macdonald-Ruijsenaars
equations}

\author[K. Styrkas and A. Varchenko]
{K. Styrkas$^{\,\diamond}$ and A. Varchenko$^{\,\star}$}

\maketitle

\begin{center}
{\it $^\diamond$International University Bremen, Bremen 28759,
Germany

\medskip
$^\star$Department of Mathematics,
University of North Carolina at Chapel Hill\\
Chapel Hill, NC 27599-3250, USA }
\end{center}

\begin{abstract}

The resonance relations are identities between coordinates of
functions $\psi(\lambda)$ with values in  tensor products of
representations of the quantum group $\U$. We show that the space of
hypergeometric solutions of the associated qKZB equations is
characterized as the space of functions of Baker-Akhiezer type,
satisfying the resonance relations. We give an alternative
representation-theoretic construction of this space, using the
traces of regularized intertwining operators for the quantum group
$\U$, and thus establish the equivalence between hypergeometric and
trace function solutions of the qKZB equations.

We define the quantum conformal blocks as distinguished Weyl
anti-invariant hypergeometric qKZB solutions with values in a tensor
product of finite-dimensional $\U$-modules. We prove that for
generic $q$ the dimension of the space of quantum conformal blocks
equals the dimension of $\U$-invariants, and when $q$ is a root of
unity is computed by the Verlinde algebra.
\end{abstract}




\tableofcontents

\section{Introduction}

\medskip
\noindent
1.1.\
The motivation for this paper was to establish correspondence
between holomorphic solutions to the quantum
Knizhnik-Zamolodchikov-Bernard (qKZB) equations, arising from two
different constructions: as trace functions of regularized
intertwining operators for the quantum group $\U$ and as
hypergeometric solutions, given in terms of contour integrals.

In this paper we consider the qKZB operators without spectral
parameters, which act in the space $\Fun \o M_{\vec\Lambda}[0]$,
where $\Fun$ is the space of functions of a complex variable
$\lambda$, and $M_{\vec\Lambda}[0]$ is the zero weight subspace of a
tensor product $M_{\Lambda_1} \o \dots \o M_{\Lambda_n}$ of Verma
modules for the quantum group $\U$. The qKZB operators $\qKZB_j$ are
expressed in terms of the dynamical $R$-matrix operators
$\dR_{\Lambda_i,\Lambda_k}(\lambda) \in \End(M_{\Lambda_i} \o
M_{\Lambda_k})$ and act on $M_{\vec\Lambda}[0]$-valued functions of
$\lambda$:
\begin{equation}\label{intro:qKZB operator}
\begin{split}
    \qKZB_j =& \( \dR_{\Lambda_{j+1},
\Lambda_j}(\lambda-\bh^{(j+2,\dots,n)}) \)^{-1} \dots
    \( \dR_{\Lambda_n,\Lambda_j}(\lambda) \)^{-1}
    \ \Gamma_j
    \\
    & \dR_{\Lambda_j,\Lambda_1}(\lambda-\bh^{(2,\dots,j-1)}
-\bh^{(j+1,\dots,n)})  \dots
    \dR_{\Lambda_j,\Lambda_{j-1}}(\lambda-\bh^{(j+1,\dots,n)}).
\end{split}
\end{equation}
where $\Gamma_j \, \psi(\lambda) = \psi(\lambda-\bh^{(j)})$, and
$\bh^{(i,\dots,k)}$ must be replaced with $\mu_i+\dots+\mu_k$ when
acting on homogeneous functions
 $\psi(\lambda)$ with values in $ M_{\Lambda_1}[\mu_1] \o \dots \o
M_{\Lambda_n}[\mu_n]$. Thus defined operators $\qKZB_j$ pairwise
commute for $j=1,\dots,n$, and the corresponding qKZB equations are
the equations for common eigenfunctions of operators $\qKZB_j$.

A distinguished hypergeometric family of solutions to the general
qKZB equations, associated with Felder's elliptic dynamical
$R$-matrix with spectral parameter, was introduced in \cite{FTV
monodromy} in terms of contour integrals. Our operators
$\dR_{\Lambda_i,\Lambda_k}(\lambda)$ are the asymptotic limits of
this elliptic $R$-matrix, and hence a suitable limit of the integral
formulae in \cite{FTV monodromy} yields a family of holomorphic
eigenfunctions for the operators $\qKZB_j$. Coordinates of these
vector-valued eigenfunctions are arranged in the hypergeometric qKZB
matrix $\dH(\lambda,x;\vec\Lambda)$, where $x$ is an additional
complex parameter. For a fixed $x\in\C$ the matrix can be regarded
as an operator $\dH(\lambda,x;\vec\Lambda): M_{\vec\Lambda}[0] \to
\Fun \o M_{\vec\Lambda}[0]$, which satisfies  the equations
\begin{equation}\label{eq:intro qKZB equations}
    \qKZB_j \ \dH(\lambda,x;\vec\Lambda) =
    \dH(\lambda,x;\vec\Lambda) \ \mathcal E_j(x),
    \qquad j=1,\dots,n,
\end{equation}
for suitable diagonal operators $\mathcal E_j(x)$. One can show that
the hypergeometric construction yields all qKZB solutions of certain
type, cf. Theorem \ref{thm:completeness of hypergeometric qKZB}.

\medskip
\noindent
1.2.\
Solutions $\psi(\lambda)$ of the qKZB equations, arising from the
hypergeometric construction, satisfy the so-called resonance
relations \cite{FV}, which are algebraic identities between values
of coordinate functions of $\psi(\lambda)$ at certain values of
$\lambda$.

In the simplest case $n=1$, the resonance relations first appeared
in the study of algebraic integrability of Schr\"odinger operators
with the Calogero-Sutherland potential as a system of axioms on a
Baker-Akhiezer function $\psi(\lambda,x) = e^{\lambda x}
P(\lambda,x)$, where $P(\lambda,x)$ is a monic polynomial in the
first variable of degree $\m$. In particular, it was shown in
\cite{CV} that for generic $x \in \C$ the one-point resonance
relations
\begin{equation*}
    \psi(-\delta,x) = \psi(\delta,x), \qquad \delta = 1,\dots,\m,
\end{equation*}
determine such a function uniquely, and that $\psi(\lambda,x)$ is an
eigenfunction for the Calogero-Sutherland differential operator in
the $x$ variable. This fact has a $q$-difference analogue: a
function of the form $\psi(\lambda,x) = q^{\lambda x}
P(q^{2\lambda},x)$, satisfying the same conditions, exists and is
unique;  moreover this function is an eigenfunction of the Macdonald
difference operator in the $x$ variable. In both cases, the function
 $\psi(\lambda,x)$ also satisfies a
difference equation with respect to $\lambda$.

The solutions of the qKZB equations are parameterized by a complex
parameter $x$ and a finite set of multi-indices $\vec m =
(m_1,\dots,m_n)$. The  solutions have the trigonometric
quasi-polynomial form
\begin{equation}\label{intro:n>1 BA type vector function}
    \psi(\lambda) = q^{\lambda (x-\m)}
    \( \sum_{k=0}^\m  q^{2 k \lambda} \, \psi^{(k)} \),
    \qquad
    \psi^{(k)} \in M_{\vec\Lambda}[0]\ .
\end{equation}
Here the "level" $\m$ is the nonnegative integer, determined by the
equation $\Lambda_1+\dots+\Lambda_n = 2\m$.

In this paper we regard the $n$-point resonance relations of
\cite{FV} as a system of axioms on a vector-valued trigonometric
quasi-polynomial function $\psi(\lambda)$ of the form
\eqref{intro:n>1 BA type vector function}. A detailed combinatorial
description of the $n$-point resonance relations is given in Section
\ref{sec:resonance}.

\begin{thm*}
Let $x,\Lambda_2,\dots,\Lambda_n \in \C$ be generic. Then for every
$\psi^{(0)} \in M_{\vec\Lambda}[0]$ there exists a unique function
of the form \eqref{intro:n>1 BA type vector function}, satisfying
the $n$-point resonance relations.
\end{thm*}
 This is proved in Section \ref{sec:resonance}, see Theorem
\ref{thm:resonance n>1}.
\bigskip

Using the above theorem, we introduce the fundamental resonance
matrix $\Psi(\lambda,x;\vec\Lambda)$
\begin{equation*}
    \Psi(\lambda,x;\vec\Lambda) = q^{\lambda (x-\m)}
    \( \Id + \sum_{k=1}^\m  q^{2 k \lambda} \, \Psi^{(k)} \),
    \qquad
    \Psi^{(k)} \in \End(M_{\vec\Lambda}[0]),
\end{equation*}
whose colums form a basis of the space of functions satisfying the
resonance conditions. For  a fixed $x\in\C$ the fundamental
resonance matrix can be regarded as an operator
$\Psi(\lambda,x;\vec\Lambda): M_{\vec\Lambda}[0] \to \Fun \o
M_{\vec\Lambda}[0]$. Any function of the form \eqref{intro:n>1 BA
type vector function}, satisfying the resonance relations, must lie
in the image of $\Psi(\lambda,x;\vec\Lambda)$; in particular, this
applies to all hypergeometric qKZB solutions. From this we derive
that for generic $x$ the column spaces of
$\dH(\lambda,x;\vec\Lambda)$ and $\Psi(\lambda,x;\vec\Lambda)$
coincide.

\medskip
\noindent
1.3.\
Representation theory of quantum groups gives rise to another family
of solutions of the qKZB equations, realized as traces of
intertwining operators between Verma modules and their tensor
products, see \cite{EV}.

The versions of qKZB operators, considered in \cite{EV}, were
associated with the so-called fusion dynamical $R$-matrix, and the
corresponding qKZB solutions were meromorphic in $\lambda$. For the
simplest case of the quantum group $\U$ and $n=1$ these meromorphic
trace functions were explicitly computed in \cite{EV}, and were
related to the corresponding hypergeometric solutions by a simple
meromorphic gauge transformation.

In this paper we modify the trace function construction and obtain
holomorphic eigenfunctions of the qKZB operators \eqref{intro:qKZB
operator}, arising from Felder's dynamical $R$-matrix. We show that
these holomorphic trace functions are precisely the hypergeometric
qKZB solutions, computed using the contour integrals.

The key ingredient in our construction is the notion of a
holomorphic intertwining operator $\Phi_m^\Lambda(\lambda):
M_{\lambda-1} \to M_{\lambda-\Lambda+2m-1} \o M^\dual_{\Lambda}$,
where $M_\Lambda^\dual$ denotes the contragredient dual Verma
module, see \cite{STV}. These operators holomorphically depend on
highest weights $\lambda$ and $\Lambda$. For any $\vec\Lambda \in
\C^n$ and $\vec m \in \Z_{\ge0}^n$, satisfying the zero weight
condition $\Lambda_1+\ldots+\Lambda_n = 2(m_1+\ldots+m_n)$, we
consider the composition of the holomorphic intertwining operators
\begin{equation*}
    \Phi^{\vec\Lambda}_{\vec m}(\lambda): M_{\lambda-1} \to
    M_{\lambda-1} \o M^\dual_{\Lambda_1} \o\dots\o
    M^\dual_{\Lambda_n},
\end{equation*}
and the corresponding $(M_{\Lambda_1}^\dual \o\dots\o
M_{\Lambda_n}^\dual)[0]$-valued trace function:
\begin{equation*}
    \F_{\vec m}(\lambda,x;\vec \Lambda) = q^{-\frac{\m (\m+1)} 2} (q-q^{-1})^{\m} (q^x-q^{-x}) \ \,
    \Tr\biggr|_{M_{\lambda-1}} \!\!\! \(\Phi_{\vec m}^{\vec \Lambda} (\lambda) \ q^{x \mathbf h}
    \).
\end{equation*}
In the above formula, the trace is a formal power series in $q^x$,
which converges to a $M_{\vec\Lambda}^\dual[0]$-valued function
$\F_{\vec m}(\lambda,x;\vec\Lambda)$,  meromorphic in  $x$ and
holomorphic in variables $\lambda$ and $\vec\Lambda$. Arranging the
coordinates of various $\F_{\vec m}(\lambda,x;\vec\Lambda)$ as rows,
we obtain the universal trace matrix $\F(\lambda,x;\vec\Lambda)$.
For a fixed $x \in \C$, it can be thought of as an operator
$\F(\lambda,x;\vec\Lambda): M_{\vec\Lambda}[0] \to \Fun \o
M_{\vec\Lambda}[0]$.

The holomorphic intertwining operators $\Phi_m^\Lambda(\lambda)$ and
their compositions have particularly nice compatibility with
inclusions and fusion of Verma modules for $\U$. That compatibility
leads to the following result.

\begin{thm*}
Let $x \in \C$, and let $\psi(\lambda)$ be a
$M_{\vec\Lambda}[0]$-valued function of $\lambda$, which belongs to
the image of $\F(\lambda,x;\vec\Lambda)$ in $\Fun \o
M_{\vec\Lambda}[0]$. Then $\psi(\lambda)$ satisfies the $n$-point
resonance relations.
\end{thm*}
This is proved in Section \ref{sec:intertwining}, see Theorem
\ref{thm:resonance relations for traces}.
\bigskip

By uniqueness of the resonance solutions, the image of
$\F(\lambda,x;\vec\Lambda)$ is a subspace of the image of
$\Psi(\lambda,x;\vec\Lambda)$. In fact, these two spaces coincide,
and we have the following relation between
$\F(\lambda,x;\vec\Lambda)$ and $\Psi(\lambda,x;\vec\Lambda)$.

\begin{thm*}
There exists an invertible operator $\Xi_{\vec\Lambda} \in
\End(M_{\vec\Lambda})$, diagonal with respect to the standard
monomial basis, such that
\begin{equation*}
    \Psi(\lambda,x;\vec\Lambda) = \F(\lambda,x;\vec\Lambda) \ \Xi_{\vec\Lambda}.
\end{equation*}
\end{thm*}

This is proved in Section \ref{sec:intertwining}, see Theorem
\ref{thm:trace vs resonance}.
\bigskip

The diagonal operators $\Xi_{\vec\Lambda}$ play the role of the
operators, relating the actions of the quantum and the dynamical
quantum groups, associated with $\sl_2$. In particular, they provide
the gauge equivalence between Drinfeld's quantum $R$-matrix and the
asymptotic $R$-matrix, obtained as the limit of
$\dR_{\Lambda_1,\Lambda_2}(\lambda)$ as $\lambda \to \infty$, see
Theorem \ref{thm:gauge equivalence for R matrices}.

\medskip
\noindent
1.4.\
We establish the connection between the holomorphic intertwining
operators and the so-called dynamical Shapovalov form, which
appears, in particular, in the study of the qKZB heat equation and
the orthogonality relations, satisfied by the hypergeometric qKZB
solutions, see \cite{{FV heat-0}}, \cite{FV heat}, \cite{TV small}.
The dynamical Shapovalov form can be thought of as the bilinear
pairing $\mathbb Q: M_{\vec\Lambda}^* \o M_{\vec\Lambda}^* \to
\Fun$, defined by $\mathbb Q(u_{\vec m}^{\vec\Lambda} \o u_{\vec
m'}^{\vec\Lambda}) = \delta_{\vec m,\vec m'} \mathbb Q_{\vec
m}^{\vec\Lambda}(\lambda)$, where
\begin{equation}\label{intro:dynamical Shapovalov}
    \mathbb Q_{\vec m}^{\vec\Lambda}(\lambda) =
    \prod_{j=1}^n \, [m_i]!  \prod_{k=1}^{m_i}  \frac
    {[\lambda-\sum_{i=j+1}^n
(\Lambda_i-2m_i)+k]
[\lambda- \sum_{i=j}^n(\Lambda_i-2m_i) - k]}{[\Lambda_i-k+1]}.
\end{equation}

We interpret the dynamical Shapovalov form as the "square norm" of
the holomorphic intertwining operators with respect to the inner
product $\<\cdot,\cdot\>$ on Verma modules, see \cite{L}.

\begin{thm*}
Let $\Phi_{\vec m}^{\vec\Lambda}(\lambda): M_{\lambda-1} \to M_\mu
\o M_{\vec\Lambda}^\dual$ and $\Phi_{\vec
m'}^{\vec\Lambda}(\lambda'): M_{\lambda'-1} \to M_\mu \o
M_{\vec\Lambda}^\dual$ be the holomorphic intertwining operators as
above. Then for any $v \in M_{\lambda-1}$ and $v' \in
M_{\lambda'-1}$ we have
\begin{equation*}
    \<\Phi_{\vec m}^{\vec\Lambda}(\lambda) v, \Phi_{\vec
    m'}^{\vec\Lambda}(\lambda') v'\> = \delta_{\vec m,\vec m'}\
    \mathbb Q_{\vec m}^{\vec\Lambda}(\lambda) \  \<v,v'\>.
\end{equation*}
\end{thm*}
This is proved in Section \ref{sec:intertwining}, see Theorem
\ref{thm:orthogonality lemma}.

\medskip

The dynamical Shapovalov form appears as the density function for
the integrals, participating in the qKZB heat equation and the
orthogonality relations for the hypergeometric qKZB solutions. These
relations are crucial for the construction in \cite{FV heat} and
\cite{EV heat} of the generalized Fourier transform for the
vector-valued spherical functions. In \cite{FV heat} the qKZB heat
equation and the orthogonality relations were proved for $n$-point
hypergeometric solutions at level $\m=1$, and in \cite{EV heat} they
were established for the one-point meromorphic trace functions,
associated with arbitrary quantum groups. We hope that the previous
theorem will be instrumental for a uniform general proof of these
relations.

\medskip
\noindent
1.5.\
The resonance relations were introduced in \cite{FV} as the
identities satisfied by the elliptic hypergeometric qKZB solutions,
and our qKZB matrix inherits the same properties.

\begin{thm*}\cite{FV}
Let $x \in \C$, and let $\psi(\lambda)$ be a
$M_{\vec\Lambda}[0]$-valued function of $\lambda$, which belongs to
the image of $\dH(\lambda,x;\vec\Lambda)$ in $\Fun \o
M_{\vec\Lambda}[0]$. Then $\psi(\lambda)$ satisfies the $n$-point
resonance relations.
\end{thm*}

The $M_{\vec\Lambda}^\dual[0]$-valued functions $\dH^{\vec
m}(\lambda,x;\vec\Lambda)$, represented by columns of
$\dH(\lambda,x;\vec\Lambda)$, can therefore be expressed as linear
combinations of columns $\F^{\vec m}(\lambda,x;\vec\Lambda)$ of the
universal trace matrix. The corresponding transformation operator is
closely related to the quantum Knizhnik-Zamolodchikov (qKZ)
equations, which we now describe.

For $x\in\C$ and $j=1,\dots, n$ the qKZ operators $\qqKZ_j(x)$
without spectral parameter, acting in $M_{\vec\Lambda}[0]$, are
expressed in terms of the quantum $R$-matrix operators
$\qqR_{\Lambda_i,\Lambda_k} \in \End(M_{\Lambda_i} \o
M_{\Lambda_k})$:
\begin{equation*}
    \qqKZ_j(x) = \qqR_{\Lambda_j,\Lambda_{j+1}}^{-1} \dots
    \qqR_{\Lambda_j,\Lambda_n}^{-1}
    \ \(q^{x \bh}\)_{j} \ \qqR_{\Lambda_1,\Lambda_j} \dots
    \qqR_{\Lambda_{j-1},\Lambda_j}\ ,
\end{equation*}
The qKZ equations are equations for common eigenvectors of the
operators $\qqKZ_j(x)$.

The qKZ operators $\qqKZ_j(x)$ can be regarded as certain limits of
the qKZB operators $\qKZB_j$, restricted to the space of
$M_{\vec\Lambda}[0]$-valued quasi-trigonometric polynomials. The
corresponding limits of the qKZB equations become the qKZ equations.

For every level $\m$ a version of the universal hypergeometric
construction \cite{FTV monodromy} yields a matrix $\{\qqH^{\vec
m'}_{\vec m}(x;\vec\Lambda)\}$, holomorphically depending on $x$.
Regarded as an operator $\qqH(x;\vec\Lambda) \in
\End(M_{\vec\Lambda}[0])$, it satisfies
\begin{equation}\label{eq:intro qKZ equations}
    \qqKZ_j(x) \ \qqH(x;\vec\Lambda) =
    \qqH(x;\vec\Lambda) \ \mathcal E_j(x;\vec\Lambda),
    \qquad j=1,\dots,n,
\end{equation}
for the same operators $\mathcal E_j(x;\vec\Lambda)$ as in
\eqref{eq:intro qKZB equations}. We obtain the following fundamental
relations between the universal  trace functions and the qKZ and
qKZB operators and solutions.

\begin{thm*}
Let $n\ge2$ and let $x \in \C$ be generic. Then for all
$j=1,\dots,n$ we have
\begin{equation*}
    \qKZB_j \  \F(\lambda,x;\vec \Lambda) =
    \F(\lambda,x;\vec \Lambda) \  \qqKZ_j(x).
\end{equation*}
\end{thm*}
This is proved in Section \ref{sec:equations}, see Theorem
\ref{thm:holonomic system}.

\begin{thm*}
Let $x \in \C$ be generic. Then
\begin{equation*}
    \dH(\lambda,x;\vec\Lambda) = \F(\lambda,x;\vec \Lambda) \
    \qqH(x;\vec\Lambda).
\end{equation*}
\end{thm*}
This is proved in Section \ref{sec:hypergeometric}, see Theorem
\ref{thm:qKZ to qKZB}.

\medskip

In other words, the universal trace function intertwines the qKZB
and qKZ operators, and transforms the hypergeometric solutions of
the qKZ equations to the hypergeometric solutions of the qKZB
equations. Thus we can regard $\F(\lambda,x;\vec \Lambda)$ as the
"quantization" operator, reconstructing the qKZB solutions from
their asymptotic limit. Note also that the previous theorem gives a
formula for the universal trace function in terms of hypergeometric
integrals,
\begin{equation*}
    \F(\lambda,x;\vec \Lambda) = \dH(\lambda,x;\vec\Lambda)
    \qqH(x;\vec\Lambda)^{-1}.
\end{equation*}

\medskip
\noindent
1.6.\
The Macdonald-Ruijsenaars operators $\mathbb M_\Theta: \Fun \o
M_{\vec\Lambda}[0] \to \Fun \o M_{\vec\Lambda}[0]$, introduced in
\cite{EV}, are defined for dominant integral highest weights
$\Theta$ by
\begin{equation}
    \mathbb M_\Theta = q^{\m\Theta} \
    \sum_{\mu \in \h^*} \Tr\biggr|_{L_\Theta[\mu]}
    \dR_{\Theta,\Lambda_1}(\lambda-\bh^{(2,\dots,n)}) \dots
    \dR_{\Theta,\Lambda_n}(\lambda) \ \mathbb T_{-\mu},
\end{equation}
where $\mathbb T_{-\mu}$ is the shift operator, $\mathbb
T_{-\mu}\psi(\lambda) = \psi(\lambda-\mu)$, and $L_\Theta$ is the
irreducible finite-dimensional $\U$-module with highest weight
$\Theta$.

The operators $\mathbb M_\Theta$ for different $\Theta$ pairwise
commute. The Macdonald-Ruijsenaars equations are the equations for
common eigenfunctions of operators $\mathbb M_\Theta$.

We use the equivalence between the hypergeometric and trace function
qKZB solutions to establish the following result.

\begin{thm*}
For $\Theta\in \Z_{\ge0}$, let $\chi_\Theta(x) = \sum_{\mu} \dim
L_\Theta[\mu] \, q^{-\mu x}$ denote the character of $L_\Theta$.
Then
\begin{equation}
\label{eq:hypergeometric Ruijsenaars}
    \mathbb M_\Theta \ \dH(\lambda,x;\vec\Lambda) =
    \chi_\Theta(x) \, \dH(\lambda,x;\vec\Lambda).
\end{equation}
\end{thm*}
This is proved in Section \ref{sec:equations}, see Theorem
\ref{thm:hypergeometric Ruijsenaars}.

One can show that the hypergeometric construction gives all
trigonometric quasi-polynomial solutions to equations
(\ref{eq:hypergeometric Ruijsenaars}), cf. Theorem
\ref{thm:completeness of hypergeometric MR}.

\medskip

When $n=1$ the Macdonald-Ruijsenaars operators reduce to the
ordinary Macdonald operators for $\mathfrak{sl}_2$ with parameters
$q$ and $t=q^\m$. The Macdonald polynomials $p_\lambda(x)$ are then
related to the hypergeometric solution $\dH(\lambda,x;2\m)$ by
\begin{equation}\label{intro:original Macdonald and hypergeometric}
    p_\lambda(x) = \frac{(-1)^{\m} q^{\m(\m+1)/2}}{(q-q^{-1})^{2\m+1}} \, \frac{[2\m]!}{[\m]!}
    \ \frac {\dH(\lambda+\m+1,x;2\m) - \dH(\lambda+\m+1,-x;2\m)}
    {\prod_{j=1}^{\m} [\lambda+j] \, \prod_{j=-\m}^{\m} [x-j]}.
\end{equation}
Thus the anti-symmetrized qKZB solutions, studied in this paper, can
be regarded as vector-valued generalizations of the Macdonald
polynomials for $\mathfrak {sl}_2$.

\medskip
\noindent
1.7.\
Summarizing the results, we conclude that for generic $x \in \C$
there exists a distinguished subspace of $M_{\vec\Lambda}[0]$-valued
trigonometric quasi-polynomials, which can be constructed in several
equivalent ways: using the resonance relations, the hypergeometric
integrals, the holomorphic trace functions, the qKZB equations, or
the MR equations, see Theorem \ref{thm:five definitions generic}. We
call this subspace the harmonic space $\Harm_{x,\vec\Lambda}$.

For integral dominant $\vec\Lambda$ we obtain the integral version
$\adm\Harm_{x,\vec\Lambda} \subset \Fun \o L_{\vec\Lambda}[0]$ of
the harmonic space, which consists of functions with values in the
tensor product of finite-dimensional $\U$-modules. The space
$\adm\Harm_{x,\vec\Lambda}$ can be characterized as the space of
$L_{\vec\Lambda}[0]$-valued trigonometric quasi-polynomial solutions
of qKZB or MR equations.

\medskip
\noindent
1.8.\
When the highest weights $\vec\Lambda$ are nonnegative integers, the
Weyl group acts on the space of functions with values in the
corresponding tensor product of finite-dimensional modules. The Weyl
anti-symmetric solutions of KZB-type equations play an important
role in the conformal field theory, see for example \cite{FW}.

It was shown in \cite{FV} that the coordinates of Weyl
anti-symmetrized hypergeometric qKZB solutions $\psi(\lambda)$ must
vanish at the special values of $\lambda$, participating in the
resonance relations. These vanishing conditions can be described in
terms of the fusion rules for the tensor category of
finite-dimensional $\sl_2$-modules. Namely, for each $\lambda$ and
$\vec m$ one constructs a chain of highest weights
$\lambda^{(0)},\dots,\lambda^{(n)}$, and the vanishing conditions
state that $\psi_{\vec m}(\lambda)=0$ when at least one of the
triples $(\lambda^{(i-1)},\lambda^i,\Lambda_i)$ violates the
$\sl_2$-fusion rules. The proof in \cite{FV} of the vanishing
conditions was based on some identities for theta functions entering
the hypergeometric integrals.

We use the equivalence between the spaces of hypergeometric
solutions and holomorphic trace functions to give a new
representation-theoretic proof of the vanishing conditions. Namely,
we prove the Weyl formula for the holomorphic trace functions
(Theorem \ref{thm:Weyl formula}), which shows that the
anti-symmetrized trace functions can be interpreted as traces of
regularized intertwining operators between finite-dimensional
$\U$-modules. The violation of the fusion rules as above then
implies that the corresponding intertwining operator $\Phi_{\vec
m}^{\vec\Lambda}(\lambda)$ must be the zero operator, and hence the
corresponding trace function vanishes.

\medskip
\noindent 1.9.
The anti-symmetrized hypergeometric qKZB solutions belong to the sum
of spaces $\adm\Harm_{x,\vec\Lambda}+\adm\Harm_{-x,\vec\Lambda}$.
For generic $x$ this sum is direct, because the "supports" of the
corresponding trigonometric quasi-polynomials are the disjoint
strings $(x-\m,\dots,x+\m)$ and $(-x-\m,\dots,-x+\m)$. In
particular, the dimension of the space of Weyl anti-symmetric qKZB
solutions is the same for all generic $x$.

The hypergeometric qKZB solutions become the true trigonometric
polynomials when the parameter $x$ is an integer, and for small
integral values of $x$ the spaces $\adm\Harm_{x,\vec\Lambda}$,
$\adm\Harm_{-x,\vec\Lambda}$ have a nontrivial intersection. The
structure of the space of Weyl anti-symmetric qKZB solutions is more
subtle in these cases.

We show that for $x=0$ all Weyl anti-symmetrized hypergeometric qKZB
solutions vanish identically.

Then we consider the case $x = \pm 1$, and study in detail the
corresponding Weyl anti-symmetric solutions $\cb^{\vec m}(\lambda)$,
which we call the quantum conformal blocks. The following result
resembles the Macdonald special value identity, and plays an
important role in the analysis of functions $\cb^{\vec
m}(\lambda)$.

\begin{thm*}
Let $\vec\Lambda \in \ZZ^n[2\m]$ and $\vec m \in \Z^n[\m]$. Then
$\vartheta^{\vec m}(1) = - \mathbb Q_{\vec m}^{\vec\Lambda}(1) \,
v_{\vec\Lambda}^{(\vec m)}$, where $\mathbb Q_{\vec
m}^{\vec\Lambda}(\lambda)$ are the diagonal entries of the matrix of
the dynamical Shapovalov form, see \eqref{intro:dynamical
Shapovalov}.
\end{thm*}

This is proved in Section 8, see Theorem \ref{thm:theta(1)}.

\medskip
\noindent
1.10.\
We denote $\Conf_{\vec\Lambda}$ the subspace of $\Fun \o
L_{\vec\Lambda}[0]$, spanned by the quantum conformal blocks
$\cb^{\vec m}(\lambda)$. One of our principal results is

\begin{thm*}
For any $\vec\Lambda \in \ZZ^n$ the dimension of the space
$\Conf_{\vec\Lambda}$ equals the dimension of $\U$-invariants in the
tensor product $L_{\Lambda_1} \o \dots \o L_{\Lambda_n}$.
\end{thm*}

This is proved in Section \ref{sec:conformal blocks}, see Theorem
\ref{thm:conformal blocks generic}.

\bigskip

The proof of the above theorem is constructive in the sense that we
explicitly describe the set of indices $\vec m$, such that
$\cb^{\vec m}(\lambda)$ is a basis of $\Conf_{\vec\Lambda}$.
Moreover, we show that for all other $\vec m$ the functions
$\cb^{\vec m}(\lambda)$ are identically zero. This property suggests
that there may be a correspondence between hypergeometric solutions
and Lusztig's canonical basis, which have similar properties. We
plan to investigate this relation in a subsequent paper.

\medskip
\noindent
1.11.\
When the highest weights $\vec\Lambda$ are integral, the qKZB and MR
equations can be considered independently on each of the cosets
$\C/\Z$. Particularly important is the restriction to the lattice
$\Z \subset \C$. In that case the operators $\qKZB_j$ and $\mathbb
M_\Theta$ may have singularities, and we consider their modified
versions $\widetilde{\qKZB}_j$ and $\widetilde{\mathbb M}_\Theta$,
which are regular on the lattice. We conjecture that the
restrictions of the Weyl anti-symmetrized hypergeometric solutions
remain the eigenfunctions of these modified qKZB and
Macdonald-Ruijsenaars operators; in some special cases this
conjecture was shown to hold, see Section 12 in \cite{FV}.

For generic $q$ the restriction of a trigonometric polynomial
$\cb^{\vec m}(\lambda)$ can be uniquely reconstructed from its
restriction to $\lambda \in \Z$. However, when $q$ becomes a root of
unity, some of the $\cb^{\vec m}(\lambda)$ may restrict to the zero
function on $\Z$. For $q = \exp\(\frac{\pi i}\ell\)$ with $\ell \in
\{ 3,4,\dots \}$, we define the discrete version
$\Conf_{\vec\Lambda}^{(\ell)}$ of the conformal block space as the
subspace of functions $\psi:\Z \to L_{\vec\Lambda}[0]$, spanned by
the restrictions of $\cb^{\vec m}(\lambda)$.

\begin{thm*}
Let $\vec\Lambda \in \ZZ^n$ be such that $\Lambda_i \le \ell-2$ for
$i=1,\dots,n$. Then the dimension of the space
$\Conf_{\vec\Lambda}^{(\ell)}$ equals the dimension of the Verlinde
algebra invariants in the tensor product $L_{\Lambda_1} \o \dots \o
L_{\Lambda_n}$.
\end{thm*}

This is proved in Section \ref{sec:conformal blocks}, see Theorem
\ref{thm:conformal blocks roots}.

\bigskip

Finally, we conjecture that the spaces $\Conf_{\vec\Lambda}$ and
$\Conf_{\vec\Lambda}^{(\ell)}$ are uniquely characterized as the
spaces of Weil anti-symmetric trigonometric polynomial solutions of
the MR equations
\begin{equation*}
    \mathbb M_\Theta \psi(\lambda) = \(\dim_q L_\Theta\)  \psi(\lambda),
\end{equation*}
see Conjecture \ref{conj:conformal generic} and Conjecture
\ref{conj:conformal root} for precise formulations.

\medskip
\noindent 1.12.\ The hypergeometric construction for the general
$sl_2$ qKZB equations, associated with the elliptic $R$ matrices
with spectral parameters, was studied in \cite{FTV monodromy},
\cite{FTV Bethe}, \cite{FV}-\cite{FV heat}. The corresponding trace
function construction, based on the intertwining operators for the
affine quantum group $\Uhat$, was developed in \cite{ESV}. It was
proved in \cite{ESV} that the constructed trace functions satisfy
the general qKZB equations. Conjecturally the  elliptic
hypergeometric solutions and trace function solutions  coincide up
to normalization.

We expect that the techniques developed in this paper can be used to
establish the desired equivalence of the two constructions. We plan
to discuss this subject in a subsequent paper.

\medskip
\noindent
1.13.\
We are grateful to P. Etingof for valuable
discussions. The second author was supported by NSF grant
DMS-0555327.

\section{Resonance relations: the elementary approach}
\label{sec:resonance}

\subsection{Basic notation}

Fix $\eta \in \C$ such that $\Im \eta >0$, and for any $x \in \C, k
\in \Z_{\ge0}$, denote
\begin{equation*}
q^x = e^{\pi i \eta x }, \qquad [x] = \frac {q^x-q^{-x}}{q-q^{-1}},
\qquad [x]_k = [x][x+1]\dots[x+k-1], \qquad [k]! = [1][2]\dots[k] .
\end{equation*}
In this paper we fix a non-negative integer $\m$, and a positive
integer $n$. We set
\begin{equation*}
    \ZZ^n = (\Z_{\ge0})^n,
    \qquad
    \ZZ^n[\m] = \{\vec m \in\ZZ^n \ \bigr| \ m_1+\dots+m_n = \m \}.
\end{equation*}
We denote by $\C^n[2\m]$ the hyperplane in $\C^n$, defined by
\begin{equation*}
    \C^n[2\m] = \{ \vec\Lambda \in \C^n \,\bigr|\, \Lambda_1+\dots+\Lambda_n =
    2\m \}.
\end{equation*}
We say that $\vec\Lambda \in \C^n[2\m]$ is generic, if
$\Lambda_2,\dots,\Lambda_n$ are generic in $\C$.

We say that $m \in \Z_{\ge0}$ is $\Lambda$-admissible, if $m-\Lambda
\notin \Z_{>0}$. More generally, $\vec m \in \ZZ^n$ is called
$\vec\Lambda$-admissible, if $m_i-\Lambda_i \notin \Z_{>0}$ for all
$i=1,\dots,n$. For any $\vec\Lambda \in \ZZ^n$  we denote
\begin{equation}
    \Adm_{\vec\Lambda} = \left\{ \vec m \in \ZZ^n\ \bigr| \
    \vec m \text{ is $\vec\Lambda$-admissible} \right\},
    \qquad
    \Adm_{\vec\Lambda}[\m] = \Adm_{\vec\Lambda} \cap \ZZ^n[\m].
\end{equation}

Let $\Fun$ denote the space of functions of a complex variable
$\lambda$. For any $x \in \C$ and $d \in \Z_{\ge0}$ we denote
$\FunBA_{x,d}$ the subspace of $\Fun$, consisting of the functions
of the form
\begin{equation}\label{eq:def FunBA}
    \psi(\lambda) = q^{\lambda (x-\m)} \( \sum_{k=0}^d  \psi^{(k)} \, q^{2 k \lambda} \),
    \qquad
    \psi^{(k)} \in \C.
\end{equation}
We call $\psi(\lambda)$ of the form \eqref{eq:def FunBA} a
\textbf{\emph{trigonometric quasi-polynomial}} of degree $d$.

\subsection{The resonance relations}
Let $V_{n,\m}$ denote a complex vector space with a basis
$\{v^{(\vec m)}\}_{\vec m \in \ZZ^n[\m]}$. Any $V_{n,\m}$-valued
function $\varphi(\lambda)$ of a complex variable $\lambda$ can be
written as
\begin{equation}\label{eq:combinatorial generating function}
    \varphi(\lambda) = \sum_{\vec m \in \ZZ^n[\m]} \varphi_{\vec
    m}(\lambda) \, v^{(\vec m)},
\end{equation}
where $\varphi_{\vec m}(\lambda)$ are scalar functions; we call them
the coordinates of $\varphi(\lambda)$.

A holomorphic $V_{n,\m}$-valued function $\varphi(\lambda)$ is said
to satisfy the $n$-point resonance relations at level $\m$ with
respect to $\vec\Lambda \in \C^n[2\m]$, if the following conditions
hold for all $\vec m \in \ZZ^n[\m]$.

\bigskip

\begin{itemize}
\item {\bf Condition $\mathbf {C_j}(\vec m)$} for $j=1,\dots,n-1$. For
any $\delta \in \{1,\dots,m_j\}$ we have
\begin{equation}\label{eq:condition C_j abstract}
    \varphi_{\vec m} (-\delta + \sum_{i=j+1}^{n}(\Lambda_i-2m_i))
    = \varphi_{\tau_j^\delta(\vec m)} (-\delta +
    \sum_{i=j+1}^{n}(\Lambda_i-2m_i)),
\end{equation}
where $\tau_j^\delta(\vec m) =
(m_1,\dots,m_{j-1},m_j-\delta,m_{j+1}+\delta,m_{j+2},\dots,m_n)$.

\medskip

\item {\bf Condition $\mathbf {C_n}(\vec m)$}. For any $\delta \in
\{1,\dots,m_n\}$ we have
\begin{equation}\label{eq:condition C_n abstract}
    \varphi_{\vec m} (-\delta)
    = \varphi_{\tau_n^\delta(\vec m)} (\delta),
\end{equation}
where $\tau_n^\delta(\vec m) =
(m_1+\delta,m_2,m_3,\dots,m_{n-2},m_{n-1},m_n-\delta)$.
\end{itemize}
\bigskip

It is easy to see that for each $\vec m$ the values of
$\varphi_{\vec m}(\lambda)$ appear exactly $\m$ times on the left,
and $\m$ times on the right side of the equations
\eqref{eq:condition C_j abstract}, \eqref{eq:condition C_n
abstract}. Note also that these equations do not directly involve
$\Lambda_1$. Another useful observation is that if $\xi(\lambda)$ is
an even scalar function, then the resonance relations for
$\varphi(\lambda)$ imply the resonance relations for
$\xi(\lambda)\varphi(\lambda)$.

\medskip

\ex Let $n=1$. Then the resonance relations are
\begin{equation*}
    \varphi_\m(-\delta) = \varphi_\m(\delta),
    \qquad \qquad \delta = 1,\dots,\m.
\end{equation*}
\medskip

\ex Let $n=2$ and $\m = 1$. Then the resonance relations are
\begin{align*}
    \varphi_{1,0}(-1+\Lambda_2) &= \varphi_{0,1}(-1+\Lambda_2),\\
    \varphi_{0,1}(-1) &= \varphi_{1,0}(1).
\end{align*}
\qed

\medskip

\begin{thm}\label{thm:resonance n>1}
Let $x \in \C$ and $\vec\Lambda \in \C^n[2\m]$ be generic. Then for
any vector $\psi^{(0)} = \sum_{\vec m} \psi_{\vec m}^{(0)} \,
v^{(\vec m)}$ there exists a unique function $\psi(\lambda)$ of the
form
\begin{equation}\label{eq:n>1 BA type}
    \psi(\lambda) = q^{\lambda (x-\m)}
    \( \sum_{k=0}^\m  q^{2 k \lambda} \, \psi^{(k)} \),
    \qquad
    \psi^{(k)} \in V_{n,\m},
\end{equation}
satisfying the $n$-point resonance relations with respect to
$\vec\Lambda$.
\end{thm}
\begin{proof}
For any $\vec m$ we form an $\m$-tuple of complex numbers by
concatenating the lists of values of $\lambda$, appearing on the
left in the resonance relations, associated with $\vec m$:
\begin{equation*}
    (\z_{\vec m}^{(1)},\dots,\z_{\vec m}^{(\m)}) =
    (-1+\sum_{j=2}^n(\Lambda_j-2m_j),\dots,-m_1+\sum_{j=2}^n(\Lambda_j-2m_j),
    \dots \ \dots,
    -1,\dots,-m_n),
\end{equation*}
and set for each $k=0,\dots,\m$
\begin{equation}\label{eq:n>1 ansatz basis}
    \xi_{\vec m;k}(\lambda) = q^{k \lambda} \, \prod_{i=k+1}^\m [\lambda-\z_{\vec
    m}^{(i)}].
\end{equation}
The coordinates of $\psi(\lambda)$ can be uniquely represented in
the form
\begin{equation*}
    \psi_{\vec m}(\lambda) = q^{\lambda x}
    \sum_{k=0}^\m c_{\vec m}^{(k)} \, \xi_{\vec m;k}(\lambda),
    \qquad c_{\vec m}^{(k)}\in\C,
\end{equation*}
and the initial condition forces $c_{\vec m}^{(0)} = (q-q^{-1})^\m
\, \psi_{\vec m}^{(0)}$. The remaining coefficients $c_{\vec
m}^{(k)}$ with $k=1,\dots,\m$ are determined from the resonance
relations, which are equivalent to a system of linear equations on
$c_{\vec m}^{(k)}$. More precisely, the equations
\eqref{eq:condition C_j abstract} and \eqref{eq:condition C_n
abstract} become respectively
\begin{equation}\label{eq:resonance equation j}
    \sum_{k=0}^\m c_{\vec m}^{(k)} \, \xi_{\vec m;k}(-\delta + \sum_{i=j+1}^{n}(\Lambda_i-2m_i))
    =
    \sum_{k=0}^\m c_{\tau_j^\delta(\vec m)}^{(k)} \, \xi_{\tau_j^\delta(\vec m);k}(-\delta +
    \sum_{i=j+1}^{n}(\Lambda_i-2m_i)),
\end{equation}
\begin{equation}\label{eq:resonance equation n}
    \sum_{k=0}^\m c_{\vec m}^{(k)} \, \xi_{\vec m;k}(-\delta)
    =
    q^{2\delta x} \sum_{k=0}^\m c_{\vec m'}^{(k)} \,
    \xi_{\vec m';k}(\delta).
\end{equation}
Since the number of resonance relations equals the number of
variables $c_{\vec m}^{(k)}$, it suffices to show that the
determinant $D$ of the corresponding matrix is nonzero for generic
$x, \Lambda_2,\dots,\Lambda_n$. It is clear that $D$ is a polynomial
in $q^{2x}$; let $D_0$ denote its constant term.

Introduce a partial order $\preccurlyeq$ on $\ZZ^n$ by writing $\vec
m \preccurlyeq \vec m'$ if and only if
\begin{equation}\label{eq:def partial order}
    \sum_{i=1}^j m_i \ge \sum_{i=1}^j m'_i \qquad \text{ for all }
    j=1,\dots,n,
\end{equation}
and list the resonance relations according to the following rules:
\begin{enumerate}
    \item
    If $i<j$, then conditions $\mathbf C_{\mathbf i}$ appear
    before conditions $\mathbf C_{\mathbf j}$,
    \item
    If $\vec m \prec \vec m'$, then conditions $\mathbf C_{\mathbf j}(\vec m)$
    appear before conditions $\mathbf C_{\mathbf j}(\vec m')$,
    \item
    Conditions $\mathbf C_{\mathbf j}(\vec m)$ appear in increasing order of $\delta$.
\end{enumerate}
Simultaneously we make a compatible ordered list of the undetermined
coefficients $c_{\vec m}^{(k)}$, by associating to the equation
\eqref{eq:resonance equation j} the variable $c_{\vec
m}^{(m_1+\dots+m_{j-1}+\delta)}$, and similarly for
\eqref{eq:resonance equation n}. Consider the matrix of the system
of resonance relations with respect to these orderings of equations
and variables, and take its limit as $q^{2x}\to 0$. It immediately
follows from \eqref{eq:n>1 ansatz basis} that the resulting matrix
is upper-triangular, with diagonal entries $ q^{ k \, \z_{\vec
m}^{(k)}} \prod_{i=k+1}^n [\z_{\vec m}^{(k)} - \z_{\vec m}^{(i)}]$.

The constant term $D_0$ is equal to the the product of the above
entries, which are all nonzero for generic
$\Lambda_2,\dots,\Lambda_n$. Therefore, $D_0 \ne 0$, and hence $D$
does not vanish identically, which implies that generically we have
$D \ne 0$.
\end{proof}

\subsection{The determinant $D$}

Computations suggest that the determinant $D$, which appeared in the
proof of Theorem \ref{thm:resonance n>1}, is always a product of its
constant term $D_0$ and several factors of the form
$(1-q^{2(x+k)})$. We formulate

\begin{conj}
One has $D \ne 0$ unless $x, \Lambda_2,\dots,\Lambda_n$ satisfy one
of the equations
\begin{align}
    \label{eq:det singularity x}
    [x+k]&=0,
    \qquad k=-\m+2,\dots,\m,\\
    \label{eq:det singularity lambda}
    [\Lambda_{i+1} + \dots + \Lambda_j-k]&=0,
    \qquad k=2,\dots,2\m-2,
\end{align}
for some $1\le i < j \le n$.
\end{conj}
Equivalently, we expect that up to a constant and a non-vanishing
exponential function in $\vec\Lambda$, the determinant $D$ factors
into the product of terms appearing on the left in \eqref{eq:det
singularity x}, \eqref{eq:det singularity lambda} with certain
multiplicities, which can be described combinatorially. This
factorization property is reminiscent of some of the determinant
formulas arising in representation theory, although the
representation-theoretic significance of $D$ is not clear yet.

We illustrate the factorization property of the determinant $D$ -
and the argument of the proof of Theorem \ref{thm:resonance n>1} -
in another example.
\medskip

\ex Let $n=2$ and $\m=2$. Then the ordered list of resonance
relations for the coordinates of a $V_{n,\m}$-valued function
$\psi(\lambda)$ is given by
\begin{align*}
    \psi_{2,0}(-1+\Lambda_2) &= \psi_{1,1}(-1+\Lambda_2),\\
    \psi_{2,0}(-2+\Lambda_2) &= \psi_{0,2}(-2+\Lambda_2),\\
    \psi_{1,1}(-3+\Lambda_2) &= \psi_{0,2}(-3+\Lambda_2),\\
    \psi_{1,1}(-1) &= \psi_{2,0}(1),\\
    \psi_{0,2}(-1) &= \psi_{1,1}(1),\\
    \psi_{0,2}(-2) &= \psi_{2,0}(2).
\end{align*}
We look for a trigonometric quasi-polynomial function
$\psi(\lambda)$ of the form
\begin{align*}
    \begin{pmatrix}
    \psi_{0,2}(\lambda)\\
    \psi_{1,1}(\lambda)\\
    \psi_{2,0}(\lambda)
    \end{pmatrix}
    = q^{\lambda x} \
    \begin{pmatrix}
    c_{0,2}^{(0)} &\, [\lambda+1][\lambda+2] &+&
    c_{0,2}^{(1)} &\, q^\lambda [\lambda+2] &+& c_{0,2}^{(2)} &\, q^{2\lambda}\\
    c_{1,1}^{(0)} &\, [\lambda+1][\lambda-\Lambda_2+3] &+&
    c_{1,1}^{(1)} &\, q^\lambda [\lambda+1] &+& c_{1,1}^{(2)} &\, q^{2\lambda}\\
    c_{2,0}^{(0)} &\, [\lambda-\Lambda_2+1][\lambda-\Lambda_2+2] &+&
    c_{2,0}^{(1)} &\, q^\lambda [\lambda-\Lambda_2+1] &+& c_{2,0}^{(2)} &\, q^{2\lambda}
    \end{pmatrix},
\end{align*}
where $c_{0,2}^{(0)},c_{1,1}^{(0)},c_{2,0}^{(0)}\in\C$ are the
prescribed initial conditions, and $(c_{2,0}^{(1)}, c_{2,0}^{(2)},
c_{1,1}^{(1)}, c_{1,1}^{(2)}, c_{0,2}^{(1)}, c_{0,2}^{(2)})$ are
undetermined coefficients. The resonance relations are linear
equations on these variables, and the determinant $D$ of the
corresponding system is given by
\begin{equation*}
    D = \det \begin{pmatrix}
    q^{\Lambda_2-1} & q^{2\Lambda_2-2} & -q^{\Lambda_2-1} \, [\Lambda_2] & -q^{2\Lambda_2-2} & 0 & 0\\
    0 & q^{2\Lambda_2-4} & 0 & 0 & q^{\Lambda_2-2} \, [\Lambda_2] & -q^{2\Lambda_2-4}\\
    0 & 0 & q^{\Lambda_2-3
    } \, [\Lambda_2-2] & q^{2\Lambda_2-6} & q^{\Lambda_2-3} \, [\Lambda_2-2] & -q^{2\Lambda_2-6}\\
    q^{2x+1} \, [\Lambda_2-3] & -q^{2x+2} & 0 & q^{-2} & 0 & 0 \\
    0 & 0 & -q^{2x+1} [2] & -q^{2x+2} &  q^{-1} & q^{-2}\\
    q^{4x+2} \, [\Lambda_2-4] &  -q^{4x+4}  & 0 & 0 & 0 & q^{-4}
    \end{pmatrix}.
\end{equation*}
In the limit $q^{2x}\to 0$ the matrix above becomes
upper-triangular, and has determinant $D_0 =
q^{4\Lambda_2-15}[\Lambda_2-2]$. This implies that $D \ne 0$ for
generic $x$ and $\Lambda_2$. In fact, one checks that
\begin{equation*}
    D = q^{4\Lambda_2-15} [\Lambda_2-2] \,  (1-q^{2x}) (1-q^{2(x+1)})^2 (1-q^{2(x+2)}) .
\end{equation*}

\medskip

\begin{rem}
The existence part of Theorem \ref{thm:resonance n>1} fails when $x$
satisfies \eqref{eq:det singularity x} with $k = 1,\dots,\m$. For
$x$ satisfying the remaining conditions \eqref{eq:det singularity x}
or $\Lambda_2,\dots,\Lambda_n$ satisfying \eqref{eq:det singularity
lambda}, the existence part is valid, but the uniqueness part fails.
\end{rem}

\subsection{The fundamental resonance matrix $\Psi(\lambda,x; \vec \Lambda)$}

For every $\vec m' \in \ZZ^n[\m]$, the initial condition $\psi^{(0)}
= v^{(\vec m')}$ in Theorem \ref{thm:resonance n>1} produces a
unique function $\Psi^{\vec m'}(\lambda,x; \vec \Lambda)$ of the
form \eqref{eq:n>1 BA type}, satisfying the resonance relations and
meromorphically depending on parameters $x,\vec\Lambda$. Let
$\Psi(\lambda,x; \vec \Lambda) = \{\Psi_{\vec m}^{\vec
m'}(\lambda,x; \vec \Lambda)\}_{\vec m,\vec m' \in \ZZ^n[\m]}$
denote the square matrix, formed by arranging the coordinates of
$\Psi^{\vec m'}(\lambda,x; \vec \Lambda)$ as columns.

We call $\Psi(\lambda,x; \vec \Lambda)$ the
\textbf{\emph{fundamental resonance matrix}} at level $\m$. For a
fixed $x\in\C$ it can be thought of as an operator
\begin{equation}\label{eq:def universal resonance}
    \Psi(\lambda,x; \vec \Lambda): V_{n,\m} \to \Fun \o V_{n,\m},
    \qquad
    v^{(\vec m')} \mapsto
    \Psi^{\vec m'}(\lambda,x; \vec \Lambda).
\end{equation}

\medskip

\ex Let $n=1$. Then the fundamental resonance matrix reduces to a
scalar, and solving the resonance relations using the method of
Theorem \ref{thm:resonance n>1}, one gets the explicit formula
\begin{equation}\label{eq:example Psi n=1}
    \Psi(\lambda,x;2\m) = C_\m \  q^{\lambda x} \ \sum_{k=0}^\m (-q)^{-k} \,
    \frac{[\m+k]!}{[\m-k]![k]!}
    \frac{q^{k\lambda} [\lambda+k+1] \dots [\lambda+\m]}{(1-q^{-2(x+1)}) \dots
    (1-q^{-2(x+k)})},
\end{equation}
where $C_\m = (q^{-1}-q)^\m \ q^{\m(\m+1)/2}$ is a normalization
constant, cf. \cite{EV}.

\medskip

\ex Let $n=2$ and $\m = 1$. Then the fundamental resonance matrix is
given by
\begin{equation}\label{eq:example Psi n=2,m=1}
\begin{split}
    \Psi(\lambda,x; \vec\Lambda) &=
    \begin{pmatrix} \Psi_{0,1}^{0,1}(\lambda,x; \vec\Lambda) & \Psi_{0,1}^{1,0}(\lambda,x; \vec\Lambda)\\
    \Psi_{1,0}^{0,1}(\lambda,x; \vec\Lambda) & \Psi_{1,0}^{1,0}(\lambda,x; \vec\Lambda)
    \end{pmatrix}
    \\
    &=
    q^{\lambda (x-1)} \left\{
    \begin{pmatrix} 1 & 0 \\ 0 & 1 \end{pmatrix} + \frac{q^{2\lambda-\Lambda_2}}{[x+1]}
    \begin{pmatrix}
    q^2 \, [-x-\Lambda_2-1] & q^{x+1} \, [2-\Lambda_2] \\
    q^{-x+1} \, [\Lambda_2] & [-x-\Lambda_2+1]
    \end{pmatrix}
    \right\}.
\end{split}
\end{equation}

\subsection{Difference equations with respect to $x$}

The axiomatic characterization of the $\Psi$-function implies that
$\Psi(\lambda,x;\vec\Lambda)$, regarded as a matrix-valued function
of $x$, is an eigenfunction for a family of difference operators
with matrix coefficients, cf. \cite{CV}.

\begin{lem}\label{thm:CV descent lemma}
Let $\vec\Lambda\in \C^n[2\m]$ be generic. Let $\varphi(\lambda) \in
\FunBA_{y,d}$ for some $d \in \Z_{\ge0}$ and generic $y \in \C$.
Assume that $\varphi(\lambda)$ satisfies the $n$-point resonance
relations with respect to $\vec\Lambda$. Then $\varphi(\lambda)$ can
be represented as
\begin{equation*}
    \varphi(\lambda) = \sum_{j=0}^{d-\m} \Psi(\lambda,y + 2j;\vec\Lambda) \  w_j ,
    \qquad w_j \in V_{n,\m}.
\end{equation*}
In particular, if $d<\m$, then $\varphi(\lambda) \equiv 0$.
\end{lem}
\begin{proof}

We argue by induction on the degree $d$ of the trigonometric
quasi-polynomial $\varphi(\lambda)$. If $d<\m$, then
$\varphi(\lambda)$ has the form \eqref{eq:n>1 BA type} with $x =
y-2$ and initial condition $\psi^{(0)} = 0$, so the uniqueness part
of Theorem \ref{thm:resonance n>1} implies that $\varphi(\lambda)
\equiv 0$.

Assume now that $d \ge \m$, and consider the function
\begin{equation*}
    \phi(\lambda) = \varphi(\lambda) - \Psi(\lambda,y;\vec\Lambda) \,
    \varphi^{(0)} \in \FunBA_{y+2,d-1}.
\end{equation*}

It satisfies the resonance relations, and by the induction
hypothesis we have
\begin{equation*}
    \phi(\lambda) = \sum_{j=1}^{d-\m} \Psi(\lambda,y + 2j;\vec\Lambda) \ w_j,
    \qquad w_j \in V_{n,\m}.
\end{equation*}
Then $\varphi(\lambda) =  \Psi(\lambda,y;\vec\Lambda) \,
\varphi^{(0)} + \sum_{j=1}^{d-\m} \Psi(\lambda,y + 2j;\vec\Lambda) \
w_j,$ which proves the lemma.
\end{proof}

\medskip

Next we consider the special class of linear difference operators
$\mathfrak D$, acting in the space of $\End(V_{n,\m})$-valued
functions of a complex variable $x$, which have the form
\begin{equation}\label{eq:difference operator type}
    \mathfrak D = \sum_\mu B_\mu(x) \,  \mathbb T_\mu,
    \qquad B_\mu(x) \in \End(V_{n,\m}),
\end{equation}
where $\mu$ runs over a finite subset of $\Z$, $\mathbb T_\mu$ is
the difference operator acting by $\mathbb T_\mu \psi(\lambda) =
\psi(\lambda-\mu)$ and $B_\mu(x)$ act by {\it right} multiplication.
In other words, for any $\End(V_{n,\m})$ - valued function
$\varPsi(x)$ we have
\begin{equation}\label{eq:difference operator right action}
    \mathfrak D \, \varPsi(x) = \sum_\mu \varPsi(x+\mu) B_\mu(x).
\end{equation}

\begin{thm}
Let $\xi(\lambda) \in \C[q^\lambda,q^{-\lambda}]$ be such that
$\xi(\lambda) = \xi(-\lambda)$. Then there exists a unique
difference operator $\mathfrak D_\xi$ as in \eqref{eq:difference
operator type}, \eqref{eq:difference operator right action}, such
that for any $\lambda \in \C$ we have
\begin{equation}\label{eq:dual Ruijsenaars}
    \mathfrak D_\xi \, \Psi(\lambda,x;\vec\Lambda) =
    \xi(\lambda) \, \Psi(\lambda,x;\vec\Lambda).
\end{equation}
\end{thm}

\begin{proof}
For each $\vec m'$ the function $\xi(\lambda) \, \Psi^{\vec
m'}(\lambda,x;\vec\Lambda)$ satisfies the conditions of Lemma
\ref{thm:CV descent lemma}, and therefore we can write
\begin{equation*}
    \xi(\lambda) \, \Psi^{\vec m'}(\lambda,x;\vec\Lambda) =
    \sum_{\mu} \Psi(\lambda,x + \mu;\vec\Lambda) \  w_\mu^{\vec m'}(x) ,
    \qquad w_\mu^{\vec m'}(x,\xi) \in V_{n,\m},
\end{equation*}
where $\mu$ runs over some finite set of integers. For each $\mu$
let $B_\mu(x,\xi) \in \End(V_{n,\m})$ denote the operator, defined
by $B_\mu(x,\xi) \, v^{(\vec m')} = w_\mu^{\vec m'}(x,\xi)$. Then
the desired equation \eqref{eq:dual Ruijsenaars} clearly holds if we
set $\mathfrak D_\xi = \sum_\mu B_\mu(x,\xi) \, \mathbb T_\mu$. To
prove the uniqueness of $\mathfrak D_\xi(x)$, one must show that the
only difference operator, annihilating
$\Psi(\lambda,x;\vec\Lambda)$, is the zero operator. This is easily
done by considering the asymptotics $\lambda \to \infty$.
\end{proof}

\medskip

\ex Let $n = 2, \, \m = 1$ and $\xi(\lambda) = q^\lambda +
q^{-\lambda}$. Then \eqref{eq:dual Ruijsenaars} is illustrated by
the identity
\begin{equation*}
    \Psi(\lambda,x-1;\vec\Lambda) + \Psi(\lambda,x+1;\vec\Lambda)
    \begin{pmatrix}
    1 - \frac{q^{2-\Lambda_2}[\Lambda_2]}{[x][x+1]} & - \frac{q^{-\Lambda_2}[2-\Lambda_2]}{[x][x+1]}
    \\
    - \frac{q^{2-\Lambda_2}[\Lambda_2]}{[x][x+1]} & 1 - \frac{q^{-\Lambda_2}[2-\Lambda_2]}{[x][x+1]}
    \end{pmatrix}
    = (q^\lambda + q^{-\lambda}) \, \Psi(\lambda,x;\vec\Lambda)
\end{equation*}
for the matrix $\Psi(\lambda,x;\vec\Lambda)$ given by
\eqref{eq:example Psi n=2,m=1}.

\section{Holomorphic trace functions}
\label{sec:intertwining}

\subsection{Representations of the quantum group $\U$}

Let $\U$ be the algebra with unit 1, generated by $\mathbf E,\mathbf
F$ and $q^{x \mathbf h}$ for $x \in\C$, with relations
\begin{equation}
\label{eq:quantum commutation}
\begin{split}
    q^{x \mathbf h} \, \mathbf E = q^{2x} \  \mathbf E \, q^{x \mathbf h},
    &\qquad
    q^{x \mathbf h} \, \mathbf F = q^{-2x} \ \mathbf F \, q^{x \mathbf h},
    \\
    q^{x\mathbf h} q^{x\mathbf h} = q^{(x+y)\mathbf h},
    &\qquad
    \mathbf E \mathbf F  - \mathbf F \mathbf E =
    \frac{q^{\mathbf h} - q^{-\mathbf h}}{q-q^{-1}}.
\end{split}
\end{equation}

If $V$ is a $\U$-module, we say that $v \in V$ has weight
$\mu\in\C$, if $q^{x \mathbf h} \, v = q^{x \mu} v$. The set of all
such $v \in V$ forms the weight subspace $V[\mu]$. In this paper we
consider $\U$-modules which admit a weight subspace decompositon,
i.e. split into the direct sum of weight subspaces.

Let $\lambda \in \C$. Denote by $M_\lambda$ the Verma $\U$-module,
generated by a vector $v_\lambda$ of weight $\lambda$. It has a
basis
$\{v^{(0)}_{\lambda}=v_\lambda,v^{(1)}_{\lambda},v^{(2)}_{\lambda},\dots\}$,
such that
\begin{equation}\label{eq:Verma action}
    \mathbf E \, v^{(m)}_{\lambda} = [\lambda-m] \, v^{(m-1)}_{\lambda},
    \qquad
    q^{x \mathbf h} \, v^{(m)}_{\lambda} = q^{(\lambda-2m)x} \, v^{(m)}_{\lambda},
    \qquad
    \mathbf F \, v^{(m)}_{\lambda} = [m+1] \, v^{(m+1)}_{\lambda}.
\end{equation}
When $\lambda \in \Z_{>0}$, there exists a $\U$-inclusion
\begin{equation}\label{eq:Verma inclusion}
    \iota(\lambda): M_{-\lambda-1} \to M_{\lambda-1},
    \qquad
    v_{-\lambda-1} \mapsto \mathbf F^{\lambda} v_{\lambda-1}.
\end{equation}
The image of $\iota(\lambda)$ is then spanned by vectors
$\{v_{\lambda-1}^{(\lambda)}, v_{\lambda-1}^{(\lambda+1)},\dots\}$
and constitutes a proper submodule of $M_{\lambda-1}$. The
corresponding quotient has dimension $\lambda$ and is denoted by
$L_{\lambda-1}$.

Let $M_\lambda^\dual= \bigoplus_{m\in\Z_{\ge0}} \C u_m^{\lambda}$ be
the restricted dual of $M_\lambda$, where $\{u_m^{\lambda}\}$ is the
dual basis to $\{v^m_{\lambda}\}$, determined by $\<u_m,v^{m'}\> =
\delta_{m,m'}$. Define the contravariant $\U$-action on $M^\dual$ by
\begin{equation}\label{}
    \<\mathbf E u,v\> = \<u,\mathbf F \, v\>,
    \qquad
    \<\mathbf F u,v\> = \<u,\mathbf E \, v\>,
    \qquad
    \<q^{x\mathbf h} u,v\> = \<u,q^{x\mathbf h} \, v\>
\end{equation}
for any $u \in M^\dual$ and $v \in M$. Equivalently, we have
explicit formulae
\begin{equation}\label{eq:dual Verma action}
    \mathbf E \, u_m^{\lambda} = [m] \, u_{m-1}^{\lambda},
    \qquad
    q^{x \mathbf h} \, u_m^{\lambda} = q^{(\lambda-2m)x} \, u_m^{\lambda},
    \qquad
    \mathbf F \, u_m^{\lambda} = [\lambda-m] \, u_{m+1}^{\lambda}.
\end{equation}
We call $M_\lambda$ the contravariant dual Verma module with highest
weight $\lambda$. When $\lambda \in \Z_{\ge0}$, the vectors
$\{u_0^\lambda, \dots, u_\lambda^{\lambda}\}$ span a submodule
$L_\lambda^\dual$, isomorphic to $L_\lambda$.

An inner product on a $\U$-module $V$ is a bilinear form
$\<\cdot,\cdot\>: V \o V \to \C$, satisfying
\begin{equation*}
    \<X v,v'\> = \<v, \varrho(X) v'\>,
    \qquad\qquad
    X \in \U, \quad v,v' \in V,
\end{equation*}
where $\varrho$ is the algebra anti-involution of $\U$, determined
by
\begin{equation*}
    \varrho(\mathbf F) = q^{-\bh} \mathbf E, \qquad
    \varrho(\mathbf E) = \mathbf F q^{\bh}, \qquad
    \varrho(q^{x\bh})=  q^{x\bh}.
\end{equation*}
For any $\Lambda \in \C$ the Verma module $M_\Lambda$ admits a
unique up to proportionality inner product. We normalize it by the
condition $\<v_\Lambda,v_\Lambda\> = 1$; then it is given by the
explicit the formula
\begin{equation*}
    \<v_\Lambda^{(m)},v_\Lambda^{(m')} \> = \delta_{m,m'}   \, q^{-m(\Lambda-m+1)}\
    \frac{ \prod_{k=1}^m [\Lambda-k+1]}{[m]!}.
\end{equation*}

A tensor product of two $\U$-modules is equipped with a $\U$-action
by means of the comultiplication $\Delta$, determined by
\begin{equation}
\label{eq:def quantum comultiplication}
    \Delta(\mathbf E) = \mathbf E \o 1 + q^{\mathbf h} \o \mathbf E,
    \qquad
    \Delta(\mathbf F) =  \mathbf F \o q^{-\mathbf h} + 1 \o \mathbf F,
    \qquad
    \Delta (q^{x \mathbf h}) = q^{x \mathbf h} \o q^{x \mathbf h}.
\end{equation}

One immediately verifies that if $\U$-modules $V,W$ are equipped
with inner products, then their product gives an inner product form
on the $\U$-module $V \o W$.

\medskip
The comultiplication $\Delta$ is coassociative, and the tensor
product of any finite collection of $\U$-modules also acquires a
$\U$-module structure. For $\vec\Lambda =
(\Lambda_1,\dots,\Lambda_n) \in \C^n$ we set
\begin{equation*}
    M_{\vec\Lambda} = M_{\Lambda_1} \o\dots\o M_{\Lambda_n},
    \qquad\qquad
    M_{\vec\Lambda}^\dual = M_{\Lambda_1}^\dual \o\dots\o
    M_{\Lambda_n}^\dual.
\end{equation*}
The modules $M_{\vec\Lambda}$ and $M_{\vec\Lambda}^\dual$ have
standard monomial bases $\{v^{(\vec m)}_{\vec\Lambda}\}$ and
$\{u_{\vec m}^{\vec\Lambda}\}$ respectively, where $\vec m =
(m_1,\dots,m_n) \in \ZZ^n$ and we denote
\begin{equation*}
    v^{(\vec m)}_{\vec\Lambda} = v_{\Lambda_1}^{(m_1)} \o\dots\o v_{\Lambda_n}^{(m_n)}.
    \qquad\qquad
    u_{\vec m}^{\vec\Lambda} = u^{\Lambda_1}_{m_1} \o\dots\o u^{\Lambda_n}_{m_n}.
\end{equation*}

\medskip

Clearly, the vectors $\{u^{\vec\Lambda}_{\vec m}\}$ with
$\vec\Lambda$-admissible indices $\vec m$ form a basis of the
irreducible submodule $L_{\vec\Lambda}^\dual$, and similarly the
images of vectors $\{v_{\vec \Lambda}^{(\vec m)}\}$ form a basis of
the irreducible quotient $L_{\vec\Lambda}$.

\subsection{The holomorphic intertwining operators and the orthogonality lemma}
Let  $\lambda,\Lambda \in \C$ and $m \in \Z_{\ge0}$. Denote by
$\Phi_m^\Lambda(\lambda):M_{\lambda-1} \to M_{\lambda-\Lambda+2m-1}
\o M_\Lambda^\dual$ the intertwining operators determined by
\begin{equation}\label{eq:n=1 operator}
\begin{split}
    \Phi_m^\Lambda(\lambda) v_{\lambda-1}
    &= [-\lambda+\Lambda-2m+1]_m \sum_{i=0}^m \frac{q^{i(\lambda-\Lambda+2m-i)}}{[-\lambda-\Lambda+2m+1]_i} \,
    v_{\lambda-\Lambda+2m-1}^{(i)} \o \mathbf E^i u_m^\Lambda = \\
    &= \sum_{i=0}^m q^{i(\lambda-\Lambda+2m-i)} [-\lambda+\Lambda-2m+i+1]_{m-i} \
    v_{\lambda-\Lambda+2m-1}^{(i)} \o \mathbf E^i u_m^\Lambda.
\end{split}
\end{equation}
Matrix elements of the operator $\Phi_m^\Lambda(\lambda)$
holomorphically depend on highest weights $\lambda,\Lambda$. For any
$\vec m = (m_1,\dots,m_n) \in \ZZ^n$ we define the intertwining
operator
\begin{equation*}
    \Phi_{\vec m}^{\vec\Lambda} (\lambda) : M_{\lambda-1} \to
    M_{\lambda-1-\sum_{i=1}^n(\Lambda_i-2m_i)} \o M_{\Lambda_1}^\dual
    \o \dots \o M_{\Lambda_n}^\dual
\end{equation*}
as the composition of operators \eqref{eq:n=1 operator}:
\begin{equation*}
    \Phi_{\vec m}^{\vec\Lambda} (\lambda) =
    ( \Phi_{m_1}^{\Lambda_1} (\lambda-\sum_{i=2}^n(\Lambda_i-2m_i)) \o 1^{n-1} ) \dots
    \(\Phi_{m_{n-1}}^{\Lambda_{n-1}}(\lambda-\Lambda_n+2m_n) \o 1\)
    \Phi_{m_n}^{\Lambda_n}(\lambda).
\end{equation*}

\bigskip

For $m \in \Z_{\ge0}$ and generic $\Lambda \in \C$, we denote
\begin{equation*}
    \mathcal Q_m^\Lambda = q^{-m(\Lambda-m+1)} \, \frac{ [m]!}{ \prod_{k=1}^m [\Lambda-k+1]},
    \qquad\qquad
    \mathbb Q_m^\Lambda(\lambda) = [m]! \,
    \prod_{k=1}^{m} \frac
    {[\lambda+k][\lambda- k - \Lambda +2m]}{[\Lambda-k+1]}.
\end{equation*}
It is easy to see that the formula $\<u_m^\Lambda,u_{m'}^\Lambda\> =
\delta_{m,m'} \, \mathcal Q_m^\Lambda$ defines an inner product on
$M_\Lambda^\dual$. The bilinear form $\mathbb Q: M_\Lambda^\dual \o
M_\Lambda^\dual \to \Fun$, defined by $\mathbb
Q(u_m^\Lambda,u_{m'}^\Lambda) = \delta_{m,m'} \, \mathbb
Q_m^\Lambda(\lambda)$, is called the dynamical Shapovalov pairing
\cite{FV heat}. If $\Lambda \in \Z_{\ge0}$, then $M_\Lambda^\dual$
does not admit nonzero inner product, but the restriction of the
above bilinear forms to the irreducible submodule $L_\Lambda^\dual$
is well-defined, and determines an inner product and dynamical
Shapovalov pairing on $L_\Lambda^\dual$.

More generally, for $\vec m \in \ZZ^n$ and generic $\vec\Lambda \in
\C^n$ we set
\begin{equation}\label{eq:def Q and dynamical Q}
    \mathcal Q_{\vec m}^{\vec\Lambda} = \prod_{j=1}^n \mathcal Q_{m_j}^{\Lambda_j},
    \qquad\qquad
    \mathbb Q_{\vec m}^{\vec\Lambda}(\lambda) = \prod_{j=1}^n \mathbb Q_{m_j}^{\Lambda_j}(\lambda -
    \bh^{(j+1,\dots,n)}),
\end{equation}
which yields an inner product and the dynamical Shapovalov pairing
on $M_{\vec\Lambda}^\dual$. The following orthogonality lemma
establishes a relation between the dynamical Shapovalov pairing and
holomorphic intertwining operators $\Phi_{\vec
m}^{\vec\Lambda}(\lambda)$. We will use it in Section
\ref{sec:conformal blocks}, see Theorem \ref{thm:conformal blocks
generic}.

\begin{thm}\label{thm:orthogonality lemma}
Let $\mu \in \C, \, \vec m, \vec m' \in \ZZ^n$, and let $\vec\Lambda
\in \C^n$ be generic. Set $\lambda = \mu +
\sum_{j=1}^n(\Lambda_i-2m_i)$ and $\lambda' = \mu +
\sum_{j=1}^n(\Lambda_i-2m'_i)$. Then for any $v \in M_{\lambda-1}$
and $v' \in M_{\lambda'-1}$ we have
\begin{equation}\label{eq:orthogonality lemma}
    \bigr< \Phi_{\vec m}^{\vec\Lambda}(\lambda) \, v,
    \Phi_{\vec m'}^{\vec\Lambda}(\lambda') \, v' \bigr>
    =\delta_{\vec m,\vec m'} \ \mathbb Q_{\vec m}^{\vec\Lambda}(\lambda) \ \<v ,
    v'\>.
\end{equation}
\end{thm}

\begin{proof}
Both sides of \eqref{eq:orthogonality lemma} are trigonometric
polynomials in $q^\mu,q^{\vec\Lambda}$ and therefore, it suffices
verify the statement when $\mu$ is generic and
$\Lambda_1,\dots,\Lambda_n$ are large enough positive integers.
Under these assumptions there exists a $\U$-module isomorphism
\begin{equation*}
    (M_{\mu-1} \o L_{\Lambda_1}^\dual) \o L_{\Lambda_2}^\dual \o \dots
    L_{\Lambda_n}^\dual \cong
    \(\bigoplus_{k=0}^{\Lambda_1}  M_{\mu-\Lambda_1+2k-1} \) \o L_{\Lambda_2}^\dual \o \dots
    L_{\Lambda_n}^\dual,
\end{equation*}
and the direct summands are orthogonal with respect to the
contravariant form. Therefore, for $m_1 \ne m'_1$ the subspaces
$\Phi_{\vec m}^{\vec\Lambda}(\lambda)(M_{\lambda-1})$ and
$\Phi_{\vec m'}^{\vec\Lambda}(\lambda')(M_{\lambda'-1})$ are
orthogonal in $M_{\mu-1} \o L_{\vec\Lambda}^\dual$. If $m_1 = m'_1$,
then a similar argument shows that the above subspaces are
orthogonal for $m_2 \ne m'_2$, etc. Therefore, we obtain the
orthogonality relations \eqref{eq:orthogonality lemma} for $\vec m
\ne \vec m'$.

The desired statement for $\vec m = \vec m'$ is easily reduced by
induction to the case $n=1$, i.e.
\begin{equation}\label{eq:orthogonality induction base}
    \bigr< \Phi_{m}^{\Lambda}(\lambda) \, v,
    \Phi_{m}^{\Lambda}(\lambda) \, v' \bigr>
    = \mathbb Q_{m}^\Lambda(\lambda) \, \<v,v'\>.
\end{equation}
The left side of the above equation represents an inner product on
$M_{\lambda-1}$, which by uniqueness must be proportional to the
standard one. Therefore, it remains to show that the coefficient of
proportionality equals $\mathbb Q^{m}_\Lambda(\lambda)$. Let $v = v'
= v_{\lambda-1}$, and as before $\mu = \lambda-\Lambda+2m$. Then
\begin{multline*}
    \bigr<  \Phi_{m}^{\Lambda}(\lambda) \, v_{\lambda-1},
    \Phi_{m}^{\Lambda}(\lambda) \, v_{\lambda-1} \bigr>
    = \sum_{i=0}^m \(q^{i(\mu-i)} [-\mu+i+1]_{m-i}\)^2 \
    \<v_{\mu-1}^{(i)},v_{\mu-1}^{(i)}\>  \,
    \<\mathbf E^i u_m^\Lambda,\mathbf E^i u_m^\Lambda\>
    \\
    =
    \sum_{i=0}^m q^{i(\mu-i)-(m-i)(\Lambda-m+i+1)} \(\frac{[m]![-\mu+i+1]_{m-i}}{[m-i]!}\)^2 \
    \frac{[\mu-1]\dots[\mu-i]}{[i]!}
    \frac {[m-i]!}{[\Lambda]\dots[\Lambda-m+i+1]}
    \\
    = q^{-m(\Lambda-m+1)} \frac{[m]! ([\mu-1]\dots[\mu-m])^2}{[\Lambda] \dots [\Lambda-m+1]}
    \ \( \sum_{i=0}^m q^{i(\mu+\Lambda-2m+1)} \frac{[-m]_i  [\Lambda-m+1]_i }{ [i]![-\mu+1]_i} \).
\end{multline*}
Using the Gauss identity for the $q$-hypergeometric function $\hyp
\( \begin{matrix} a, b \\ c
\end{matrix}; z \)$, we see that
\begin{equation*}
    \sum_{i=0}^m q^{i(\mu+\Lambda-2m+1)} \frac{[-m]_i  [\Lambda-m+1]_i }{ [i]![-\mu+1]_i }
    =
    \hyp \( \begin{matrix} -m, \Lambda-m+1 \\ -\mu+1 \end{matrix}; 1
    \)
    =
    q^{m(\Lambda-m+1)}\, \frac{[-\mu-\Lambda+m]_m}{[-\mu+1]_m},
\end{equation*}
and therefore we obtain
\begin{equation*}
    \bigr< \Phi_{m}^{\Lambda}(\lambda) \, v_{\lambda-1},
    \Phi_{m}^{\Lambda}(\lambda) \, v_{\lambda-1} \bigr> =
    [m]! \frac{ [\mu-1]\dots[\mu-m] [\mu+\Lambda-m] \dots[\mu+\Lambda-2m+1]}
    {[\Lambda]\dots[\Lambda-m+1]} =    \mathbb
    Q_{m}^\Lambda(\lambda).
\end{equation*}
This establishes \eqref{eq:orthogonality induction base} and
concludes the proof of the theorem.
\end{proof}

\subsection{Two inclusion lemmas}

The operators $\Phi_m^\Lambda(\lambda)$ have nice compatibility
properties with respect to the inclusions of Verma modules.

\begin{lem}\label{thm:diagram1}
Let $m, \mu \in \Z_{\ge0}$ and $\Lambda \in\C$. Then we have a
commutative diagram
\begin{equation*}
\begin{gathered}\xymatrix{
    &
    M_{-\mu+\Lambda-2m-1}
    \ar[rd]^(0.6){\quad\Phi_{m+\mu}^\Lambda(-\mu+\Lambda-2m)}
    \ar[ld]_(0.6){\Phi_m^\Lambda(-\mu+\Lambda-2m)\quad}
    &
    \\
    M_{-\mu-1} \o M_{\Lambda}^\dual
    \ar[rr]_{\iota(\mu)\o 1}
    &&
    M_{\mu-1} \o M_{\Lambda}^\dual
}\end{gathered}.
\end{equation*}
\end{lem}

\begin{proof}
Denote for convenience $\lambda = -\mu+\Lambda-2m$. We compute
\begin{equation*}
\begin{split}
    \Phi_{m+\mu}^\Lambda(\lambda) v_{\lambda-1}
    &= \sum_{i=0}^{m+\mu} \frac{ q^{i(\mu-i)} [-\mu+i+1]_{m+\mu-i}} {[i]!} \,
    \mathbf  F^i v_{\mu-1} \o \mathbf E^i u_{m+\mu}^\Lambda
    \\
    &=
    \sum_{i=\mu}^{m+\mu}  \frac{q^{i(\mu-i)} [-\mu+i+1]_{m+\mu-i}}{[i]!} \,
    \mathbf F^i v_{\mu-1} \o \mathbf E^i u_{m+\mu}^\Lambda
    \\
    &=
    \sum_{j=0}^m \frac{q^{(j+\mu)(-j)}}{[j+\mu]!} \frac{[m]!}{[j]!} \,
    \mathbf F^{j+\mu} v_{\mu-1} \o \mathbf E^{j+\mu} u_{m+\mu}^\Lambda
    \\
    &=
    \sum_{j=0}^m \frac{q^{j(-\mu-j)}}{[j]!} \, \frac{[m+\mu]!} {[j+\mu]!} \,
    \mathbf F^j \iota(\mu) v_{-\mu-1} \o \mathbf E^{j} u_m^\Lambda
    \\
    &=
    \sum_{j=0}^m \frac{q^{j(-\mu-j)}}{[j]!} \, [\mu+j+1]_{m-j} \
    \iota(\mu) \mathbf F^j v_{-\mu-1} \o \mathbf E^{j} u_m^\Lambda
    \\
    &= (\iota(\mu) \o 1) \ \Phi_m^\Lambda(\lambda) v_{\lambda-1}.
\end{split}
\end{equation*}
This implies that the intertwining operators
$\Phi_m^\Lambda(\lambda)$ and $(\iota(\mu) \o 1) \
\Phi_{m+\mu}^\Lambda(\lambda)$ agree on the generator
$v_{\lambda-1}$ of the module $M_{\lambda-1}$, and therefore
coincide.
\end{proof}

\begin{lem}\label{thm:diagram2}
Let $m, \lambda \in \Z_{\ge0}$. Then we have a commutative diagram
\begin{equation*}
\begin{gathered}\xymatrix{
    M_{-\lambda-1}
    \ar[rd]_(0.4){\Phi_{\lambda+m}^\Lambda(-\lambda)\ }
    \ar[rr]^{\iota(\lambda)}
    &&
    M_{\lambda-1}
    \ar[ld]^(0.4){\Phi_m^\Lambda(\lambda)}
    \\
    &
    M_{\lambda - \Lambda+2m-1} \o M_{\Lambda}^\dual
    &
}\end{gathered}.
\end{equation*}
\end{lem}

\begin{proof}
Denote for convenience $\mu = \lambda - \Lambda+2m$. We have
\begin{equation}
\label{eq:stuff1}
\begin{split}
\Phi_{m}^\Lambda&(\lambda) \mathbf F^{\lambda} v_{\lambda-1} =
[-\mu+1]_m \sum_{i=0}^{m} \frac{q^{i(\mu-i)} }{[i]! \, [-\mu+1]_i}
\( \sum_{j=0}^{\lambda} q^{j(\lambda-j)}\qbinom{\lambda}{j} \mathbf
F^{i+j} v_{\mu-1} \o \mathbf F^{\lambda-j} q^{-j \mathbf h} \mathbf
E^i u_m^\Lambda \)
\\ & =
[-\mu+1]_m \sum_{i=0}^{m} \sum_{j=0}^{\lambda}
 \frac{q^{(i+j)(\mu-i-j)} }{[i]! \, [-\mu+1]_i } \, \frac{[\lambda-j+1]_j}{[j]!} \,
\frac{[m]!}{[m-i]!} \, \mathbf F^{i+j} v_{\mu-1} \o \mathbf
F^{\lambda-j} u_{m-i}^\Lambda
\\ & =
[-\mu+1]_m \sum_{k=0}^{m+\lambda}  q^{k(\mu-k)} \frac{
[\lambda-k+1]_k}{[k]!} \( \sum_i \frac{[-k]_i \, [-m]_i \,
[\Lambda-m+1]_i}{[i]! \, [-\mu+1]_i \, [\lambda-k+1]_i} \) \mathbf
F^k v_{\mu-1} \o \mathbf F^{\lambda-k} u_m^\Lambda
\\ & =
[-\mu+1]_m \sum_{k=0}^{m+\lambda}\frac{ q^{k(\mu-k)} }{[k]!} \,
[\lambda-k+1]_k\ \hyper \(\begin{matrix} -k, -m, \Lambda-m+1 \\
-\mu+1,\lambda-k+1 \end{matrix} \, ; 1\) \mathbf F^k v_{\mu-1} \o
\mathbf F^{\lambda-k} u_m^\Lambda ,
\end{split}
\end{equation}
where $\hyper$ denotes the generalized $q$-hypergeometric function.
The $q$-version of Saalsch\"utz's theorem  (see e.g. \cite{Sl})
states that for any nonnegative integer $k$ we have
\begin{equation*}
\hyper \( \begin{matrix} -k, \, a, \, b \\ c,\, 1+a+b-c-k
\end{matrix}\, ;\ 1\) = \frac{[c-a]_k \, [c-b]_k}{[c]_k \,
[c-a-b]_k},
\end{equation*}
and therefore
\begin{equation}
\label{eq:stuff2}
\begin{split}
 [\lambda-k+1]_k\ & \hyper \(\begin{matrix} -k, -m, \Lambda-m+1 \\ -\mu+1,\lambda-k+1 \end{matrix} \, ; 1\) =
 [\lambda-k+1]_k \, \frac{[-\mu+m+1]_k \, [-\mu -\Lambda + m]_k}{[-\mu+1]_k \, [-\mu - \Lambda + 2m]_k} =
\\ & =
(-1)^k [-\lambda]_k \, \frac{[-\mu+m+1]_k \, [-\lambda-m]_k}
{[-\mu+1]_k \, [-\lambda]_k} = \frac{[-\mu+m+1]_k}{[-\mu+1]_k} \,
[\lambda+m-k+1]_{k}.
\end{split}
\end{equation}
Moreover, we have
\begin{equation}
\label{eq:stuff3}
\begin{split}
[\lambda+m-k+1]_{k} \  & \mathbf F^{\lambda-k} \, u_m^\Lambda =
[\lambda+m-k+1]_{k} \, [\Lambda-m-\lambda+k+1]_{\lambda-k} \,
u_{m+\lambda-k}^\Lambda =
\\ & = [\Lambda-m-\lambda+k+1]_{\lambda-k}\  \mathbf E^k
u_{m+\lambda}^\Lambda = [-\mu+m+k+1]_{\lambda-k} \ \mathbf E^k
u_{m+\lambda}^\Lambda.
\end{split}
\end{equation}
Finally, combining \eqref{eq:stuff1}, \eqref{eq:stuff2} and
\eqref{eq:stuff3}, we compute
\begin{equation*}
\begin{split}
    \Phi_m^\Lambda(\lambda) \iota(\lambda) v_{-\lambda-1} & =
    [-\mu+1]_m \sum_{k=0}^{m+\lambda}\frac{ q^{k(\mu-k)} }{[k]!} \,
    \frac{[-\mu+m+1]_\lambda}{[-\mu+1]_k} \, \mathbf F^k v_{\mu-1} \o
    \mathbf E^k u_{m+\lambda}^\Lambda =
    \\ & =
    \sum_{k=0}^{m+\lambda}\frac{ q^{k(\mu-k)} }{[k]!} \,
    \frac{[-\mu+1]_{m+\lambda}}{[-\mu+1]_k} \ \mathbf F^k v_{\mu-1} \o
    \mathbf E^k u_{m+\lambda}^\Lambda =
    \Phi_{m+\lambda}^\Lambda(-\lambda) v_{-\lambda-1},
\end{split}
\end{equation*}
which shows that the intertwining operators
$\Phi_{m+\lambda}^\Lambda(-\lambda)$ and
$\Phi_m^\Lambda(\lambda)\iota(\lambda)$ agree on the generator
$v_{-\lambda-1}$ of the Verma module $M_{-\lambda-1}$, and therefore
coincide.
\end{proof}

\subsection{The universal trace matrix $\F(\lambda,x;\vec \Lambda)$ }

Let $\vec\Lambda \in \C^n[\m]$ for some $\m \in \Z_{\ge0}$. For any
$\vec m \in \ZZ^n[\m]$ we define the
$M_{\vec\Lambda}^\dual[0]$-valued holomorphic trace function
$\F_{\vec m}(\lambda,x;\vec \Lambda)$ by
\begin{equation}\label{eq:definition trace function}
    \F_{\vec m}(\lambda,x;\vec \Lambda) = q^{-\frac{\m (\m+1)} 2} (q-q^{-1})^{\m} (q^x-q^{-x}) \ \,
    \Tr\biggr|_{M_{\lambda-1}} \!\!\! \(\Phi_{\vec m}^{\vec \Lambda} (\lambda) \ q^{x \mathbf h}
    \).
\end{equation}
In the above formula, the trace is an infinite power series in
$q^{-2x}$, which converges in the domain $|q^{2x}| \gg 1$ to a
meromorphic function $\F_{\vec m}(\lambda,x;\vec \Lambda)$.

Using the standard basis $\{u^{\vec\Lambda}_{\vec m'}\}$ of
$M_{\vec\Lambda}^\dual[0]$, we define the matrix elements
$\{\F_{\vec m}^{\vec m'}(\lambda,x;\vec \Lambda)\}$ by
\begin{equation*}
    \F_{\vec m}(\lambda,x;\vec \Lambda) = \sum_{\vec m' \in
    \ZZ^n[\m]} \F_{\vec m}^{\vec m'}(\lambda,x;\vec \Lambda) \,
    u^{\vec\Lambda}_{\vec m'},
\end{equation*}
Let $\F(\lambda,x;\vec \Lambda) = \{\F_{\vec m}^{\vec
m'}(\lambda,x;\vec \Lambda)\}_{\vec m,\vec m' \in \ZZ^n[\m]}$ denote
the square matrix, formed by arranging coordinates of $\F_{\vec
m}(\lambda,x;\vec \Lambda)$ as rows. The columns of this matrix
represent vectors
\begin{equation*}
    \F^{\vec m'}(\lambda,x;\vec \Lambda) =
    \sum_{\vec m\in
    \ZZ^n[\m]} \F_{\vec m}^{\vec m'}(\lambda,x;\vec \Lambda) \, v^{(\vec m)}_{\vec\Lambda}.
\end{equation*}
We call $\F(\lambda,x;\vec \Lambda)$ the \textbf{\emph{universal
trace matrix}} at level $\m$. For a fixed $x$ it can be thought of
as an operator
\begin{equation*}
    \F(\lambda,x;\vec \Lambda): M_{\vec\Lambda}[0] \to \Fun \o
    M_{\vec\Lambda}[0],
    \qquad
    v^{(m')}_{\vec\Lambda} \to \F^{\vec m'}(\lambda,x;\vec \Lambda).
\end{equation*}

\begin{thm}
\label{thm:resonance relations for traces} Let $\vec\Lambda \in
\C^n[2\m]$, and let $x \in \C$ be generic. Then the vector-valued
functions $\{\F^{\vec m'}(\lambda,x;\vec\Lambda)\}$ satisfy the
$n$-point resonance relations with respect to $\vec\Lambda$.
\end{thm}

The case of one-point resonance relations was treated in \cite{ES}.
Therefore, we assume that $n>1$, and prove that the columns of the
universal trace matrix satisfy the $n$-point resonance relations.

Let $\vec m \in \ZZ^n$. For $j\in\{1,\dots,n-1\}$ and $\delta \in
\{1,\dots,m_j\}$, we use Lemma \ref{thm:diagram1} and Lemma
\ref{thm:diagram2} to construct the commutative diagram
$$
\begin{gathered}\xymatrix@C=-25pt{
    &
    M_{-\delta -1 +\sum_{i=j+1}^n(\Lambda_i-2m_i)}
    \ar@{.>}[d]^{\Phi_{m_{j+2},\dots,m_n}^{\Lambda_{j+2},\dots,\Lambda_n}(-\delta + \sum_{i={j+1}}^n(\Lambda_i-2m_i))}
    &
    \\
    &
    M_{-\delta - 1 + \Lambda_{j+1}-2m_{j+1}} \o M_{\Lambda_{j+2}}^\dual \o \dots \o M_{\Lambda_n}^\dual
    \ar[ld]_(0.6){\Phi_{m_{j+1}}^{\Lambda_{j+1}}(-\delta+\Lambda_{j+1}-2m_{j+1})\o 1^{n-j-1}\qquad}
    \ar[rd]^(0.6){\qquad\qquad\Phi_{m_{j+1}+\delta}^{\Lambda_{j+1}}(-\delta+\Lambda_{j+1}-2m_{j+1})\o 1^{n-j-1}}
    &
    \\
    M_{-\delta-1} \o M_{\Lambda_{j+1}}^\dual \o \dots \o M_{\Lambda_n}^\dual
    \ar[rd]_(0.35){\Phi_{m_j}^{\Lambda_j}(-\delta) \o 1^{n-j}\quad}
    \ar[rr]^{\iota(\delta)\o 1^{n-j}}
    &&
    M_{\delta-1} \o M_{\Lambda_{j+1}}^\dual \o \dots \o M_{\Lambda_n}^\dual
    \ar[ld]^(0.35){\qquad\Phi_{m_j-\delta}^{\Lambda_j}(\delta)\o 1^{n-j}}
    \\
    &
    M_{-\delta-1-\Lambda_j+2m_j} \o M_{\Lambda_j}^\dual \o \dots \o M_{\Lambda_n}^\dual
    \ar@{.>}[d]^{\Phi_{m_1,\dots,m_{j-1}}^{\Lambda_1,\dots,\Lambda_{j-1}}(-\delta-\Lambda_j+2m_j) \o 1^{n-j+1}}
    &
    \\
    &
    M_{-\delta-1 - \sum_{i=1}^j(\Lambda_i-2m_i)} \o M_{\Lambda_1}^\dual \o \dots \o M_{\Lambda_n}^\dual
    &
}\end{gathered}.
$$

\noindent Multiplying by $q^{x \bh}$ and taking the trace, we prove
the resonance relations \eqref{eq:condition C_j abstract}. Next, we
consider the case $j=n$, and consider the commutative diagram
$$
\begin{gathered}\xymatrix@C=-70pt{
    M_{-\delta-1}
    \ar[rd]_(0.4){\Phi_{m_n}^{\Lambda_n}(-\delta)}
    \ar[rr]^{\iota(\delta)}
    &&
    M_{\delta-1}
    \ar[ld]^(0.4){\Phi_{m_n-\delta}^{\Lambda_n}(\delta)}
    \\
    &
    M_{-\delta-1-\Lambda_n+2m_n} \o M_{\Lambda_n}^\dual
    \ar@{.>}[d]^{\Phi_{m_2,\dots,m_{n-1}}^{\Lambda_2,\dots,\Lambda_{n-1}}(\delta-\Lambda_n+2m_n)\o 1}
    &
    \\
    &
    M_{-\delta-1-\sum_{i=2}^n(\Lambda_i-2m_i)} \o M_{\Lambda_2}^\dual \o \dots \o M_{\Lambda_n}^\dual
    \ar[ld]_(0.7){\Phi_{m_1}^{\Lambda_1}(-\delta-\sum_{i=2}^n(\Lambda_i-2m_i)) \o 1^{n-1}\qquad\qquad}
    \ar[rd]^(0.7){\qquad\qquad\Phi_{m_1+\delta}^{\Lambda_1}(-\delta-\sum_{i=2}^n(\Lambda_i-2m_i))\o 1^{n-1}}
    &
    \\
    M_{-\delta-1} \o M_{\Lambda_1}^\dual  \o \dots \o M_{\Lambda_n}^\dual
    \ar[rr]_{\iota(\delta) \o 1^n}
    &&
    M_{\delta-1} \o M_{\Lambda_1}^\dual \o \dots \o M_{\Lambda_n}^\dual }
\end{gathered}.
$$
In particular, we see that the image of $M_{\delta-1}$ under
$\Phi_{\tau_n^\delta(\vec m)}^{\vec\Lambda}(\delta)$ is contained in
the submodule $M_{-\delta-1} \o M_{\vec\Lambda}^\dual$. Hence the
trace function $\F_{\tau_n^\delta(\vec m)}(\delta,x;\vec\Lambda)$
can be computed from the restriction of $\Phi_{\tau_n^\delta(\vec
m)}^{\vec\Lambda}(\delta)$ to the submodule $M_{-\delta-1}$, and
therefore is equal to $\F_{\vec m}(-\delta,x;\vec\Lambda)$. This
proves \eqref{eq:condition C_n abstract} and the theorem.

\subsection{Resonance matrix vs. trace functions}

For any $\vec\Lambda\in\C^n$ let $\Xi_{\vec\Lambda}\in
\End(M_{\vec\Lambda})$ denote the invertible operator, acting
diagonally in the standard monomial basis $\{v_{\vec\Lambda}^{(\vec
m)}\}$, with diagonal elements

\begin{equation}\label{eq:def diagonal operator Xi}
    \Xi^{\vec m}_{\vec \Lambda} =
    q^{\sum_{j=1}^n m_j \sum_{i=1}^{j-1}(\Lambda_i - m_i)}.
\end{equation}

\begin{thm}\label{thm:trace vs resonance}
Let $x \in \C$ be generic. Then for any $\vec\Lambda \in \C^n[2\m]$
we have
\begin{equation}\label{eq:trace functions via resonance matrix}
    \Psi(\lambda,x;\vec\Lambda) =
    \F(\lambda,x;\vec \Lambda) \ \Xi_{\vec\Lambda}.
\end{equation}
\end{thm}

\begin{proof}
Using the definitions and easy induction on $k \in\Z_{\ge0}$, we
obtain
\begin{equation*}
    \Phi_m^\Lambda(\lambda) \mathbf F^k v_{\lambda-1} =
    \frac {q^{-m\lambda} } {(q-q^{-1})^m} \
    q^{m(\Lambda-2m) + \frac{m(m+1)}2}\
    \( \mathbf F^k v_{\lambda-\Lambda+2m-1}\o q^{-k(\Lambda-2m)} u_m^\Lambda + O(q^{2\lambda}) \),
\end{equation*}
where $O(q^{2\lambda})$ stands for a vector-valued polynomial in
$q^{2\lambda}$ of degree at most $\m$ and vanishing at zero.
Iterating and using the identity
\begin{multline*}
    \sum_{j=1}^n \( m_j\sum_{i=j}^n(\Lambda_i - 2m_i) + \frac {m_j(m_j+1)}2 \)
    \\=
    -\sum_{j=1}^n  \( m_j \sum_{i=1}^{j-1}(\Lambda_i - m_i) \) +
    \frac{\(\sum_{i=1}^n m_i\)\(1+\sum_{i=1}^n m_i\)}2 +
    \(\sum_{i=1}^n m_i\)\(\sum_{i=1}^n (\Lambda_i-2m_i) \),
\end{multline*}
for every $\vec m \in \ZZ^n[\m]$ we get
\begin{equation*}
    \Phi_{\vec m}^{\vec\Lambda}(\lambda) \mathbf F^k v_{\lambda-1} =
    \frac {q^{-\m\lambda + \m(\m+1)/2} } {(q-q^{-1})^\m} \ \(\Xi^{\vec m}_{\vec
    \Lambda}\)^{-1}
    \( \mathbf F^k v_{\lambda-1} \o u_{\vec m}^{\vec\Lambda} + O(q^{2\lambda})
    \).
\end{equation*}
Multiplying by $q^{x \bh}$ and taking the trace, we obtain
\begin{equation*}
\begin{split}
    \F_{\vec m}(\lambda,x;\vec\Lambda)
    &= (q^x-q^{-x}) \, q^{-\m \lambda} \ \(\Xi^{\vec m}_{\vec
    \Lambda}\)^{-1}
    \sum_{k=0}^\infty q^{(\lambda-2k-1) \, x} \, \( u_{\vec m}^{\vec\Lambda} + O(q^{2\lambda}) \)
    \\
    &=
    q^{ \lambda (x-\m) } \ \(\Xi^{\vec m}_{\vec \Lambda}\)^{-1}
    \( u_{\vec m}^{\vec\Lambda} + O(q^{2\lambda}) \).
\end{split}
\end{equation*}
Therefore, we have the expansion
\begin{equation}\label{eq:initial term of trace function}
    \F(\lambda,x;\vec \Lambda)  = q^{\lambda(x-\m)} \
    \( \(\Xi_{\vec\Lambda}\)^{-1} + O(q^{2\lambda}) \).
\end{equation}
The two sides of \eqref{eq:trace functions via resonance matrix} are
matrix-valued trigonometric quasi-polynomials with the same initial
term of the expansion. The columns of both matrices satisfy the
resonance relations, and for generic $\vec\Lambda$ the relation
\eqref{eq:trace functions via resonance matrix} must hold by the
uniqueness property of Theorem \ref{thm:resonance n>1}. Since both
sides are holomorphic in $\vec\Lambda$, the equality must hold
identically.
\end{proof}

\ex Let $n=1$. Then $\Xi_{2\m} = 1$, and we have $\F(\lambda,x;2\m)
= \Psi(\lambda,x;2\m)$, see \eqref{eq:example Psi n=1}.

\medskip

\ex One computes that for $\Lambda_1+\Lambda_2 = 2$
\begin{equation*}
\begin{split}
    \F(\lambda,x) &=
    \begin{pmatrix} \F_{0,1}^{0,1}(\lambda,x; \vec\Lambda) & \F_{0,1}^{1,0}(\lambda,x; \vec\Lambda)\\
    \F_{1,0}^{0,1}(\lambda,x; \vec\Lambda) & \F_{1,0}^{1,0}(\lambda,x; \vec\Lambda)
    \end{pmatrix}
    \\ &=
    \frac{q^{(\lambda-1) x}}{[x+1]}
    \left\{  \begin{pmatrix}
    q^{-\Lambda_1} [\lambda+1] & 0 \\
    q^{\lambda-1} [\Lambda_2] & q^{-\Lambda_2}[\lambda-\Lambda_2+1]
    \end{pmatrix}
    - q^{2x}
    \begin{pmatrix}
    [\lambda-\Lambda_2+1] & -q^{\lambda-\Lambda_2+1}[\Lambda_1] \\ 0 & [\lambda-1]
    \end{pmatrix} \right\} .
\end{split}
\end{equation*}
This matrix, multiplied by $\Xi_{\vec\Lambda} =
\begin{pmatrix} q^{\Lambda_1} & 0 \\ 0 & 1 \end{pmatrix}$ on the right, equals the resonance matrix
\eqref{eq:example Psi n=2,m=1}.

\section{The hypergeometric construction}
\label{sec:hypergeometric}

\subsection{The hypergeometric pairing}
Let $\vec\Lambda \in \C^n, \m \in \Z_{\ge0}$ and $\vec t =
(t_1,\dots,t_\m)$. Denote
\begin{equation*}
    \Omega(\vec t; \vec \Lambda) =
    \( \prod_{i=1}^\m
    \prod_{k=1}^n
    \frac1{(t_i - q^{\Lambda_k})(q^{\Lambda_k}t_i - 1)} \)
    \( \prod_{i<j}
    \frac{(t_i - t_j)^2}{(q t_i - q^{-1} t_j)(q^{-1} t_i - q t_j)}
    \).
\end{equation*}
For any $\vec\sigma = (\sigma_0,\sigma_1,\dots,\sigma_n) \in
\ZZ^{n+1}[\m]$ we introduce the string
\begin{equation*}
    \vec t_{\vec\sigma} =
    (0,\dots,0,q^{-\Lambda_1+2(\sigma_1-1)},\dots,q^{-\Lambda_1+2},q^{-\Lambda_1},\dots,q^{-\Lambda_n+2(\sigma_n-1)},\dots,
    q^{-\Lambda_n+2},q^{-\Lambda_n}),
\end{equation*}
and for any string $\vec x = (x_1,\dots,x_\m)$ define the multiple
residue $\Res_{\vec t = \vec x}$ by
\begin{equation*}
    \Res_{\vec t = \vec x} \ A(\vec t) =  \Res_{t_1 = x_1} \( \Res_{t_2 =
    x_2} \(
    \dots \Res_{t_\m = x_\m} \ A(\vec t) \dots \) \).
\end{equation*}
The \textbf{\emph{hypergeometric pairing}} between symmetric
polynomials $f(\vec t), g(\vec t)$ is defined by
\begin{equation}\label{eq:def hypergeometric pairing}
    \mathcal I_{\vec\Lambda}(f,g) =  \sum_{\vec\sigma}
    \frac{q^{\sigma_0(\sigma_0-1)/2}}{[\sigma_0]!}f(\vec t_{\vec\sigma})g(\vec t_{\vec\sigma})
    \
    \Res_{\vec t = \vec t_{\vec\sigma}} \frac{\Omega(\vec t; \vec\Lambda)}{t_1\dots
    t_\m}.
\end{equation}
It is clear that the hypergeometric pairing is symmetric. The above
definition is an algebraic version of a contour integral formula.
\begin{lem}
In the domain $\Re \Lambda_1 \ll \dots \ll \Re \Lambda_n \ll 0$ we
have
\begin{equation*}
    \frac{1}{(2\pi\sqrt{-1})^\m} \oint_{|t_1| = \dots = |t_\m| = 1}
    f(\vec t) \ \Omega(\vec t; \vec \Lambda) \
    \bigwedge_{i=1}^\m \frac{dt_i}{t_i}
    =
    \m ! \ \sum_{\vec\sigma \in\ZZ^{n+1}[\m]} \frac{q^{\sigma_0(\sigma_0-1)/2}}{[\sigma_0]!}f(\vec t_{\vec\sigma}) \
    \Res_{\vec t = \vec t_{\vec\sigma}} \frac{\Omega(\vec t; \vec\Lambda)}{t_1\dots
    t_\m}.
\end{equation*}
for any symmetric polynomial $f(\vec t)$.
\end{lem}
\begin{proof}
We illustrate the argument for the case $n=1, \m=2$; the general
case is analogous. Poles of the integrand are determined by the
rational expression
\begin{equation*}
    \frac{\Omega(t_1,t_2;\Lambda_1)}{t_1 t_2} =
    \frac {(t_1-t_2)^2}{(t_1-q^{\Lambda_1})(q^{\Lambda_1}t_1-1)(t_2-q^{\Lambda_1})(q^{\Lambda_1}t_2-1)
    (qt_1 - q^{-1}t_2)(q^{-1}t_1 - qt_2) \, t_1 t_2}.
\end{equation*}
The condition $\Re \Lambda_1 \ll 0$ means that $|q^{\Lambda_1}| \gg
1$, and we see that the poles inside the contour $|t_2|=1$ have
first order, and are located at $t_2 = 0, t_2 = q^{-\Lambda_1}$ and
$t_2 = q^2 t_1$. Applying the residue theorem for the variable
$t_2$, and repeating the argument for $t_1$, we get
\begin{multline*}
    \frac 1{(2\pi\sqrt{-1})^2}
    \oint_{|t_1|=|t_2|=1}
    \frac{f(t_1,t_2) \, \Omega(t_1,t_2;\Lambda_1)}{t_1 t_2} \, dt_2 \,dt_1  \ = \\
    \Res_{0,0} + \Res_{q^{-\Lambda_1},0} +
    \Res_{0, q^{-\Lambda_1}} + \Res_{q^{-\Lambda_1+2},q^{-\Lambda_1}} + \Res_{0,q^2 t_1} +
    \Res_{q^{-\Lambda_1},q^2 t_1},
\end{multline*}
where $\Res_{a,b}$ is the shortened notation for $\Res_{t_1=a}
\Res_{t_2 = b} \, \frac{f(t_1,t_2) \, \Omega(t_1,t_2;\Lambda_1)}{t_1
t_2}$. Using the symmetry of the integrand under the permutation
$t_1\leftrightarrow t_2$, we see that
\begin{equation*}
    \Res_{q^{-\Lambda_1},0} = \Res_{0, q^{-\Lambda_1}},
    \qquad
    \Res_{q^{-\Lambda_1},q^2 t_1} =
    \Res_{q^{-\Lambda_1+2},q^{-\Lambda_1}},
    \qquad
    \Res_{0,q^2 t_1} = \frac{q-q^{-1}}{q+q^{-1}} \, \Res_{0,0},
\end{equation*}
and the statement follows.
\end{proof}

\ex Let $n=2$ and $\m=1$. Then the hypergeometric pairing is given
by
\begin{equation*}
    \mathcal I_{\vec\Lambda}(f,g) = \frac{f(0)g(0)}{q^{\Lambda_1+\Lambda_2}} +
    \frac{f(q^{-\Lambda_1}) g(q^{-\Lambda_1})}
    {(1-q^{\Lambda_1-\Lambda_2})(1-q^{-2\Lambda_2})(1-q^{\Lambda_1+\Lambda_2})}  +
    \frac{f(q^{-\Lambda_2}) g(q^{-\Lambda_2})}
    {(1-q^{\Lambda_2-\Lambda_1})(1-q^{-2\Lambda_1})(1-q^{\Lambda_1+\Lambda_2})}.
\end{equation*}

\subsection{The weight functions}

For any $\vec m \in \ZZ^n[\m]$ define the weight functions
$\omega_{\vec m}(\vec t; \vec \Lambda)$ by
\begin{equation*}
\begin{split}
    \omega_{\vec m}(\vec t; \vec \Lambda) &=
    \frac{[m_1]!\dots[m_n]!}{[\m]!} \
    \sum_{I_1,\dots,I_n}
    \( \prod_{1\le k<l\le n} \prod_{i \in I_k}\prod_{j \in I_l}
    \frac{q t_i - q^{-1} t_j}{t_i - t_j} \)
    \\
    &\times
    \( \prod_{k=1}^n \prod_{i \in I_k} q^{\sum_{l=1}^{k}m_l} t_i
    \prod_{l=k+1}^n (q^{\Lambda_l} t_i - 1)
    \prod_{l=1}^{k-1} (t_i - q^{\Lambda_l} ) \).
\end{split}
\end{equation*}
where the summations are performed over all partitions
$\{I_1,\dots,I_n\}$ of the set $\{1,\dots,\m\}$ into $n$ disjoint
subsets, such that $\# I_k = m_k$. One can show that the functions
$\omega_{\vec m}(\vec t;\vec\Lambda)$ depend polynomially on
variables $\vec t$, and are invariant under the permutation group
$S_\m$.

Similarly, define the weight functions $\mathcal W_{\vec m}(\vec
t;\lambda;\vec\Lambda)$, depending on a parameter $\lambda$, by
\begin{equation*}
\begin{split}
    &\mathcal W_{\vec m}(\vec t; \lambda; \vec \Lambda) =
    \frac{[m_1]!\dots[m_n]!}{[\m]!} \
    \sum_{I_1,\dots,I_n}
    \( \prod_{1\le k<l\le n} \prod_{i \in I_k}\prod_{j \in I_l}
    \frac{q t_i - q^{-1} t_j}{t_i - t_j} \)
    \times\\
    &
    \prod_{k=1}^n \( \prod_{i \in I_k}(q^{-\lambda + m_k + \sum_{l=k}^{n}(\Lambda_l-2m_l)} t_i -
    q^{\lambda + m_k + \sum_{l=1}^{k}(\Lambda_l-2m_l)})
    \prod_{l=k+1}^n \prod_{i\in I_l} (q^{\Lambda_k} t_i - 1)
    \prod_{l=1}^{k-1} \prod_{i\in I_l} (t_i - q^{\Lambda_k} ) \),
\end{split}
\end{equation*}
and the dual weight functions ${\mathcal W}^{\vec m}(\vec t; x; \vec
\Lambda)$, depending on a parameter $x$, by
\begin{equation*}
\begin{split}
    &{\mathcal W}^{\vec m}(\vec t; x; \vec \Lambda) =
    \frac{[m_1]!\dots[m_n]!}{[\m]!} \
    \sum_{I_1,\dots,I_n}
    \( \prod_{1\le l<k \le n} \prod_{i \in I_k}\prod_{j \in I_l}
    \frac{q t_i - q^{-1} t_j}{t_i - t_j} \)
    \times\\
    &
    \prod_{k=1}^n \( \prod_{i \in I_k}(q^{-x + m_k + \sum_{l=1}^k(\Lambda_l-2m_l)} t_i -
    q^{x + m_k + \sum_{l=k}^n(\Lambda_l-2m_l)})
    \prod_{l=1}^{k-1} \prod_{i\in I_l} (q^{\Lambda_k} t_i - 1)
    \prod_{l=k+1}^n \prod_{i\in I_l} (t_i - q^{\Lambda_k} ) \).
\end{split}
\end{equation*}
The functions ${\mathcal W}_{\vec m}(\vec t; \lambda; \vec \Lambda)$
and ${\mathcal W}^{\vec m}(\vec t; x; \vec \Lambda)$ are symmetric
polynomials in variables $\vec t$.

{\bf Example.} Let $n=2$ and $\m=1$. Then
\begin{align*}
    \omega_{0,1}(t;\vec\Lambda) &= q\, t \, (q^{\Lambda_1}t-1),
    \\
    \omega_{1,0}(t;\vec\Lambda) &= q \, t \, (t-q^{\Lambda_2}),
    \\
    \mathcal W_{0,1}(t;\lambda;\vec\Lambda) &= q^{-1} \,
    (q^{-\lambda+\Lambda_2}t - q^{\lambda+\Lambda_1+\Lambda_2})(q^{\Lambda_1}t-1),
    \\
    \mathcal W_{1,0}(t;\lambda;\vec\Lambda) &= q^{-1} \,
    (q^{-\lambda+\Lambda_1+\Lambda_2}t - q^{\lambda+\Lambda_1})(t-q^{\Lambda_2}),
    \\
    {\mathcal W}^{0,1}(t;x;\vec\Lambda) &= q^{-1} \,
    (q^{-x+\Lambda_1+\Lambda_2}t - q^{x+\Lambda_2})(t-q^{\Lambda_1}),
    \\
    {\mathcal W}^{1,0}(t;x;\vec\Lambda) &= q^{-1} \,
    (q^{-x+\Lambda_1}t - q^{x+\Lambda_1+\Lambda_2})(q^{\Lambda_2}t-1).
\end{align*}

\subsection{The hypergeometric qKZ matrix $\qqH(x;\vec\Lambda)$}

Let $x\in\C$ and $\vec\Lambda \in \C^n$. For $\vec m,\vec m' \in
\ZZ^n[\m]$ set
\begin{equation}\label{eq:KZ hypergeometric integral}
    \qqH_{\vec m}^{\vec m'} (x;\vec\Lambda) = \mathcal I_{\vec\Lambda}
    \( \omega_{\vec m}(\cdot \,; \vec \Lambda) ,
    {\mathcal W}^{\vec m'}(\cdot \, ; x; \vec \Lambda) \).
\end{equation}

\begin{prop}\label{thm:qKZ is triangular}
Let $x \in \C$. Then $\qqH_{\vec m}^{\vec m'}(x;\vec\Lambda) = 0$
unless $\vec m \preccurlyeq \vec m'$, and for $\vec m' = \vec m$ we
have
\begin{equation}\label{eq:qKZ matrix diagonal}
    \qqH_{\vec m}^{\vec m}(x;\vec\Lambda) =
    q^{\m(\m-1)/2} \ \prod_{i=1}^n [m_i]!  \prod_{j=0}^{m_i-1}
    \frac{q^{x+\sum_{k=i+1}^n(\Lambda_k-2m_k) - j} - q^{-x + \sum_{k=1}^{i-1}(\Lambda_k-2m_k) -\Lambda_i + j}}{q^{\Lambda_i-j}-q^{-\Lambda_i+j} }.
\end{equation}
\end{prop}
\begin{proof} One can check that for any $\vec\sigma = (\sigma_0,\dots,\sigma_n) \in \ZZ^{n+1}[\m]$ one has
\begin{equation*}
    \omega_{\vec m}(\vec t_{\vec\sigma}; \vec\Lambda) = 0
    \qquad \text{ unless } \sigma_0 = 0 \text{ and } \vec m  \preccurlyeq
    (\sigma_1,\dots,\sigma_n).
\end{equation*}
Similarly, for any $\vec\sigma = (0,\sigma_1,\dots,\sigma_n) \in
\ZZ^{n+1}[\m]$ we have
\begin{equation*}
    \mathcal W^{\vec m'}(\vec t_{\vec\sigma}; x; \vec\Lambda) = 0
    \qquad \text{ unless } \vec m' \succcurlyeq (\sigma_1,\dots,\sigma_n).
\end{equation*}
Therefore, the sum of residues in the definition \eqref{eq:def
hypergeometric pairing} of the hypergeometric pairing is reduced to
the sum over the subset of indices $\vec \sigma =
(0,\sigma_1,\dots,\sigma_n)$, satisfying $\vec m \preccurlyeq
(\sigma_1,\dots,\sigma_n) \preccurlyeq \vec m'$. The desired
vanishing property of $\qqH_{\vec m}^{\vec m'}(x;\vec\Lambda)$
follows immediately. Finally, if $\vec m' = \vec m$, the only
nonzero term in the hypergeometric pairing corresponds to
$\vec\sigma = (0,m_1,\dots,m_n)$, and a straightforward computation
of the residue yields \eqref{eq:qKZ matrix diagonal}.
\end{proof}

For each $\vec m' \in \ZZ^n$ define the vector $\qqH^{\vec m'}
(x;\vec\Lambda) \in M_{\vec\Lambda}$ by
\begin{equation}\label{eq:KZ rows and columns}
    \qqH^{\vec m'}(x;\vec\Lambda) = \sum_{\vec m\in \ZZ^n[\m]}
    \qqH_{\vec m}^{\vec m'} (x;\vec\Lambda) \, v^{(\vec m)}_{\vec\Lambda}.
\end{equation}
Let $\qqH (x;\vec\Lambda) = \{\qqH_{\vec m}^{\vec m'}
(x;\vec\Lambda)\}_{\vec m,\vec m' \in \ZZ^n[\m]}$ denote the square
matrix, formed by vectors $\qqH^{\vec m'}(x;\vec\Lambda)$ arranged
as columns. We call $\qqH (x;\vec\Lambda)$ the
\textbf{\emph{hypergeometric qKZ matrix}} at level $\m$. For any
fixed $x \in \C$ it can be regarded as a weight-preserving operator
\begin{equation}\label{}
    \qqH (x;\vec\Lambda) \in \End(M_{\vec\Lambda}),
    \qquad v_{\vec\Lambda}^{\vec m'} \mapsto \qqH^{\vec
    m'}(x;\vec\Lambda).
\end{equation}
Matrix elements of $\qqH(x;\vec\Lambda)$ are trigonometric
polynomials in $x$, and are meromorphic in the $\vec\Lambda$
variables. It follows from Proposition \ref{thm:qKZ is triangular}
that $\qqH(x;\vec\Lambda)$ is invertible for generic $x \in \C$.

\medskip \noindent
{\bf Example.} Let $n=2$ and $\m=1$. We have
\begin{equation}\label{eq:example hypergeometric qKZ n=2,m=1}
    \qqH(x; \vec\Lambda) =
    \begin{pmatrix} \qqH_{0,1}^{0,1}(x; \vec\Lambda) & \qqH_{0,1}^{1,0}(x; \vec\Lambda)\\
    \qqH_{1,0}^{0,1}(x; \vec\Lambda) & \qqH_{1,1}^{1,1}(x; \vec\Lambda)
    \end{pmatrix} =
    \begin{pmatrix} \frac{q^x - q^{-x-\Lambda_2+\Lambda_1}}{q^{\Lambda_2}-q^{-\Lambda_2}} & 0\\
    q^{-x} & \frac{q^{x+\Lambda_2} - q^{-x-\Lambda_1}}{q^{\Lambda_1}-q^{-\Lambda_1}}
    \end{pmatrix}.
\end{equation}
\qed
\subsection{The hypergeometric qKZB matrix $\dH (\lambda,x;\vec\Lambda)$}

Let $\lambda,x\in\C$ and $\vec\Lambda \in \C^n$. For any $\vec
m,\vec m' \in \ZZ^n[\m]$, we set
\begin{equation}\label{eq:KZB hypergeometric integral}
    \dH_{\vec m}^{\vec m'} (\lambda,x;\vec\Lambda) = q^{\lambda x} \, \mathcal I_{\vec\Lambda}
    \( \mathcal W_{\vec m}(\cdot \,; \lambda; \vec \Lambda) , {\mathcal W}^{\vec m'}(\cdot \, ; x; \vec \Lambda)
    \).
\end{equation}
For each $x\in\C$ and $\vec m' \in \ZZ^n$ define the vector
$\dH^{\vec m'} (\lambda,x;\vec\Lambda) \in \Fun \o M_{\vec\Lambda}$
by
\begin{equation}\label{eq:KZB rows and columns}
    \dH^{\vec m'}(\lambda,x;\vec\Lambda) = \sum_{\vec m\in \ZZ^n[\m]}
    \dH_{\vec m}^{\vec m'} (\lambda,x;\vec\Lambda) \, v^{(\vec m)}_{\vec\Lambda}.
\end{equation}

Let $\dH (\lambda,x;\vec\Lambda) = \{\dH_{\vec m}^{\vec m'}
(\lambda,x;\vec\Lambda)\}_{\vec m,\vec m' \in \ZZ^n[\m]}$ denote the
square matrix, formed by vectors $\dH^{\vec
m'}(\lambda,x;\vec\Lambda)$ arranged as columns. We call $\dH
(\lambda,x;\vec\Lambda)$ the \textbf{\emph{hypergeometric qKZB
matrix}} at level $\m$. For any fixed $x \in \C$ it can be regarded
as a weight-preserving operator
\begin{equation}\label{eq:def hypergeometric qKZB operator}
    \dH (\lambda,x;\vec\Lambda): M_{\vec\Lambda} \to \Fun \o
    M_{\vec\Lambda},
    \qquad
    v_{\vec\Lambda}^{\vec m} \mapsto \dH^{\vec
    m}(\cdot,x;\vec\Lambda).
\end{equation}

The following theorem is a version of the more general result in
\cite{FV}.

\begin{thm}
Let $\vec\Lambda\in\C^n[2\m]$, and let $x \in \C, \, \vec
m'\in\ZZ^n[\m]$. Then the function $\dH^{\vec
m'}(\lambda,x;\vec\Lambda)$ satisfies the $n$-point resonance
relations with respect to $\vec\Lambda$.
\end{thm}
\begin{proof}
Conditions $\mathbf {C_j}(\vec m)$ for $j<n$ follow from the zero
weight condition and the identity
\begin{multline*}
    \mathcal W_{\vec m} \(\vec t; -\delta -
    \frac{\sum_{i=1}^{j}(\Lambda_i-2m_i)+ \sum_{i=j+1}^{n}(\Lambda_i-2m_i)}2;\vec\Lambda \)
    \\=
    \mathcal W_{\tau_j^\delta(\vec m)} \(\vec t;-\delta -
    \frac{\sum_{i=1}^{j}(\Lambda_i-2m_i)+ \sum_{i=j+1}^{n}(\Lambda_i-2m_i)}2;\vec\Lambda \).
\end{multline*}
Verification of the conditions $\mathbf {C_n}(\vec m)$ is more
technical. For complete details, we refer the reader to \cite{FV},
where the resonance relations for the hypergeometric qKZB matrix
were established in a more general elliptic case.
\end{proof}

The hypergeometric matrix has remarkable symmetries with respect to
variables $\lambda,x$.

\begin{lem}\label{thm:symmetry of hypergeometric qKZB}
For any $\vec m, \vec m' \in \ZZ^n[\m]$ we have
\begin{equation}\label{eq:symmetry of hypergeometric qKZB}
    \dH_{\vec m}^{\vec m'}(\lambda,x;\vec\Lambda) =
    \dH_{\vec m'}^{\vec m}(-x,-\lambda;
    \vec\Lambda)
    =
    \dH_{\Opp(\vec m)}^{\Opp(\vec m')}(-\lambda,-x;
    \Opp(\vec\Lambda)),
\end{equation}
where $\Opp(x_1,\dots,x_n) = (x_n,\dots,x_1)$.
\end{lem}
\begin{proof} Denote $\vec t^{-1} =
(t_1^{-1},\dots,t_n^{-1})$. Straightforward verification shows that
\begin{equation*}
\begin{split}
    \mathcal W_{\vec m}(\vec t; -\lambda; \vec\Lambda) &=
    t^{\m n} \ \mathcal W_{\vec m}(\vec t^{-1}; \lambda; \vec\Lambda),
    \\
    \mathcal W^{\vec m}(\vec t; -x; \vec\Lambda) &=
    t^{\m n} \ \mathcal W^{\vec m}(\vec t^{-1}; x; \vec\Lambda),
    \\
    \Omega(\vec t; \vec \Lambda) &= t^{-2 \m  n} \ \Omega(\vec t^{-1}; \vec
    \Lambda),
\end{split}
\end{equation*}
and the first of equalities \eqref{eq:symmetry of hypergeometric
qKZB} is obtained by the change of variables $\vec t \mapsto \vec
t^{-1}$ in integration over the torus $|t_1| \dots = |t_n| = 1$. The
second one follows from the symmetry of the hypergeometric pairing
and the identity
\begin{equation*}
    \mathcal W_{\vec m}(\vec t; \lambda; \vec\Lambda) =
    \mathcal W^{\Opp(\vec m)}(\vec t; \lambda; \Opp(\vec\Lambda)).
\end{equation*}
\end{proof}

{\bf Example.} Let $n=2$ and $\m = 1$. Then the hypergeometric KZB
matrix is given by
\begin{multline}\label{eq:example hypergeometric qKZB n=2,m=1}
    \dH(\lambda, x; \vec\Lambda) =
    \begin{pmatrix}
    \dH_{0,1}^{0,1}(\lambda,x; \vec\Lambda) & \dH_{0,1}^{1,0}(\lambda,x; \vec\Lambda)\\
    \dH_{1,0}^{0,1}(\lambda,x; \vec\Lambda) & \dH_{1,1}^{1,1}(\lambda,x; \vec\Lambda)
    \end{pmatrix}
    =
    q^{\lambda x + \Lambda_1+\Lambda_2-2}
    \\
    \times \left\{ q^{-\lambda}
    \begin{pmatrix}
    \frac{q^{x-\Lambda_1} - q^{-x-\Lambda_2}}{q^{\Lambda_2}-q^{-\Lambda_2}} & 0\\
    q^{-x} & \frac{q^{x+\Lambda_2} - q^{-x-\Lambda_1}}{q^{\Lambda_1}-q^{-\Lambda_1}}
    \end{pmatrix}
    + q^{\lambda}
    \begin{pmatrix} \frac{q^{-x+\Lambda_1} - q^{x-\Lambda_2}}{q^{\Lambda_2}-q^{-\Lambda_2}} & q^x\\
    0 & \frac{q^{-x-\Lambda_2} - q^{x-\Lambda_1}}{q^{\Lambda_1}-q^{-\Lambda_1}}
    \end{pmatrix} \right\}.
\end{multline}

\subsection{Trace functions vs. hypergeometric matrices}

Our next goal is to prove that the universal trace matrix
\eqref{eq:definition trace function} relates the hypergeometric qKZ
and qKZB matrices, studied in the previous subsection.

\begin{thm}\label{thm:qKZ to qKZB}
Let $\vec\Lambda \in \C^n[2\m]$, and let $x \in \C$. Then
\begin{equation}\label{eq:qKZ to qKZB}
    \dH(\lambda,x;\vec\Lambda) = \F(\lambda,x;\vec \Lambda) \
    \qqH(x;\vec\Lambda).
\end{equation}
\end{thm}

\begin{proof}

It follows from Theorem \ref{thm:resonance relations for traces}
that the columns of the matrix $\F(\lambda,x;\vec\Lambda) \,
\qqH(x;\vec\Lambda)$ satisfy the resonance relations, and therefore
for generic $x, \vec\Lambda$ it is uniquely determined by the
initial term of its expansion. Using \eqref{eq:initial term of trace
function}, we see that
\begin{equation}\label{eq:initial term FH}
    \F(\lambda,x;\vec\Lambda) \, \qqH(x;\vec\Lambda) =
    q^{\lambda(x-\m)} \( \(\Xi_{\vec\Lambda}\)^{-1} \, \qqH(x;\vec\Lambda) + O(q^{2\lambda})
    \).
\end{equation}
On the other hand, it is clear from the definitions that for any
$\vec m$ we have the expansion
\begin{equation*}
    \mathcal W_{\vec m}(\vec t;\lambda;\vec\Lambda) =
    q^{- \m\lambda} \( \(\Xi_{\vec\Lambda}^{\vec m}\)^{-1} \, \omega_{\vec m}(\vec t;\vec\Lambda) + O(q^{2\lambda})
    \),
\end{equation*}
which implies that $\dH_{\vec m}(\lambda,x;\vec\Lambda)$ is a
vector-valued function with expansion
\begin{equation}\label{eq:initial term HH}
    \dH_{\vec m}(\lambda,x;\vec\Lambda) =
    q^{\lambda(x-\m)} \( \(\Xi_{\vec\Lambda}^{\vec m}\)^{-1} \, \qqH_{\vec m}(x;\vec\Lambda) + O(q^{2\lambda})
    \).
\end{equation}
The uniqueness property for the solutions of resonance relations and
\eqref{eq:initial term FH}, \eqref{eq:initial term HH} imply that
for generic $x, \vec\Lambda$ we have $\F(\lambda,x;\vec\Lambda) \,
\qqH(x;\vec\Lambda) = \dH(\lambda,x;\vec\Lambda)$. Since both sides
are meromorphic in $\vec\Lambda$, the desired equality holds
identically.

\end{proof}

\section{Braiding properties of intertwining operators and the $R$-matrices}
\label{sec:R matrices}

\subsection{The quantum $R$-matrix} For any $\Lambda_1,\Lambda_2 \in \C$ define
$\qqR_{\Lambda_1,\Lambda_2}\in \End(M_{\Lambda_1} \o M_{\Lambda_2})$
by
\begin{equation}\label{eq:construction quantum R-matrix}
    \qqR_{\Lambda_1,\Lambda_2} = q^{-\frac{\Lambda_1\Lambda_2}2} \,
    \mathcal R \, \bigr|_{M_{\Lambda_1} \o M_{\Lambda_2}},
\end{equation}
where
\begin{equation}
    \mathcal R = \( \sum_{k\ge0} q^{-\frac{k(k-1)}2} \frac{(q^{-1}-q)^k}{[k]!}
    \mathbf E^k \o \mathbf F^k  \) q^{\frac{\bh \o \bh}2} .
\end{equation}
We refer to the family of operators $\qqR_{\Lambda_1,\Lambda_2}$ as
the quantum $R$-matrix. We define
\begin{equation*}
    \qqRch_{\Lambda_1,\Lambda_2} =
    \qqR_{\Lambda_2,\Lambda_1} \mathrm P : M_{\Lambda_1} \o M_{\Lambda_2} \to M_{\Lambda_2} \o
    M_{\Lambda_1},
\end{equation*}
where $\mathrm P$ is the permutation of tensor factors, i. e.
$\mathrm P(x \o y) = y \o x$.

For $m_1,m_2,n_1,n_2 \in \Z_{\ge0}$ define the matrix elements
$\(\qqR_{\Lambda_1,\Lambda_2}\)^{m_1,m_2}_{n_1,n_2}$ by the relation
\begin{equation}\label{eq:quantum matrix elements}
    \qqR_{\Lambda_1,\Lambda_2} (v^{(m_1)}_{\Lambda_1} \o v^{(m_2)}_{\Lambda_2}) =
    \sum_{n_1,n_2} \(\qqR_{\Lambda_1,\Lambda_2}\)^{m_1,m_2}_{n_1,n_2} v^{(n_1)}_{\Lambda_1} \o
    v^{(n_2)}_{\Lambda_2}\ .
\end{equation}
For $\vec\Lambda = (\Lambda_1,\dots,\Lambda_n)$ and distinct $i,j
\in \{1,\dots,n\}$, the action of the quantum $R$-matrix in the
$i$-th and $j$-th tensor factors yields an operator
$\qqR_{\Lambda_i,\Lambda_j}^{(i,j)}\in \End(M_{\vec\Lambda})$, which
we will simply denote $\qqR_{\Lambda_i,\Lambda_j}$; the relevant
superscripts can always be reconstructed from the context.

\medskip

The quantum $R$-matrix has the following properties (for more
details see e.g. \cite{L,CP}):
\begin{itemize}
\item
The quantum Yang-Baxter equation:
\begin{equation*}
    \qqR_{\Lambda_1,\Lambda_2}\
    \qqR_{\Lambda_1,\Lambda_3}\
    \qqR_{\Lambda_2,\Lambda_3}
    =
    \qqR_{\Lambda_2,\Lambda_3}\
    \qqR_{\Lambda_1,\Lambda_3}\
    \qqR_{\Lambda_1,\Lambda_2}.
\end{equation*}
\item
The weight-preserving property:
\begin{equation*}
    \(\qqR_{\Lambda_1,\Lambda_2}\)_{n_1,n_2}^{m_1,m_2} = 0,
    \qquad \text{ if } m_1+m_2 \ne n_1+n_2.
\end{equation*}
\item
The fusion compatibility:
\begin{equation}\label{eq:quantum R fusion}
\begin{split}
    \Upsilon_{\Lambda_1,\Lambda_2} \ \qqR_{\Lambda_1+\Lambda_2,\Lambda_3}
    &=
    \qqR_{\Lambda_2,\Lambda_3}\
    \qqR_{\Lambda_1,\Lambda_3}\
    \Upsilon_{\Lambda_1,\Lambda_2},
    \\
    \Upsilon_{\Lambda_2,\Lambda_3} \ \qqR_{\Lambda_1,\Lambda_2+\Lambda_3}
    &=
    \qqR_{\Lambda_1,\Lambda_2}\
    \qqR_{\Lambda_1,\Lambda_3}\
    \Upsilon_{\Lambda_2,\Lambda_3},
\end{split}
\end{equation}
where linear maps $\Upsilon_{\Lambda',\Lambda''}:
M_{\Lambda'+\Lambda''} \to M_{\Lambda'} \o M_{\Lambda''}$ are
defined by
\begin{equation}\label{eq:def quantum fusion inclusion}
    \Upsilon_{\Lambda',\Lambda''} \ v_{\Lambda'+\Lambda''}^{(m)} =
    \sum_{m'+m'' = m} q^{m''(\Lambda'-m')} \
    v_{\Lambda'}^{(m')} \o v_{\Lambda''}^{(m'')}.
\end{equation}
\item
The rationality property:
\begin{equation*}
    \(\qqR_{\Lambda_1,\Lambda_2}\)_{n_1,n_2}^{m_1,m_2} \in
    \C(q^{\Lambda_1},q^{\Lambda_2}).
\end{equation*}
\item
The vanishing property:
\begin{equation*}
    \(\qqR_{\Lambda_1,\Lambda_2}\)_{n_1,n_2}^{m_1,m_2} = 0, \qquad
    \text{if $\Lambda_i \in \{n_i,n_i+1, \dots, m_i-1\}$
    for $i=1$ or $2$}.
\end{equation*}
\end{itemize}

\medskip

The vanishing property implies that for $\Lambda_1,\Lambda_2 \in
\Z_{\ge0}$ we obtain induced admissible endomorphisms
$\qqR_{\Lambda_1,\Lambda_2}^{adm} \in \End(L_{\Lambda_1} \o
L_{\Lambda_2})$, represented by the matrix elements with
$\vec\Lambda$-admissible indices. The admissible maps
$\Upsilon^{adm}_{\Lambda',\Lambda''}:L_{\Lambda'+\Lambda''} \to
L_{\Lambda'} \o L_{\Lambda''}$ are defined as in \eqref{eq:def
quantum fusion inclusion} with summation over admissible indices,
and an analogue of \eqref{eq:quantum R fusion} holds for
$\qqR_{\Lambda_1,\Lambda_2}^{adm}$ and
$\Upsilon^{adm}_{\Lambda',\Lambda''}$.

The fusion compatibility property recovers admissible operators for
arbitrary $\Lambda_1,\Lambda_2\in \Z_{\ge0}$ from the normalization
condition: with respect to the basis $\{v_1^{(0)} \o v_1^{(0)},
v_1^{(0)} \o v_1^{(1)}, v_1^{(1)} \o v_1^{(0)}, v_1^{(1)} \o
v_1^{(1)}\}$ we have
\begin{equation*}
    \qqR_{1,1}^{adm} =
    \begin{pmatrix}
    1 & 0 & 0 & 0 \\
    0 & q^{-1} &  1-q^{-2} &   0 \\
    0 & 0 & q^{-1}  &   0 \\
    0 & 0 & 0 & 1
    \end{pmatrix}.
\end{equation*}
The rationality property then uniquely determines all the matrix
elements $R$ for $\Lambda_1,\Lambda_2\in\C$ as rational functions of
$q^{\Lambda_1},q^{\Lambda_2}$.

\subsection{The dynamical $R$-matrix} The elliptic $R$-matrix
$R_{\Lambda_1,\Lambda_2}^{\, \tau,\eta} (z, \lambda)\in
\End(M_{\Lambda_1} \o M_{\Lambda_2})$ was constructed in \cite{FV
elliptic} using representation theory of the elliptic quantum group
$E_{\tau,\eta}(\mathfrak{sl}_2)$. In this paper we consider
operators, which are asymptotic limits of
$R_{\Lambda_1,\Lambda_2}^{\, \tau,\eta} (z, \lambda)$. We set
\begin{equation}\label{eq:construction dynamical R-matrix}
    \dR_{\Lambda_1,\Lambda_2}(\lambda) =
    \lim_{z\to+\infty}
    \lim_{\tau \to i\infty}
    R_{\Lambda_1,\Lambda_2}^{\,\tau,\eta}(2\eta z, 2\eta \lambda).
\end{equation}

Thus, the dynamical $R$-matrix is the distinguished family
$\dR_{\Lambda_1,\Lambda_2}(\lambda)$ of endomorphisms of tensor
products $M_{\Lambda_1} \o M_{\Lambda_2}$ for
$\Lambda_1,\Lambda_2\in\C$, meromorphically depending on a complex
parameter $\lambda$. We also consider the operators
\begin{equation*}
    \dRch_{\Lambda_1,\Lambda_2}(\lambda) =
    \dR_{\Lambda_2,\Lambda_1}(\lambda) \, \mathrm P: M_{\Lambda_1} \o
    M_{\Lambda_2} \to M_{\Lambda_2} \o
    M_{\Lambda_1}.
\end{equation*}
\begin{rem}
The operators $\dRch_{\Lambda_1,\Lambda_2}(\lambda)$ are
isomorphisms of tensor products of modules for a suitable
degeneration of the elliptic quantum group
$E_{\tau,\eta}(\mathfrak{sl}_2)$. \qed
\end{rem}

For $m_1,m_2,n_1,n_2 \in \Z_{\ge0}$ define the matrix elements
$\(\dR_{\Lambda_1,\Lambda_2}(\lambda)\)^{m_1,m_2}_{n_1,n_2}$  by
\begin{equation*}\label{eq:dynamical matrix elements}
    \dR_{\Lambda_1,\Lambda_2}(\lambda) (v^{(m_1)}_{\Lambda_1} \o v^{(m_2)}_{\Lambda_2}) =
    \sum_{n_1,n_2} \(\dR_{\Lambda_1,\Lambda_2}(\lambda)\)_{n_1,n_2}^{m_1,m_2} v^{(n_1)}_{\Lambda_1} \o
    v^{(n_2)}_{\Lambda_2}\ .
\end{equation*}
We use the dynamical notation $\dR_{\Lambda_i,\Lambda_j}(\lambda -
\bh^{(k)})$ for the dynamical $R$-matrix acting in $i$-th and $j$-th
tensor factors of $M_{\vec\Lambda}$, with the convention that
$\bh^{(k)}$ must be replaced by the weight in $k$-th tensor
component. For example,
\begin{equation*}
    \dR_{\Lambda_1,\Lambda_3}(\lambda-\bh^{(2)})\
    v_{\Lambda_1}^{(m_1)} \o v_{\Lambda_2}^{(m_2)} \o
    v_{\Lambda_3}^{(m_3)} =
    \sum_{n_1,n_3}
    \(\dR_{\Lambda_1,\Lambda_3}(\lambda-\Lambda_2+2m_2)\)_{n_1,n_3}^{m_1,m_3}
    v_{\Lambda_1}^{(n_1)} \o v_{\Lambda_2}^{(m_2)} \o
    v_{\Lambda_3}^{(n_3)}.
\end{equation*}

\medskip

The dynamical $R$-matrix has the following properties, derived from
the results in \cite{FV elliptic}:

\begin{itemize}
\item
the quantum dynamical Yang-Baxter equation:
\begin{equation*}
    \dR_{\Lambda_1,\Lambda_2}(\lambda-\bh^{(3)}) \
    \dR_{\Lambda_1,\Lambda_3}(\lambda)\
    \dR_{\Lambda_2,\Lambda_3}(\lambda-\bh^{(1)})
    =
    \dR_{\Lambda_2,\Lambda_3}(\lambda)\
    \dR_{\Lambda_1,\Lambda_3}(\lambda-\bh^{(2)})\
    \dR_{\Lambda_1,\Lambda_2}(\lambda).
\end{equation*}

\item
The weight-preserving property
\begin{equation*}
    \(\dR_{\Lambda_1,\Lambda_2}(\lambda)\)_{n_1,n_2}^{m_1,m_2} = 0,
    \qquad \text{ if } m_1+m_2 \ne n_1+n_2.
\end{equation*}
\item
The fusion compatibility:
\begin{equation}\label{eq:dynamical R fusion}
\begin{split}
    \gamma_{\Lambda_1,\Lambda_2} \ \dR_{\Lambda_1+\Lambda_2,\Lambda_3}(\lambda)
    &=
    \dR_{\Lambda_2,\Lambda_3}(\lambda)\
    \dR_{\Lambda_1,\Lambda_3}(\lambda-\bh^{(2)})\
    \gamma_{\Lambda_1,\Lambda_2},
    \\
    \gamma_{\Lambda_2,\Lambda_3} \ \dR_{\Lambda_1,\Lambda_2+\Lambda_3}(\lambda)
    &=
    \dR_{\Lambda_1,\Lambda_2}(\lambda-\bh^{(3)})\
    \dR_{\Lambda_1,\Lambda_3}(\lambda)\
    \gamma_{\Lambda_2,\Lambda_3},
\end{split}
\end{equation}
where $\gamma_{\Lambda',\Lambda''}: M_{\Lambda'+\Lambda''} \to
M_{\Lambda'} \o M_{\Lambda''}$ are the linear maps determined by
\begin{equation}\label{eq:def fusion inclusion}
    \gamma_{\Lambda',\Lambda''} \ v_{\Lambda'+\Lambda''}^{(m)} =
    \sum_{m'+m'' = m}   v_{\Lambda'}^{(m')} \o
    v_{\Lambda''}^{(m'')}.
\end{equation}
\item
The rationality property:
\begin{equation*}
    \(\dR_{\Lambda_1,\Lambda_2}(\lambda)\)_{n_1,n_2}^{m_1,m_2} \in
    \C(q^{\Lambda_1},q^{\Lambda_2},q^{2\lambda}).
\end{equation*}
\item
The vanishing property:
\begin{equation*}
    \(\dR_{\Lambda_1,\Lambda_2}(\lambda)\)_{n_1,n_2}^{m_1,m_2} = 0
    \qquad
    \text{if $\Lambda_i \in \{n_i,n_i+1, \dots, m_i-1\}$
    for $i=1$ or $2$}.
\end{equation*}

\end{itemize}

The vanishing property implies that for $\Lambda_1,\Lambda_2 \in
\Z_{\ge0}$ we obtain induced admissible endomorphisms
$\dR_{\Lambda_1,\Lambda_2}^{adm}(\lambda) \in \End(L_{\Lambda_1} \o
L_{\Lambda_2})$, represented by the matrix elements with
$\vec\Lambda$-admissible indices. An analogue of \eqref{eq:dynamical
R fusion} holds for $\dR_{\Lambda_1,\Lambda_2}^{adm}(\lambda)$ and
$\gamma^{adm}_{\Lambda',\Lambda''}$, which are defined as in
\eqref{eq:def fusion inclusion} with summation over admissible
indices.

\medskip

The above properties uniquely determine the dynamical $R$-matrix
from the normalization condition: the fundamental matrix
$\dR_{1,1}^{adm}(\lambda)\in \End(L_1 \o L_1)$ is given by
\begin{equation*}
    \dR_{1,1}^{adm}(\lambda) =
    \begin{pmatrix}
    1 & 0 & 0 & 0 \\
    0 & \frac {q^{-1} [\lambda+1]}{[\lambda]}  &  -\frac {q^{-\lambda-1}}{[\lambda]}  &   0 \\
    0 & \frac {q^{\lambda-1}}{[\lambda]}  & \frac {q^{-1} [\lambda-1]}{[\lambda]}   &   0 \\
    0 & 0 & 0 & 1
    \end{pmatrix},
\end{equation*}
with respect to the basis $\{v_1^{(0)} \o v_1^{(0)}, v_1^{(0)} \o
v_1^{(1)}, v_1^{(1)} \o v_1^{(0)}, v_1^{(1)} \o v_1^{(1)}\}$.

\subsection{Gauge equivalence of the $R$-matrices}

The dynamical $R$-matrix $\dR_{\Lambda_1,\Lambda_2}(\lambda)$ has a
limit as $\lambda \to -\infty$. We set
\begin{equation}
\label{eq:asymptotical R-matrix}
    \qR_{\Lambda_1,\Lambda_2} =
    \lim_{\lambda \to -\infty} \dR_{\Lambda_1,\Lambda_2}(\lambda).
\end{equation}
In other words, the matrix elements
$\(\qR_{\Lambda_1,\Lambda_2}\)_{n_1,n_2}^{m_1,m_2}$ are obtained
from the matrix elements
$\(\dR_{\Lambda_1,\Lambda_2}(\lambda)\)_{n_1,n_2}^{m_1,m_2}$ by
taking the limit $q^{2\lambda}\to 0$ of the corresponding rational
functions of $q^{2\lambda}$, regarded as an independent variable.
The properties of $\qR_{\Lambda_1,\Lambda_2}\in\End(M_{\Lambda_1} \o
M_{\Lambda_2})$, inherited from the dynamical $R$-matrix, are the
same as those of $\qqR_{\Lambda_1,\Lambda_2}$, but with the modified
fusion compatibility property
\begin{equation}\label{eq:asymptotic R fusion}
\begin{split}
    \gamma_{\Lambda_1,\Lambda_2} \ \qR_{\Lambda_1+\Lambda_2,\Lambda_3}
    &=
    \qR_{\Lambda_2,\Lambda_3}\
    \qR_{\Lambda_1,\Lambda_3}\
    \gamma_{\Lambda_1,\Lambda_2},
    \\
    \gamma_{\Lambda_2,\Lambda_3} \ \qR_{\Lambda_1,\Lambda_2+\Lambda_3}
    &=
    \qR_{\Lambda_1,\Lambda_2}\
    \qR_{\Lambda_1,\Lambda_3}\
    \gamma_{\Lambda_2,\Lambda_3},
\end{split}
\end{equation}
and the normalization condition for $\qR_{1,1}^{adm}\in \End(L_1 \o
L_1)$
\begin{equation*}
    \qR_{1,1}^{adm} =
    \begin{pmatrix}
    1 & 0 & 0 & 0 \\
    0 & q^{-2} &  1-q^{-2} &   0 \\
    0 & 0 & 1  &   0 \\
    0 & 0 & 0 & 1
    \end{pmatrix}.
\end{equation*}
As in the case of the quantum $R$-matrix, the inductive fusion
procedure recovers the matrix elements
$\(\qR_{\Lambda_1,\Lambda_2}\)_{n_1,n_2}^{m_1,m_2}$ for positive
integral $\Lambda_1,\Lambda_2$ and admissible $(m_1,m_2)$,
$(n_1,n_2)$; the rationality property then determines all the matrix
elements for arbitrary $\Lambda_1,\Lambda_2$.

\medskip

We call the collection of operators $\qR_{\Lambda_1,\Lambda_2}$ the
asymptotic $R$-matrix. Below we show that it is gauge equivalent to
the quantum $R$-matrix $\qqR_{\Lambda_1,\Lambda_2}$, and that the
equivalence is given by the diagonal operators
$\Xi_{\Lambda_1,\Lambda_2}$ defined by \eqref{eq:def diagonal
operator Xi}.

\begin{prop}\label{thm:gauge equivalence for R matrices}
For any $\Lambda_1,\Lambda_2\in\C$ we have
\begin{equation}\label{eq:gauge equivalence for R matrices}
    \qqRch_{\Lambda_1,\Lambda_2} \
    \Xi_{\Lambda_1,\Lambda_2} =
    \Xi_{\Lambda_2,\Lambda_1} \
    \qRch_{\Lambda_1,\Lambda_2}.
\end{equation}
\end{prop}
\begin{proof}
The desired formula is equivalent to an infinite collection of
identities between matrix elements of $R$-matrices. Since the matrix
elements of both $R$-matrices are rational functions of
$q^{\Lambda_1},q^{\Lambda_2}$, it suffices to check each of these
identities for all sufficiently large positive integers
$\Lambda_1,\Lambda_2$. To achieve this goal, we establish the
admissible version of \eqref{eq:gauge equivalence for R matrices} by
using double induction in $\Lambda_1$ and $\Lambda_2$.

The base of induction is $\Lambda_1 = \Lambda_2 = 1$, and a
straightforward computation settles the case of the gauge
equivalence of the fundamental $R$-matrices, cf. the example below.

To prove the inductive step, we start with the following identities:
\begin{equation}\label{eq:Xi intertwining property}
\begin{split}
    \Xi_{\Lambda_1,\Lambda_2,\Lambda_3}\
    \gamma_{\Lambda_1,\Lambda_2}
    &=
    \Upsilon_{\Lambda_1,\Lambda_2}\
    \Xi_{\Lambda_1+\Lambda_2,\Lambda_3},
    \\
    \Xi_{\Lambda_1,\Lambda_2,\Lambda_3} \
    \gamma_{\Lambda_2,\Lambda_3}
    &=
    \Upsilon_{\Lambda_2,\Lambda_3}\
    \Xi_{\Lambda_1,\Lambda_2+\Lambda_3},
\end{split}
\end{equation}
which are verified directly from the definitions. Assuming now that
that \eqref{eq:gauge equivalence for R matrices} is true for pairs
$(\Lambda_1,\Lambda_3)$ and $(\Lambda_2,\Lambda_3)$, we get from the
explicit formula \eqref{eq:def diagonal operator Xi} and the
weight-preserving properties of the $R$-matrices
\begin{equation}\label{eq:gauge assumption}
\begin{split}
    \qqRch_{\Lambda_1,\Lambda_3}\
    \Xi_{\Lambda_1,\Lambda_3,\Lambda_2}
    &=
    \Xi_{\Lambda_3,\Lambda_1,\Lambda_2}\
    \qRch_{\Lambda_1,\Lambda_3},
    \\
    \qqRch_{\Lambda_2,\Lambda_3}\
    \Xi_{\Lambda_1,\Lambda_2,\Lambda_3}
    &=
    \Xi_{\Lambda_1,\Lambda_3,\Lambda_2}\
    \qRch_{\Lambda_2,\Lambda_3}.
\end{split}
\end{equation}
Since $\Upsilon_{\Lambda_1,\Lambda_2}$ is injective, the computation
\begin{align}
    \Upsilon_{\Lambda_1,\Lambda_2} \
    \qqRch_{\Lambda_1+\Lambda_2,\Lambda_3}\
    \Xi_{\Lambda_1+\Lambda_2,\Lambda_3}
    &=
    \qqRch_{\Lambda_1,\Lambda_3}\
    \qqRch_{\Lambda_2,\Lambda_3}\
    \Upsilon_{\Lambda_1,\Lambda_2}\
    \Xi_{\Lambda_1+\Lambda_2,\Lambda_3}
    \tag{due to \eqref{eq:quantum R fusion}}
    \\
    &=
    \qqRch_{\Lambda_1,\Lambda_3}\
    \qqRch_{\Lambda_2,\Lambda_3}\
    \Xi_{\Lambda_1,\Lambda_2,\Lambda_3}\
    \gamma_{\Lambda_1,\Lambda_2}
    \tag{due to \eqref{eq:Xi intertwining property}}
    \\
    &=
    \Xi_{\Lambda_3,\Lambda_1,\Lambda_2}\
    \qRch_{\Lambda_1,\Lambda_3}\
    \qRch_{\Lambda_2,\Lambda_3}\
    \gamma_{\Lambda_1,\Lambda_2}
    \tag{due to \eqref{eq:gauge assumption}}
    \\
    &=
    \Xi_{\Lambda_3,\Lambda_1,\Lambda_2}\
    \gamma_{\Lambda_1,\Lambda_2}\
    \qRch_{\Lambda_1+\Lambda_2,\Lambda_3}\
    \tag{due to \eqref{eq:asymptotic R fusion}}
    \\
    &=
    \Upsilon_{\Lambda_1,\Lambda_2} \
    \Xi_{\Lambda_3,\Lambda_1+\Lambda_2}\
    \qRch_{\Lambda_1+\Lambda_2,\Lambda_3},
    \tag{due to \eqref{eq:Xi intertwining property}}
\end{align}
shows that \eqref{eq:gauge equivalence for R matrices} holds for the
pair $(\Lambda_1+\Lambda_2,\Lambda_3)$. Note also that all of the
above formulae remain valid when the operators are replaced by their
admissible versions.

We conclude that the admissible version of \eqref{eq:gauge
equivalence for R matrices} holds when $\Lambda_1$ is any positive
integer and $\Lambda_2=1$. A similar inductive argument in
$\Lambda_2$ completes the proof.
\end{proof}

\begin{rem}
The operators $\Xi_{\vec\Lambda}$ are characterized by the property
$\Upsilon_{\vec\Lambda} = \Xi_{\vec\Lambda} \ \gamma_{\vec\Lambda}$,
where $\Upsilon_{\vec\Lambda}$ and $\gamma_{\vec\Lambda}$ are
inclusions of modules for the quantum group and a suitable version
of the dynamical quantum group, respectively. Thus they can be
regarded as "intertwining" operators between the two categories of
representations.
\end{rem}

\ex Consider the level $\m=1$ weight subspace, spanned by
$\{v_{\Lambda_1}^{(0)} \o v_{\Lambda_2}^{(1)}, v_{\Lambda_1}^{(1)}
\o v_{\Lambda_2}^{(0)}\}$. The dynamical $R$-matrix and its
asymptotic version are given by
\begin{equation*}
    \dR_{\Lambda_1,\Lambda_2}(\lambda) =
    \begin{pmatrix}
    \frac {q^{-\Lambda_1} [\lambda+1]}{[\lambda-\Lambda_1+1]} &
    -\frac {q^{-\lambda-1}[\Lambda_1]}{[\lambda-\Lambda_1+1]} \\
    \frac {q^{\lambda-\Lambda_1-\Lambda_2+1} [\Lambda_2]}{[\lambda-\Lambda_1+1]} &
    \frac {q^{-\Lambda_2} [\lambda-\Lambda_1-\Lambda_2+1]}{[\lambda-\Lambda_1+1]}
    \end{pmatrix},
    \qquad
    \qR_{\Lambda_1,\Lambda_2} =
    \begin{pmatrix}
    q^{-2\Lambda_1} && 1 - q^{-2\Lambda_1} \\
    0 && 1
    \end{pmatrix},
\end{equation*}
and the quantum $R$-matrix is given by
\begin{equation}\label{eq:example quantum R matrix}
    \qqR_{\Lambda_1,\Lambda_2} =
    \begin{pmatrix}
    q^{-\Lambda_1} & q^{-\Lambda_2}(q^{\Lambda_1} - q^{-\Lambda_1}) \\
    0 & q^{-\Lambda_2}
    \end{pmatrix}.
\end{equation}
The gauge equivalence between $\qqR_{\Lambda_1,\Lambda_2}$ and
$\qR_{\Lambda_1,\Lambda_2}$ is illustrated by the matrix identity
\begin{equation*}
    \begin{pmatrix}
    q^{-\Lambda_1} & q^{-\Lambda_2}(q^{\Lambda_1} - q^{-\Lambda_1}) \\
    0 & q^{-\Lambda_2}
    \end{pmatrix}
    \begin{pmatrix}
    1 & 0 \\
    0 & q^{\Lambda_2}
    \end{pmatrix}
    =
    \begin{pmatrix}
    q^{\Lambda_1} & 0 \\
    0 & 1
    \end{pmatrix}
    \begin{pmatrix}
    q^{-2\Lambda_1} && 1 - q^{-2\Lambda_1} \\
    0 && 1
    \end{pmatrix}.
\end{equation*}

\subsection{The fusion lemma}

Let $\pi_{\Lambda_1,\Lambda_2}: M_{\Lambda_1}^\dual \o
M_{\Lambda_2}^\dual \to M_{\Lambda_1+\Lambda_2}^\dual$ be the unique
$\U$-module surjection, normalized by
$\pi_{\Lambda_1,\Lambda_2}(u^{\Lambda_1}_0 \o u^{\Lambda_2}_0) =
u^{\Lambda_1+\Lambda_2}_0$. More explicitly, for
$m_1,m_2\in\Z_{\ge0}$ we have
\begin{equation}\label{eq:definition pi for two}
    \pi_{\Lambda_1,\Lambda_2}
    (u^{\Lambda_1}_{m_1} \o u^{\Lambda_2}_{m_2}) =
    q^{(\Lambda_1-m_1)m_2} \
    u^{\Lambda_1+\Lambda_2}_{m_1+m_2}.
\end{equation}

\begin{lem} \label{thm:highest component fusion}
For any $m_1,m_2 \in \Z_{\ge0}$ we have a commutative diagram
\begin{equation}
\label{eq:two component fusion}
\begin{gathered}\xymatrix{
    &
    M_{\lambda-1}
    \ar[ld]_{\Phi_{m_1,m_2}^{\Lambda_1,\Lambda_2}(\lambda)}
    \ar[rd]^{\Phi_{m_1+m_2}^{\Lambda_1+\Lambda_2}(\lambda)}
    &
    \\
    M_{\lambda-\Lambda_1-\Lambda_2+2m_1+2m_2-1} \o
    M_{\Lambda_1}^\dual \o M_{\Lambda_2}^\dual
    \ar[rr]_{1\o\pi_{\Lambda_1,\Lambda_2}}
    &&
    M_{\lambda-\Lambda_1-\Lambda_2+2m_1+2m_2-1} \o M_{\Lambda_1+\Lambda_2}^\dual
}\end{gathered}.
\end{equation}
\end{lem}

\begin{proof}
Set $\mu = \lambda - \Lambda_2 + 2m_2, \ \nu = \lambda - \Lambda_1 -
\Lambda_2 + 2m_1 + 2m_2$. We compute
\begin{equation}\label{eq:fusion computation}
\begin{split}
    & \Phi_{m_1,m_2}^{\Lambda_1,\Lambda_2} (\lambda) v_{\lambda-1} =
    (\Phi_{m_1}^{\Lambda_1}(\mu) \o 1)
    \( \sum_{i=0}^{m_2} q^{i(\mu-i)} \frac{[-\mu+i+1]_{m_2-i}}{[i]!} \
    \mathbf F^i v_{\mu-1} \o \mathbf E^i u_{m_2}^{\Lambda_2} \)
    \\
    &= \sum_{i=0}^{m_2} q^{i(\mu-i)} \frac{[-\mu+i+1]_{m_2-i}}{[i]!} \ \Delta(\mathbf F^i) \(
    [-\nu+1]_{m_1} v_{\nu-1} \o u_{m_1}^{\Lambda_1} + \dots \) \o \mathbf E^i u_{m_2}^{\Lambda_2} \\
    &=  v_{\nu-1} \o [-\nu+1]_{m_1} \, \( \sum_{i=0}^{m_2}  q^{i(\mu-i)} \frac{[-\mu+i+1]_{m_2-i}}{[i]!} \,
    \mathbf F^i u_{m_1}^{\Lambda_1} \o \mathbf E^i u_{m_2}^{\Lambda_2} \)
    + \dots,
\end{split}
\end{equation}
where we omitted the terms involving $\mathbf F^k v_{\nu-1}$ with
$k>0$. We now study the leading term. Denote for brevity
\begin{equation*}
    u = [-\nu+1]_{m_1} \, [-\mu+1]_{m_2} \ \sum_{i=0}^{m_2}  \frac{q^{i(\mu-i)}}{[-\mu+1]_i [i]!} \,
    \mathbf F^i u_{m_1}^{\Lambda_1} \o \mathbf E^i u_{m_2}^{\Lambda_2}.
\end{equation*}
Explicit formula \eqref{eq:definition pi for two} for the map
$\pi_{\Lambda_1,\Lambda_2}$ yields
\begin{equation*}
    \pi_{\Lambda_1,\Lambda_2}(\mathbf F^i u_{m_1}^{\Lambda_1} \o \mathbf E^i u_{m_2}^{\Lambda_2}) =
    q^{(\Lambda_1-m_1-i)(m_2-i)}[-m_2]_i \, [m_1-\Lambda_1]_i \ u_{m_1+m_2}^{\Lambda_1+\Lambda_2},
\end{equation*}
and combining it with the summation formula
\begin{equation*}
    \sum_{i=0}^{m_2} q^{i(-m_2+m_1-\Lambda_1+\mu)} \frac{[-m_2]_i [m_1-\Lambda_1]_i}{[-\mu+1]_i [i]!} =
    q^{m_2(m_1-\Lambda_1)} \frac{[-\mu+\Lambda_1-m_1+1]_{m_2}}{[-\mu+1]_{m_2}},
\end{equation*}
which is a special case of the $q$-hypergeometric Gauss identity, we
conclude that
\begin{equation*}
\begin{split}
    \pi_{\Lambda_1,\Lambda_2} (u) &= [-\nu+1]_{m_1} \,
    [-\mu+\Lambda_1-m_1+1]_{m_2} \, u_{m_1+m_2}^{\Lambda_1+\Lambda_2}
    \\
    &=
    [-\nu+1]_{m_1} \, [-\nu+m_1+1]_{m_2} \, u_{m_1+m_2}^{\Lambda_1+\Lambda_2}
    = [-\nu+1]_{m_1+m_2} \, u_{m_1+m_2}^{\Lambda_1+\Lambda_2}.
\end{split}
\end{equation*}
Now \eqref{eq:fusion computation} yields
\begin{equation*}
    (1 \o \pi_{\Lambda_1,\Lambda_2}) \ \Phi_{m_1,m_2}^{\Lambda_1,\Lambda_2}(\lambda) v_{\lambda-1} =
    [-\nu+1]_{m_1+m_2} \, v_{\nu-1} \o  u_{m_1+m_2}^{\Lambda_1+\Lambda_2} + \dots,
\end{equation*}
which shows that the operators $(1 \o \pi_{\Lambda_1,\Lambda_2})
\Phi_{m_1,m_2}^{\Lambda_1,\Lambda_2}(\lambda)$ and
$\Phi_{m_1+m_2}^{\Lambda_1+\Lambda_2}(\lambda)$ have the same
leading terms. For generic $\lambda$ the leading term uniquely
determines the intertwining operator, proving the equality of all
other coefficients in the expansions of both operators. Since all
these coefficients are holomorphic in $\lambda$, we conclude that
\eqref{eq:two component fusion} holds for every value of $\lambda$.
\end{proof}

\begin{rem}
Lemma \ref{thm:highest component fusion} shows that for $\U$ our
regularized intertwining operators $\Phi_m^\Lambda(\lambda)$ can be
constructed inductively using fusion. We start with
$\Phi_0^1(\lambda)$ and $\Phi_1^1(\lambda)$, which determine the
operators associated with the fundamental representation
$L_1^\dual$. For nonnegative integers $\Lambda_1, \Lambda_2$ the
operators associated with $L^\dual_{\Lambda_1}$ and
$L^\dual_{\Lambda_2}$ yield the operators associated with
$L^\dual_{\Lambda_1+\Lambda_2}$, which determines
$\Phi^\Lambda_m(\lambda)$ for all large integers $\Lambda$, and
hence also for arbitrary $\Lambda \in \C$.

It is an open problem to construct similar regularized (holomorphic)
intertwining operators for higher rank quantum groups. Lemma
\ref{thm:highest component fusion} suggests that in the higher rank
case it suffices to consider the fundamental representation, and
then extend the operators to other representations by using
fusion.\qed
\end{rem}

\subsection{The braiding lemma}
Let $\qqR_{\Lambda_1,\Lambda_2}^\dual: M_{\Lambda_1}^\dual \o
M_{\Lambda_2}^\dual \to M_{\Lambda_1}^\dual \o M_{\Lambda_1}^\dual$
be the dual linear map to $\qqR_{\Lambda_1,\Lambda_2}$. In other
words, we have
\begin{equation}\label{eq:def dual quantum R}
    \qqR_{\Lambda_1,\Lambda_2}^\dual \
    u^{\Lambda_1}_{m_1} \o u^{\Lambda_2}_{m_2} =
    \sum_{n_1,n_2} \(\qqR_{\Lambda_1,\Lambda_2}\)^{n_1,n_2}_{m_1,m_2}\
    u^{\Lambda_1}_{n_1} \o u^{\Lambda_2}_{n_2},
\end{equation}
where $\(\qqR_{\Lambda_1,\Lambda_2}\)^{n_1,n_2}_{m_1,m_2}$ are
defined in \eqref{eq:quantum matrix elements}. Equivalently,
$\qqR_{\Lambda_1,\Lambda_2}^\dual = q^{-\frac{\Lambda_1\Lambda_2}2}
\, \mathcal R' \, \bigr|_{M_{\Lambda_1}^\dual \o
M_{\Lambda_2}^\dual}$, where
\begin{equation}\label{eq:universal R-matrix}
    \mathcal R' =  q^{\frac{\bh \o \bh}2}  \ \( \sum_{k\ge0} q^{-\frac{k(k-1)}2} \frac{(q^{-1}-q)^k}{[k]!}
    \mathbf F^k \o \mathbf E^k  \).
\end{equation}
\ex Consider the weight subspace of level $\m = 1$, spanned by
$\{u^{\Lambda_1}_{0} \o u^{\Lambda_2}_{1}, u^{\Lambda_1}_{1} \o
u^{\Lambda_2}_{0} \}$. Then $\qqR_{\Lambda_1,\Lambda_2}^\dual$ is
represented by the matrix, transposed to \eqref{eq:example quantum R
matrix}:
\begin{equation*}
    \qqR_{\Lambda_1,\Lambda_2}^\dual =
    \begin{pmatrix}
    q^{-\Lambda_1} & 0 \\
    q^{-\Lambda_2}(q^{\Lambda_1} - q^{-\Lambda_1}) & q^{-\Lambda_2}
    \end{pmatrix}.
\end{equation*}

\bigskip \noindent The element $\mathcal R'$, defined in
\eqref{eq:universal R-matrix}, has the property that
\begin{equation*}
    \mathcal R' \ \Delta(X) = \Delta^{op}(X) \ \mathcal R',
    \qquad X \in \U,
\end{equation*}
where $\Delta^{op}(X) = \mathrm P \  \Delta(X)$. Therefore, we
obtain $\U$-module isomorphisms
\begin{equation*}
    \qqRch_{\Lambda_1,\Lambda_2}^\dual = \mathrm P \,
    \qqR_{\Lambda_1,\Lambda_2}^\dual: M_{\Lambda_1}^\dual \o
    M_{\Lambda_2}^\dual \to M_{\Lambda_2}^\dual \o
    M_{\Lambda_1}^\dual,
\end{equation*}
with respect to the $\U$-action in tensor products, defined using
the comultiplication \eqref{eq:def quantum comultiplication}. The
compositions $\qqRch_{\Lambda_1,\Lambda_2}^\dual
\Phi^{\Lambda_1,\Lambda_2}_{m_1,m_2}(\lambda): M_{\lambda-1} \to
M_{\lambda-\Lambda_1-\Lambda_2+2m_1+2m_2} \o M_{\Lambda_2}^\dual \o
M_{\Lambda_1}^\dual$ are $\U$-intertwining operators, and can be
written as linear combinations of operators
$\Phi^{\Lambda_2,\Lambda_1}_{n_2,n_1}(\lambda)$. We now show that
the corresponding connection matrix is precisely equal to
$\dR_{\Lambda_1,\Lambda_2}(\lambda)$.


\begin{lem}\label{thm:n=2 braiding}
For any $\Lambda_1, \Lambda_2 \in \C$ and $m_1,m_2 \in \Z_{\ge0}$ we
have
\begin{equation}\label{eq:coordinate braiding}
    \qqRch_{\Lambda_1,\Lambda_2}^\dual\
    \Phi^{\Lambda_1,\Lambda_2}_{m_1,m_2}(\lambda) =
    \sum_{n_1,n_2} \(\dR_{\Lambda_1,\Lambda_2}(\lambda)\)^{n_1,n_2}_{m_1,m_2}\
    \Phi^{\Lambda_2,\Lambda_1}_{n_2,n_1}(\lambda).
\end{equation}

\end{lem}

\begin{proof}
It is convenient to rewrite \eqref{eq:coordinate braiding} using the
formal generating functions
\begin{equation}\label{eq:formal intertwiner}
    \Phi^{\vec\Lambda}(\lambda) = \sum_{\vec m}
    \Phi^{\vec\Lambda}_{\vec m}(\lambda) \o v_{\vec\Lambda}^{(\vec m)}.
\end{equation}
These generating functions have a useful property
\begin{equation}\label{eq:krivoi fusion}
    \pi_{\Lambda_1,\Lambda_2} \
    \Phi^{\Lambda_1,\Lambda_2}(\lambda) =
    \gamma_{\Lambda_1,\Lambda_2}\
    \Phi^{\Lambda_1+\Lambda_2}(\lambda),
\end{equation}
established by the following computation based on Lemma
\ref{thm:highest component fusion}:
\begin{multline*}
    \sum_{m_1,m_2} \pi_{\Lambda_1,\Lambda_2}
    \Phi^{\Lambda_1\Lambda_2}(\lambda) \o v^{(m_1)}_{\Lambda_1} \o
    v^{(m_2)}_{\Lambda_2} =
    \sum_{m_1,m_2} \Phi_{m_1+m_2}^{\Lambda_1+\Lambda_2}(\lambda) \o v^{(m_1)}_{\Lambda_1} \o
    v^{(m_2)}_{\Lambda_2}
    \\
    =
    \sum_{m}\sum_{m_1+m_2 = m} \Phi_{m}^{\Lambda_1+\Lambda_2}(\lambda) \o  v^{(m_1)}_{\Lambda_1} \o
    v^{(m_2)}_{\Lambda_2} =
    \sum_{m} \Phi_{m}^{\Lambda_1+\Lambda_2}(\lambda) \o
    \gamma_{\Lambda_1,\Lambda_2} v^{(m)}_{\Lambda_1+\Lambda_2}.
\end{multline*}
Returning to the desired relation \eqref{eq:coordinate braiding}, we
see that it is equivalent to
\begin{equation}\label{eq:n=2 braiding}
    \qqRch_{\Lambda_1,\Lambda_2}^\dual \
    \Phi^{\Lambda_1,\Lambda_2}(\lambda) =
    \dRch_{\Lambda_2,\Lambda_1}(\lambda) \
    \Phi^{\Lambda_2,\Lambda_1}(\lambda),
\end{equation}
which we prove for $\Lambda_1,\Lambda_2 \in \Z_{>0}$ using double
induction. Assuming that \eqref{eq:n=2 braiding} holds for pairs
$(\Lambda_1,\Lambda_3)$ and $(\Lambda_2,\Lambda_3)$, and using the
fusion compatibility
\begin{equation}\label{eq:pi fusion}
\begin{split}
    \qqRch_{\Lambda_1+\Lambda_2,\Lambda_3}^\dual \
    \pi_{\Lambda_1,\Lambda_2}
    &=
    \pi_{\Lambda_1,\Lambda_2}\
    \qqRch_{\Lambda_1,\Lambda_3}^\dual\
    \qqRch_{\Lambda_2,\Lambda_3}^\dual,
    \\
    \qqRch_{\Lambda_1,\Lambda_2+\Lambda_3}^\dual\
    \pi_{\Lambda_2,\Lambda_3}
    &=
    \pi_{\Lambda_2,\Lambda_3}\
    \qqRch_{\Lambda_1,\Lambda_3}^\dual\
    \qqRch_{\Lambda_1,\Lambda_2}^\dual,
\end{split}
\end{equation}
we compute
\begin{align*}
    \gamma_{\Lambda_1,\Lambda_2} \
    \qqRch_{\Lambda_1+\Lambda_2,\Lambda_3}^\dual \
    &
    \Phi^{\Lambda_1+\Lambda_2,\Lambda_3}(\lambda)
    =
    \qqRch_{\Lambda_1+\Lambda_2,\Lambda_3}^\dual \
    \pi_{\Lambda_1,\Lambda_2} \
    \Phi^{\Lambda_1,\Lambda_2,\Lambda_3}(\lambda)
    \tag{due to \eqref{eq:krivoi fusion}}
    \\
    &=
    \pi_{\Lambda_1,\Lambda_2}\
    \qqRch_{\Lambda_1,\Lambda_3}^\dual\
    \qqRch_{\Lambda_2,\Lambda_3}^\dual\
    \Phi^{\Lambda_1,\Lambda_2,\Lambda_3}(\lambda)
    \tag{due to \eqref{eq:pi fusion}}
    \\
    &=
    \pi_{\Lambda_1,\Lambda_2}\
    \dRch_{\Lambda_3,\Lambda_1}(\lambda - \bh^{(2)})\
    \dRch_{\Lambda_3,\Lambda_2}(\lambda)\
    \Phi^{\Lambda_3,\Lambda_1,\Lambda_2}(\lambda)\
    \tag{by assumption}
    \\
    &=
    \dRch_{\Lambda_3,\Lambda_1}(\lambda - \bh^{(2)})\
    \dRch_{\Lambda_3,\Lambda_2}(\lambda)\
    \gamma_{\Lambda_1,\Lambda_2}\
    \Phi^{\Lambda_3,\Lambda_1+\Lambda_2}(\lambda)\
    \tag{due to \eqref{eq:krivoi fusion}}
    \\
    &=
    \gamma_{\Lambda_1,\Lambda_2}\
    \dRch_{\Lambda_3,\Lambda_1+\Lambda_2}(\lambda) \
    \Phi^{\Lambda_3,\Lambda_1+\Lambda_2}(\lambda),
    \tag{due to \eqref{eq:dynamical R fusion}}
\end{align*}
and since $\gamma_{\Lambda_1,\Lambda_2}$ is injective, we see that
the braiding relation is also valid for
$(\Lambda_1+\Lambda_2,\Lambda_3)$. Similarly, if \eqref{eq:n=2
braiding} is satisfied for pairs $(\Lambda_1,\Lambda_2)$ and
$(\Lambda_1,\Lambda_3)$, then it holds for
$(\Lambda_1,\Lambda_2+\Lambda_3)$. When $\Lambda_1,\Lambda_2 \in
\Z_{\ge0}$, analogous conclusions are obviously true for the
admissible operators.

\smallskip

To complete the proof, we note that braiding relation holds in the
admissible case for $\Lambda_1 = \Lambda_2 =1$ (which is the result
of direct computation, see example below), and by induction for
arbitrary $\Lambda_1,\Lambda_2 \in \Z_{>0}$. In the general case,
\eqref{eq:n=2 braiding} is equivalent to a family of equalities of
rational functions, which are now established for
$\Lambda_1,\Lambda_2$ in the positive integral cone, and therefore
must hold identically.
\end{proof}

{\bf Example.} Let $\m = 1$. Then \eqref{eq:n=2 braiding} reduces to
the matrix identity
\begin{multline*}
    \begin{pmatrix}
    -[\lambda-\Lambda_2+1] & q^{\lambda-\Lambda_2+1}[\Lambda_1]  \\
    0 & -[\lambda-\Lambda_1-\Lambda_2+1]
    \end{pmatrix}
    \begin{pmatrix}
    q^{-\Lambda_1} & q^{-\Lambda_2}(q^{\Lambda_1}-q^{-\Lambda_1}) \\
    0 &  q^{-\Lambda_2}
    \end{pmatrix}
    \\
    =
    \begin{pmatrix}
    \frac{q^{-\Lambda_1}[\lambda+1]}{[\lambda-\Lambda_1+1]} &
    -\frac{q^{-\lambda-1}[\Lambda_1]}{[\lambda-\Lambda_1+1]}
    \\
    \frac{q^{\lambda-\Lambda_1-\Lambda_2+1}[\Lambda_2]}{[\lambda-\Lambda_1+1]} &
    \frac{q^{-\Lambda_1} [\lambda-\Lambda_1-\Lambda_2+1]}{[\lambda-\Lambda_1+1]}
    \end{pmatrix}
    \begin{pmatrix}
    -[\lambda-\Lambda_1-\Lambda_2+1] & 0  \\
    q^{\lambda-\Lambda_1+1}[\Lambda_2] & -[\lambda-\Lambda_1+1]
    \end{pmatrix}.
\end{multline*}

\section{Difference equations satisfied by the trace functions}
\label{sec:equations}

\subsection{The qKZ and qKZB operators}
For any $x \in \C$ define $\(q^{x\bh}\)_j \in \End(M_{\vec\Lambda})$
by
\begin{equation*}
    \(q^{x\bh}\)_j  \ v^{(\vec m)}_{\vec\Lambda}  =
    q^{(\Lambda_j-2m_j)x}\ v^{(\vec m)}_{\vec\Lambda},
\end{equation*}
and $\Gamma_j \in \End(\Fun \o M_{\vec\Lambda})$ by
\begin{equation*}
    \Gamma_j \ \( \xi(\lambda) \o v^{(\vec m)}_{\vec\Lambda} \) =
    \xi (\lambda - \Lambda_1 + 2k_1) \o v^{(\vec m)}_{\vec\Lambda},
\end{equation*}
for any $\vec m \in \ZZ^n$ and $\xi\in\Fun$. We also set
\begin{equation*}
    \(q^{x\bh}\)_\circlearrowleft =
    \mathrm P_\circlearrowleft \ \(q^{x\bh}\)_1,
    \qquad\qquad
    \Gamma_\circlearrowleft = \mathrm P_\circlearrowleft \ \Gamma_1,
\end{equation*}
where $\mathrm P_\circlearrowleft$ is the permutation map, defined
by $\mathrm P_\circlearrowleft \( x_1 \o x_2 \o \dots \o x_n \) =
    x_2 \o \dots \o x_n \o x_1$.

For any $x \in \C$ we define the qKZ operators by
\begin{equation*}
    \qqKZ_j(x) = \( \qqRch_{\Lambda_j,\Lambda_{j+1}} \)^{-1} \dots
    \( \qqRch_{\Lambda_j,\Lambda_n}\)^{-1}
    \ \(q^{x\bh}\)_\circlearrowleft \ \qqRch_{\Lambda_1,\Lambda_j} \dots
    \qqRch_{\Lambda_{j-1},\Lambda_j},
    \qquad j=1,\dots,n.
\end{equation*}
We introduce the qKZB operators $\qKZB_j$, acting in $\Fun \o
M_{\vec\Lambda}$, by

\begin{multline*}
    \qKZB_j \psi(\lambda) = \( \dRch_{\Lambda_j,\Lambda_{j+1}}(\lambda-\bh^{(j+2,\dots,n)}) \)^{-1} \dots
    \( \dRch_{\Lambda_j,\Lambda_n}(\lambda) \)^{-1}
    \ \Gamma_\circlearrowleft
    \\
    \dRch_{\Lambda_1,\Lambda_j}(\lambda-\bh^{(2,\dots,j-1)}-\bh^{(j+1,\dots,n)})  \dots
    \dRch_{\Lambda_{j-1},\Lambda_j}(\lambda-\bh^{(j+1,\dots,n)}) \, \psi(\lambda).
\end{multline*}

One of the fundamental properties of the universal trace operator
$\F(\lambda,x;\vec \Lambda)$ is that it intertwines the qKZ and qKZB
operators.

\begin{thm}\label{thm:holonomic system}
Let $\vec\Lambda \in \C^n[2\m]$, and let $x \in \C$. Then
\begin{equation}\label{eq:mixed KZ generic}
    \qKZB_j \  \F(\lambda,x;\vec \Lambda) =
    \F(\lambda,x;\vec \Lambda) \  \qqKZ_j(x),
    \qquad j=1,\dots,n.
\end{equation}
\end{thm}

\begin{proof}
For any $j \in \{1,\dots,n-1\}$, let $\sigma_j$ denote the
permutation of a set of $n$ elements, transposing elements in
positions $j$ and $j+1$. First, we observe that
\begin{equation}\label{eq:R Psi = Psi R}
    \dRch_{\Lambda_j,\Lambda_{j+1}}(\lambda-\bh^{(j+2,\dots,n)})
    \ \F(\lambda,x;\vec\Lambda)
    =
    \F(\lambda,x;\sigma_j(\vec\Lambda))
    \ \qqRch_{\Lambda_j,\Lambda_{j+1}}.
\end{equation}
Indeed, this is equivalent to
\begin{equation*}
    \dRch_{\Lambda_j,\Lambda_{j+1}}(\lambda-\bh^{(j+2,\dots,n)}) \
    \Tr\biggr|_{M_{\lambda-1}} \!
    \Phi^{\vec\Lambda}(\lambda) \, q^{x\bh}
    =
    \qqRch_{\Lambda_{j},\Lambda_{j+1}}^\dual \
    \Tr\biggr|_{M_{\lambda-1}} \!
    \Phi^{\sigma_j(\vec\Lambda)}(\lambda) \, q^{x\bh},
\end{equation*}
where $\dRch_{\Lambda_j,\Lambda_{j+1}}(\lambda-\bh^{(j+1,\dots,n)})$
and $\qqRch_{\Lambda_{j+1},\Lambda_j}^*$ act respectively in
$M_{\vec\Lambda}$ and $M_{\vec\Lambda}^\dual$. Moving both
$R$-matrices inside the trace and invoking Lemma \ref{thm:n=2
braiding}, we establish \eqref{eq:R Psi = Psi R}.

Let $\sigma_\circlearrowleft = \sigma_{n-1} \dots \sigma_{1}$ denote
the cyclic permutation, moving the first element to the last
position. Our second observation is that
\begin{equation}\label{eq:Gamma Psi = Psi Q}
    \Gamma_\circlearrowleft \ \F(\lambda,x;\vec\Lambda) =
    \F(\lambda,x;\sigma_\circlearrowleft(\vec\Lambda)) \ \(q^{x\bh}\)_\circlearrowleft.
\end{equation}
Indeed, using the cyclic property of trace, we compute:
\begin{equation*}
\begin{split}
    \F_{m_1,m_2,\dots,m_n}(\lambda-\Lambda_1+2m_1,x;\vec\Lambda)
    &=
    \Tr\biggr|_{M_{\lambda-\Lambda_1+2m_1-1}} \!\!\!\!\!\!\!
    (\Phi^{\Lambda_1}_{m_1}(\lambda) \o 1^{n-1})\
    \Phi^{\Lambda_2,\dots,\Lambda_n}_{m_2,\dots,m_n}(\lambda-\Lambda_1+2m_1) \, q^{x\bh}
    \\
    &=
    \mathrm P_\circlearrowleft \, \Tr\biggr|_{M_{\lambda-1}} \!\!
    \( \Phi^{\Lambda_2,\dots,\Lambda_n}_{m_2,\dots,m_n}(\lambda-\Lambda_1+2m_1)
    \,  q^{x\bh}\o 1 \)  \Phi^{\Lambda_1}_{m_1}(\lambda)
    \\
    &=
    \mathrm P_\circlearrowleft \, \Tr\biggr|_{M_{\lambda-1}} \!\!
    (q^{-x\bh})_{M_{\Lambda_1}^\dual} \ \Phi^{\Lambda_2,\dots,\Lambda_n,\Lambda_1}_{m_2,\dots,m_n,m_1}(\lambda) \,
    q^{x\bh}
    \\
    &=
    \(q^{x\bh}\)_\circlearrowleft^\dual \
    \F_{m_2,\dots,m_n,m_1}(\lambda,x;\sigma_\circlearrowleft(\vec\Lambda)),
\end{split}
\end{equation*}
which implies \eqref{eq:Gamma Psi = Psi Q}. The theorem is now
proved by repeatedly applying \eqref{eq:R Psi = Psi R} and
\eqref{eq:Gamma Psi = Psi Q}.
\end{proof}

\subsection{The qKZ equations}
For a fixed $x\in \C$ the operators $\qqKZ_j(x)$ pairwise commute,
and the problem of finding their common eigenvectors yields the qKZ
equations on $\phi \in M_{\vec\Lambda}$:
\begin{equation}\label{eq:qKZ system}
    \qqKZ_j(x) \, \phi = \eps_j \  \phi,
    \qquad\qquad j=1,\dots,n.
\end{equation}
The qKZ operators are triangular in a suitable basis, and their
eigenvalues $\eps_j$ can be described explicitly. For any $\vec m
\in \ZZ^n$ set
\begin{equation}\label{eq:qKZ eigenvalues}
    \eps_j^{\vec m}(x;\vec\Lambda) =
    q^{(\Lambda_j-2m_j)x + \sum_{i=j+1}^{n} (m_i \Lambda_j +m_j \Lambda_i - 2m_i m_j)
    - \sum_{i=1}^{j-1} (m_i \Lambda_j +m_j \Lambda_i - 2m_i m_j)}.
\end{equation}

\begin{lem}\label{thm:qKZ eigenvalues}
Let $x \in \C$ and $\vec\Lambda \in \C^n$. Suppose $\varphi \in
M_{\vec\Lambda}$ is a nonzero solution of the system \eqref{eq:qKZ
system}. Then there exists $\vec m \in \ZZ^n$ such that $\eps_j =
\eps_j^{\vec m}(x;\vec\Lambda)$ for all $j=1,\dots,n$.
\end{lem}

\begin{proof}
The qKZ operators are weight-preserving, and we may assume that
$\varphi$ is a linear combination of $v_{\vec\Lambda}^{(\vec m)}$
with $\vec m\in \ZZ^n[\m]$ for some fixed $\m \in \Z_{\ge0}$. Let
$\vec m$ be minimal in the sense of partial ordering $\preccurlyeq$
given by \eqref{eq:def partial order}, such that
$v_{\vec\Lambda}^{(\vec m)}$ occurs in the expansion of $\varphi$
with nonzero coefficient; we may assume that this coefficient is
equal to 1. Triangularity of the universal $R$-matrix implies that
for any $\vec k$
\begin{alignat*}{3}
    \qqR_{\Lambda_j,\Lambda_i} v^{(\vec k)}_{\vec \Lambda} &=
    q^{-k_i \Lambda_j -k_j \Lambda_i + 2k_i k_j} v^{(\vec k)}_{\vec \Lambda} + \text{ higher order terms},
    & \qquad i = 1,\dots,j-1,
    \\
    \qqR_{\Lambda_i,\Lambda_j}^{-1} v^{(\vec k)}_{\vec \Lambda} &=
    q^{k_i \Lambda_j +k_j \Lambda_i - 2k_i k_j} v^{(\vec k)}_{\vec \Lambda} + \text{ higher order terms},
    & \qquad i=j+1,\dots,n,
\end{alignat*}
where "higher order terms" stands for some linear combinations of
$v_{\vec\Lambda}^{(\vec l)}$ with $l \succ k$. Therefore,
\begin{equation*}
    \qqKZ_j(x) \ v^{(\vec m)}_{\vec \Lambda} =
    \eps_j^{\vec m}(x;\vec\Lambda) \  v^{(\vec m)}_{\vec \Lambda}
    + \text{ higher order terms}.
\end{equation*}
Comparing coefficients before $v^{(\vec m)}_{\vec \Lambda}$ in
\eqref{eq:qKZ system}, we obtain the desired statement.
\end{proof}

\begin{thm}\label{thm:TV KZ solutions}
Let $x \in \C$ and $\vec\Lambda \in \C^n$. Then for any $\vec m' \in
\ZZ^n$ we have
\begin{equation*}
    \qqKZ_j(x) \ \qqH^{\vec m'}(x;\vec\Lambda) =
    \eps_j^{\vec m'}(x;\vec\Lambda) \ \qqH^{\vec m'}(x;\vec\Lambda),
    \qquad j=1,\dots,n,
\end{equation*}
where $\eps_j^{\vec m}(x;\vec\Lambda)$ are as in \eqref{eq:qKZ
eigenvalues}.
\end{thm}
\begin{proof}
This statement follows from the more general results in \cite{TV
Asterisque}, see Theorem 6.2 and subsequent remarks there.
\end{proof}

\ex Let $n=2$ and $\m=1$. The qKZ equations for $j=1$ reduce to the
identity
\begin{equation*}
    \begin{pmatrix} q^{(1-x)\Lambda_1} & 0 \\
    q^{-x\Lambda_1} (1-q^{2\Lambda_2}) & q^{x(2-\Lambda_1)+\Lambda_2}  \end{pmatrix}
    \ \qqH(x; \vec\Lambda) = \qqH(x; \vec\Lambda) \
    \begin{pmatrix} q^{(1-x)\Lambda_1}& 0 \\
    0 & q^{x(2-\Lambda_1)+\Lambda_2}\end{pmatrix},
\end{equation*}
for the hypergeometric qKZ matrix $\qqH(x; \vec\Lambda)$, described
in \eqref{eq:example hypergeometric qKZ n=2,m=1}.\qed

\subsection{The qKZB equations}

If $\vec\Lambda \in \C^n[2\m]$, then the weight subspace
$M_{\vec\Lambda}[0] \ne 0$, and the qKZB operators restricted to
$M_{\vec\Lambda}[0]$ pairwise commute. The problem of finding their
common eigenfunctions yields the qKZB equations on $\varphi \in \Fun
\o M_{\vec\Lambda}[0]$:
\begin{equation}\label{eq:qKZB system}
    \qKZB_j \varphi(\lambda) = \eps_j \ \varphi(\lambda), \qquad
    j=1,\dots,n.
\end{equation}
In this paper we study solutions for \eqref{eq:qKZB system}, which
are trigonometric quasi-polynomials. The following lemma shows that
they only exist when the eigenvalues $\eps_j$ are given by the
formula \eqref{eq:qKZ eigenvalues} for the qKZ eigenvalues.

\begin{lem}\label{thm:qKZB eigenvalues}
Let $\vec\Lambda \in \C^n[2\m]$. Suppose $\psi(\lambda)$ is a
nonzero $M_{\vec\Lambda}[0]$-valued trigonometric quasi-polynomial
solution of the system \eqref{eq:qKZB system}. Then there exists $x
\in \C$ and $\vec m \in \ZZ^n[\m]$ such that $\eps_j = \eps_j^{\vec
m}(x;\vec\Lambda)$ for all $j=1,\dots,n$.
\end{lem}

\begin{proof} The idea of the proof is to observe is that the principal term of
the expansion of $\psi(\lambda)$ satisfies a version of the qKZ
equations, and use an analogue of Lemma \ref{thm:qKZ eigenvalues}.
For $y\in\C$ define the operators $\qKZ_j(y) \in
\End(M_{\vec\Lambda})$, associated with the asymptotic $R$-matrix
$\qR_{\Lambda_1,\Lambda_2}$ by
\begin{equation*}
    \qKZ_j(y) = \( \qRch_{\Lambda_j,\Lambda_{j+1}} \)^{-1} \dots
    \( \qRch_{\Lambda_j,\Lambda_n}\)^{-1}
    \ \(q^{x\bh}\)_\circlearrowleft (y) \ \qRch_{\Lambda_1,\Lambda_j} \dots
    \qRch_{\Lambda_{j-1},\Lambda_j}.
\end{equation*}
Write the trigonometric quasi-polynomial $\psi(\lambda)$ as
\begin{equation*}
    \psi(\lambda) = q^{\lambda y}
    \( \sum_{k=0}^d  q^{2 k \lambda} \, \psi^{(k)} \),
    \qquad\qquad
    \psi^{(k)}\in M_{\vec\Lambda}[0], \ \psi^{(0)}\ne 0,
\end{equation*}
for appropriate $y \in \C$ and $d \in \Z_{\ge0}$. Then the leading
coefficient $\psi^{(0)}$ satisfies the equations
\begin{equation*}
    \qKZ_j(y) \, \psi^{(0)} = \eps_j \ \psi^{(0)}, \qquad j=1,\dots,n.
\end{equation*}
The operators $\qR_{\Lambda_1,\Lambda_2}$ are triangular, and more
precisely
\begin{alignat*}{3}
    \qR_{\Lambda_j,\Lambda_i} v^{(\vec k)}_{\vec \Lambda} &=
    q^{-2 k_i(\Lambda_j-k_j)} v^{(\vec k)}_{\vec \Lambda} + \text{ higher order terms},
    & \qquad i = 1,\dots,j-1,
    \\
    \qR_{\Lambda_i,\Lambda_j}^{-1} v^{(\vec k)}_{\vec \Lambda} &=
    q^{2 k_j(\Lambda_i-k_i)} v^{(\vec k)}_{\vec \Lambda} + \text{ higher order terms},
    & \qquad i=j+1,\dots,n.
\end{alignat*}
The argument as in the proof of Lemma \ref{thm:qKZ eigenvalues}
implies that there exists $\vec m \in \ZZ^n[\m]$, such that
\begin{equation}\label{eq:asymptotic eigenvalues}
    \eps_j = q^{(\Lambda_j-2m_j)y +  2\sum_{i=j+1}^{n} m_j(\Lambda_i-m_i) -2\sum_{i=1}^{j-1} m_i(\Lambda_j-m_j)}.
\end{equation}
We now set $x = y + \m$, and from the identity
\begin{multline*}
    (\Lambda_j-2m_j)\sum_{i=1}^n m_i + 2 \sum_{i=j+1}^{n} m_j(\Lambda_i-m_i)
    - 2 \sum_{i=1}^{j-1} m_i(\Lambda_j-m_j) =
    \\
    \sum_{i=j+1}^{n} (\Lambda_i m_j + \Lambda_j m_i - 2m_i m_j)
    - \sum_{i=1}^{j-1} (\Lambda_i m_j + \Lambda_j m_i - 2m_i m_j)
    + m_j \sum_{i=1}^n (\Lambda_i - 2m_i)
\end{multline*}
and the zero weight condition we get the equivalence between
\eqref{eq:asymptotic eigenvalues} and \eqref{eq:qKZ eigenvalues}.
\end{proof}

\begin{thm}\label{thm:FTV KZB solutions}
Let $\vec\Lambda \in \C^n[2\m]$. Then for $\vec m \in \ZZ^n[\m]$ we
have
\begin{equation*}
    \qKZB_j \ \dH^{\vec m}(\lambda,x;\vec\Lambda) =
    \eps_j^{\vec m}(x;\vec\Lambda)
\ \dH^{\vec m}(\lambda,x;\vec\Lambda), \qquad j=1,\dots,n.
\end{equation*}
where $\eps_j^{\vec m}(x;\vec\Lambda)$ are as in \eqref{eq:qKZ
eigenvalues}.
\end{thm}
\begin{proof}
The statement follows from Theorem 31 in \cite{FTV monodromy}.
\end{proof}

\medskip
\begin{rem} Using our
Theorem \ref{thm:qKZ to qKZB} and Theorem \ref{thm:holonomic system},
one can immediately derive Theorem \ref{thm:FTV KZB solutions} from
the somewhat easier Theorem \ref{thm:TV KZ solutions}, derived from
the results in \cite{TV Asterisque}.\qed
\end{rem}

\medskip

\ex Let $n=2$ and $\m=1$. The qKZB equations \eqref{eq:intro qKZB
equations} reduce to the matrix identity, satisfied when
$\Lambda_1+\Lambda_2 = 2$:
\begin{multline*}
    \begin{pmatrix} q^{\Lambda_1}\frac{[\lambda+1]}{[\lambda-\Lambda_1+1]} &
    -q^{\lambda+1}\frac {[\Lambda_1]}{[\lambda-\Lambda_1+1]} \\
    q^{-\lambda+\Lambda_1+\Lambda_2-1}\frac {[\Lambda_2]}{[\lambda-\Lambda_1+1]} &
    q^{\Lambda_2}\frac {[\lambda-\Lambda_1-\Lambda_2+1]}{[\lambda-\Lambda_1+1]}  \end{pmatrix}
    \begin{pmatrix} \dH_{0,1}^{0,1}(\lambda-\Lambda_1,x; \vec\Lambda)  & \dH^{1,0}_{0,1}(\lambda-\Lambda_1,x; \vec\Lambda)
    \\ \dH_{1,0}^{0,1}(\lambda-\Lambda_1+2,x; \vec\Lambda) & \dH_{1,0}^{1,0}(\lambda-\Lambda_1+2,x; \vec\Lambda) \end{pmatrix}
    \\
    =
    \begin{pmatrix} \dH_{0,1}^{0,1}(\lambda,x; \vec\Lambda)  & \dH^{1,0}_{0,1}(\lambda,x; \vec\Lambda)
    \\ \dH_{1,0}^{0,1}(\lambda,x; \vec\Lambda) & \dH_{1,0}^{1,0}(\lambda,x; \vec\Lambda) \end{pmatrix}
    \begin{pmatrix} q^{(1-x)\Lambda_1}& 0 \\
    0 & q^{x(2-\Lambda_1)+\Lambda_2} \end{pmatrix}.
\end{multline*}

\subsection{The Macdonald-Ruijsenaars equations}

For every $\Theta \in \Z_{\ge0}$ define the Macdonald-Ruijsenaars
(MR) operators $\mathbb M_\Theta: \Fun \o M_{\vec\Lambda}[0] \to
\Fun \o M_{\vec\Lambda}[0]$  by
\begin{equation*}
    \mathbb M_\Theta \,\psi(\lambda) = q^{\m\Theta} \
    \sum_{\mu \in \h^*} \Tr\biggr|_{L_\Theta[\mu]}
    \dR_{\Theta,\Lambda_1}(\lambda-\bh^{(2,\dots,n)}) \dots
    \dR_{\Theta,\Lambda_n}(\lambda) \ \psi(\lambda-\mu).
\end{equation*}
In particular, it follows from this definition that the operator
$\mathbb M_{\Theta=0}$ is the identity operator. The operators
$\mathbb M_\Theta$ for various $\Theta \in \Z$ pairwise commute
\cite{EV}. The problem of finding their common eigenfunctions yields
the Macdonald-Ruijsenaars equations on $\psi \in \Fun \o
M_{\vec\Lambda}[0]$:
\begin{equation}\label{eq:general MR equations}
    \mathbb M_\Theta \psi(\lambda) = \mathcal X_\Theta \ \psi(\lambda),
    \qquad\qquad \Theta\in \Z_{\ge0}, \mathcal X_\Theta \in \C.
\end{equation}

For any $\Theta \in \Z_{\ge0}$ define the character $\chi_\Theta(x)$
of the irreducible $\U$-module $L_\Theta$ by
\begin{equation*}
    \chi_\Theta(x) = \sum_{\mu\in\h^*} \dim L_\Theta[\mu] \, q^{-\mu x}.
\end{equation*}

\begin{thm}\label{thm:hypergeometric Ruijsenaars}
Let $x \in \C$, and let $\vec\Lambda\in\C^n[2\m]$. Then for every
$\Theta \in \Z_{\ge0}$ one has
\begin{equation}\label{eq:MR for hypergeometric}
    \mathbb M_\Theta \ \dH(\lambda,x;\vec\Lambda) =
    \chi_\Theta(x) \, \dH(\lambda,x;\vec\Lambda).
\end{equation}
\end{thm}
\begin{proof}
It suffices to prove that in the trace function convergence domain
$|q^{2x}| \gg 1$ we have
\begin{equation}\label{eq:MR for trace}
    \mathbb M_\Theta \, \F(\lambda,x;\vec\Lambda)  =  \chi_\Theta(x) \,
    \F(\lambda,x;\vec\Lambda).
\end{equation}
Indeed, Theorem \ref{thm:qKZ to qKZB} then shows that \eqref{eq:MR
for hypergeometric} holds when $|q^{2x}| \gg 1$, and since
$\dH(\lambda,x;\vec\Lambda)$ is a trigonometric quasi-polynomial in
$x$, the desired equations are valid for all $x \in \C$.

\medskip
Consider the operators $\varPhi_{\Lambda}(\lambda): M_{\lambda-1} \o
M_{\Lambda} \to \bigoplus_\mu M_{\lambda-\mu-1} \o
M_{\Lambda}[\mu]$, which are dualized versions of the operators
\eqref{eq:formal intertwiner}, obtained by composing
$\Phi^\Lambda(\lambda)$ with the evaluation pairing $M_\Lambda^\dual
\o M_\Lambda \to \C$. If $\Theta \in \Z_{\ge0}$, the admissible
operator $\varPhi_\Theta(\lambda): M_{\lambda-1} \o L_{\Theta} \to
 \bigoplus_\mu M_{\lambda-\mu-1} \o L_{\Theta}$ is defined in the same
way. We form the following commutative diagram:
\begin{equation*}
\begin{gathered}\xymatrix@C=100pt{
    M_{\lambda-1} \o L_{\Theta} \o M_{\vec\Lambda}[0]
    \ar[r]^{\varPhi_\Theta(\lambda)}
    \ar[d]_{q^{x \bh} \o q^{-x \bh} \o \operatorname{inclusion}}
    &
    M_{\lambda-\bh^{(0)}-1} \o L_{\Theta} \o M_{\vec\Lambda}[0]
    \ar[d]^{q^{x \bh} \o 1 \o \operatorname{inclusion}}
    \\
    M_{\lambda-1} \o L_{\Theta} \o M_{\vec\Lambda}
    \ar[r]^{\varPhi_\Theta(\lambda)}
    \ar[d]_{\varPhi_{\Lambda_n}(\lambda) \, \qqR_{\Theta,\Lambda_n} }
    &
    M_{\lambda-\bh^{(0)}-1} \o L_{\Theta} \o M_{\vec\Lambda}
    \ar[d]^{\dR_{\Theta,\Lambda_n}(\lambda)\, \varPhi_{\Lambda_n}(\lambda-\bh^{(0)})}
    \\
    M_{\lambda-\bh^{(n)}-1} \o L_{\Theta} \o M_{\vec\Lambda}
    \ar[r]^{\varPhi_\Theta(\lambda-\bh^{(n)})}
    \ar@{.>}[d]
    &
    M_{\lambda-\bh^{(n,0)}-1} \o L_{\Theta} \o M_{\vec\Lambda}
    \ar@{.>}[d]
    \\
    M_{\lambda-\bh^{(2,\dots,n)}-1} \o L_{\Theta} \o M_{\vec\Lambda}
    \ar[r]_{\varPhi_\Theta(\lambda-\bh^{(2,\dots,n)})}
    \ar[d]_{\varPhi_{\Lambda_1}(\lambda-\bh^{(2,\dots,n)})\, \qqR_{\Theta,\Lambda_1}}
    &
    M_{\lambda-\bh^{(2,\dots,n,0)}-1} \o L_{\Theta} \o M_{\vec\Lambda}
    \ar[d]^{\dR_{\Theta,\Lambda_1}(\lambda) \,
    \varPhi_{\Lambda_1}(\lambda-\bh^{(2,\dots,n,0)})}
    \\
    M_{\lambda-\bh^{(1,\dots,n)}-1} \o L_{\Theta} \o M_{\vec\Lambda}
    \ar[r]_{\varPhi_\Theta(\lambda-\bh^{(1,\dots,n)})}
    \ar[d]_{\operatorname{projection}}
    &
    M_{\lambda-\bh^{(1,\dots,n,0)}-1} \o L_{\Theta} \o M_{\vec\Lambda}
    \ar[d]^{\operatorname{projection}}
    \\
    M_{\lambda-1} \o L_{\Theta} \o M_{\vec\Lambda}[0]
    \ar[r]_{\varPhi_\Theta(\lambda)}
    &
    M_{\lambda-\bh^{(0)}-1} \o L_{\Theta} \o M_{\vec\Lambda}[0]
}\end{gathered},
\end{equation*}
where we used the standard dynamical notation $\bh^{(i,\dots,j)}$
with $L_\Theta$ labeled by $0$, and the summation over weights is
implicit, for example $M_{\lambda-\bh^{(0)}-1} \o L_{\Theta}$ means
$\bigoplus_\mu M_{\lambda-\mu-1} \o L_{\Theta}[\mu]$.

\smallskip

Commutativity of the top square of the diagram is just a
reformulation of the weight-preserving property of the intertwining
operators $\Phi^\Theta(\lambda)$. Similarly, the squares in the
middle commute due to the braiding relation \eqref{eq:n=2 braiding}.
Commutativity of the bottom square is obvious.

\smallskip

Let $\mathcal A$ and $\mathcal B$ denote respectively the
compositions of operators in the left and right vertical columns in
the above diagram, multiplied by the scalar $q^{-\frac{\m (\m+1)} 2}
(q-q^{-1})^{\m} (q^x-q^{-x})$, which appear in the definition of the
trace function $\F(\lambda,x;\vec\Lambda)$. We observe that

\begin{align*}
    \Tr\biggr|_{M_{\lambda-1}\o L_\Theta} \mathcal A &=
    \F(\lambda,x;\vec\Lambda) \  \(\Tr\biggr|_{L_{\Theta}}
    \qqR_{\Theta,\Lambda_1} \dots \qqR_{\Theta,\Lambda_n} q^{-x\bh^{(0)}}
    \),
    \\
    \Tr\biggr|_{M_{\lambda-\bh^{(0)}-1}\o L_\Theta} \mathcal B &=
    \ q^{-\m\Theta} \ \mathbb M_\Theta \, \F(\lambda,x;\vec\Lambda).
\end{align*}

On the other hand, it is known (see \cite{EV,STV}) that for generic
$\lambda$ the operator $\varPhi_\Theta(\lambda)$ is a linear
isomorphism, and therefore the two traces above must be the same.
Hence
\begin{equation*}
    \mathbb M_\Theta \, \F(\lambda,x;\vec\Lambda) = q^{\m\Theta} \
    \F(\lambda,x;\vec\Lambda) \ \(\Tr\biggr|_{L_{\Theta}}
    \qqR_{\Theta,\Lambda_1} \dots \qqR_{\Theta,\Lambda_n} q^{-x\bh}
    \).
\end{equation*}
Finally, using the triangularity of the quantum $R$-matrix, one
easily shows that
\begin{equation*}
    q^{\m\Theta} \
    \Tr\biggr|_{L_{\Theta}} \, \qqR_{\Theta,\Lambda_1} \dots \qqR_{\Theta,\Lambda_n} q^{-x\bh}
    = \chi_\Theta(x) \ \Id_{M_{\vec\Lambda}},
\end{equation*}
which establishes \eqref{eq:MR for trace}, and thus completes the
proof of the theorem.
\end{proof}

\ex Let $n=2, \m = 1$, and let $\Theta = 1$. Then we have
\begin{equation*}
    \mathbb M_1 =
    \begin{pmatrix}
    \frac{[\lambda-\Lambda_2+1]}{[\lambda-\Lambda_2+2]} &
    -\frac{[\Lambda_1]}{[\lambda] \, [\lambda-\Lambda_2+2]}
    \\
    0 & \frac{[\lambda-1]}{[\lambda]} &
    \end{pmatrix}\
    \mathbb T_{+1} +
    \begin{pmatrix}
    \frac{[\lambda+1]}{[\lambda]} & 0
    \\
    -\frac{[\Lambda_2]}{[\lambda] \, [\lambda-\Lambda_2]}
    & \frac{[\lambda-\Lambda_2+1]}{[\lambda-\Lambda_2]}
    \end{pmatrix}\
    \mathbb T_{-1},
\end{equation*}
and one can check that the hypergeometric qKZB matrix
\eqref{eq:example hypergeometric qKZB n=2,m=1} satisfies the
equation
\begin{equation*}
    \mathbb M_1 \ \dH(\lambda,x;\vec\Lambda)
     = (q^x+q^{-x}) \, \dH(\lambda,x;\vec\Lambda).
\end{equation*}

\bigskip

\begin{rem}
The operators $\mathbb M_{\Theta}$ for $\Theta = 2,3\dots$ can be
expressed in terms of $\mathbb M_1$, since for any $\Theta$ we have
$\mathbb M_{\Theta} \mathbb M_1 = \mathbb M_{\Theta+1} + \mathbb
M_{\Theta-1}$. More precisely, let $p_n(t)$ be the $n$-th Chebyshev
polynomial of the second kind, so that
\begin{equation*}
    p_n(\cos \a) = \frac{\sin(n+1)\a}{\sin\a}.
\end{equation*}
Then $\mathbb M_\Theta = p_\Theta(\mathbb M_1/2)$, for example
$\mathbb M_2 = (\mathbb M_1)^2 - 1, \ \mathbb M_3 = (\mathbb M_1)^3
- 2 \, \mathbb M_1$, etc. \qed
\end{rem}

\subsection{Completeness of the hypergeometric solutions}

In this subsection we consider the qKZB and MR equations in the
space of the formal (without any convergence assumptions)
expressions of the form
\begin{equation}\label{eq:funBA series}
    \psi(\lambda) = q^{\lambda (x-\m)} \( \sum_{j=0}^\infty  \psi^{(j)} \, q^{2 j \lambda} \),
    \qquad
    \psi^{(j)} \in M_{\vec\Lambda}[0], \ \psi^{(0)} \ne 0.
\end{equation}
The qKZB and MR operators admit similar power series expansions, and
therefore act in the above space. We refer to their eigenfunctions
as formal solutions of the corresponding equations. We show that
generically all formal solutions are in fact trigonometric
quasi-polynomials, and come from the hypergeometric construction.

\begin{thm}\label{thm:completeness of hypergeometric qKZB}
Let $x \in \C$ and $\vec\Lambda\in\C^n[2\m]$ be generic. Suppose
that $\psi(\lambda)$ as in \eqref{eq:funBA series} is a formal
solution of the qKZB equations. Then there exists $\vec m \in
\ZZ^n[\m]$ such that the corresponding eigenvalues are equal to
$\eps_j^{\vec m}(x;\vec\Lambda)$, and $\psi(\lambda)$ is
proportional to $\dH^{\vec m}(\lambda,x;\vec\Lambda)$.

In particular, $\psi^{(j)} = 0$ for $j>\m$, and $\psi^{(\m)} \ne 0$.
\end{thm}
\begin{proof}
Let $\vec m$ be minimal such that $v_{\vec\Lambda}^{(\vec m)}$
occurs in the expansion of $\psi^{(0)}$ with nonzero coefficient.
Arguing as in the proof of Lemma \ref{thm:qKZB eigenvalues}, we see
that the qKZB eigenvalues are given by $\eps_j = \eps_j^{\vec
m}(x;\vec\Lambda)$. Subtracting from $\psi(\lambda)$ a suitable
multiple of $\dH^{\vec m}(\lambda,x;\vec\Lambda)$, we get a function
$\varphi(\lambda)$ satisfying
\begin{equation*}
    \qKZB_j\varphi(\lambda) = \eps_j^{\vec m}(x;\vec\Lambda) \,
    \varphi(\lambda),
\end{equation*}
such that the expansion of $\varphi(\lambda)$ does not contain the
term $q^{\lambda(x-\m)} \, v_{\vec\Lambda}^{(\vec m)}$; we claim
that $\varphi(\lambda) \equiv 0$. Indeed, otherwise the expansion of
$\varphi(\lambda)$ would begin with a nonzero multiple of
$q^{\lambda(x-\m+2k)} v_{\vec\Lambda}^{(\vec m')}$ for some $\vec m'
\in \ZZ^n[\m]$ and $k \in \Z_{\ge0}$, such that either $k>0$ or
$\vec m' \succ \vec m$. We would then get
\begin{equation}\label{eq:KZB simple spectrum}
    \eps_j^{\vec m}(x;\vec\Lambda) =
    \eps_j^{\vec m'}(x+2k;\vec\Lambda),
    \qquad j=1,\dots,n,
\end{equation}
which is impossible for generic $x$. Thus $\varphi(\lambda) \equiv
0$, and $\psi(\lambda)$ is a multiple of $\dH^{\vec
m}(\lambda,x;\vec\Lambda)$.
\end{proof}

\medskip

\begin{rem}\label{rem:simple spectrum}
The key fact used in the proof is the simplicity of the spectrum of
the qKZB system (i.e. that \eqref{eq:KZB simple spectrum} holds only
when $k=0, \vec m' = \vec m$), which for special $\vec\Lambda \in
\C^n[2\m]$ needs not be true. Nevertheless, the argument as above
shows that if \eqref{eq:KZB simple spectrum} has a {\it finite}
number of solutions, then $\psi(\lambda)$ is a finite linear
combination of corresponding $\dH^{\vec
m'}(\lambda,x+2k;\vec\Lambda)$, and in particular is a trigonometric
quasi-polynomial. The only case when the spectrum becomes infinitely
degenerate is when $\vec\Lambda = 2\vec m$, in which case
$\eps_j^{\vec m}(x;\vec\Lambda)$ are the same for all $x$, and there
exist {\it infinite} formal power series solutions, for example
$\psi(\lambda) = \sum_{k=0}^\infty \dH^{\vec
m}(\lambda,x+2k;\vec\Lambda)$.

\end{rem}

\medskip

\begin{thm}\label{thm:completeness of hypergeometric MR}
Let $x \in \C$ and $\vec\Lambda\in\C^n[2\m]$ be generic. Suppose
that $\psi(\lambda)$ as in \eqref{eq:funBA series} is a solution of
the MR equations \eqref{eq:general MR equations}. Then the
corresponding eigenvalues are equal to $\chi_\Theta(x)$, and
$\psi(\lambda)$ is a linear combination of $\dH^{\vec
m}(\lambda,x;\vec\Lambda)$ with $\vec m \in \ZZ^n[\m]$.

In particular, $\psi^{(j)} = 0$ for $j>\m$, and $\psi^{(\m)} \ne 0$.
\end{thm}
\begin{proof}
Let $\mathbf M_\Theta$ denote the limit as $q^{2\lambda} \to 0$ of
the Macdonald-Ruijsenaars operator $\mathbb M_\Theta$. It is easy to
derive from the triangularity of $\qR_{\Lambda_1,\Lambda_2}$ that
\begin{equation*}
    \Tr\bigr|_{L_\Theta[\mu]}
    \qR_{\Theta,\Lambda_1} \dots \qR_{\Theta,\Lambda_n} =
    q^{-\m(\Theta+\mu)} \dim L_\Theta[\mu],
\end{equation*}
and we obtain the explicit formula for $\mathbf M_\Theta$:
\begin{equation}\label{eq:asymptotic MR operator}
    \mathbf M_\Theta  = \sum_{\mu} q^{-\m \, \mu} \dim L_\Theta[\mu] \ \mathbb
    T_{-\mu}.
\end{equation}
The leading term $q^{\lambda(x-\m)} \, \psi^{(0)}$ of the expansion
of $\psi(\lambda)$ satisfies the equations
\begin{equation*}
    \mathbf M_\Theta \ \(q^{\lambda(x-\m)} \, \psi^{(0)} \)
    = \mathcal X_\Theta \ \(q^{\lambda(x-\m)} \, \psi^{(0)}\),
\end{equation*}
and a straightforward computation using \eqref{eq:asymptotic MR
operator} shows that $\mathcal X_\Theta = \chi_\Theta(x)$.

For generic $x$ a formal solution of the MR equations is uniquely
determined by the initial term of its expansion, and therefore the
dimension of the solution space is at most $\dim M_{\vec\Lambda}[0]
= \#\ZZ^n[\m]$, cf. Lemma 5.4 in \cite{EV}. On the other hand,
Theorem \ref{thm:hypergeometric Ruijsenaars} shows that $\dH^{\vec
m}(\lambda,x;\vec\Lambda)$ satisfy the desired MR equations.
Comparing the dimensions, we see that the linearly independent
functions $\{\dH^{\vec m}(\lambda,x;\vec\Lambda)\}_{\vec m \in
\ZZ^n[\m]}$ form a basis of the space of formal solutions.
\end{proof}

\section{The harmonic space}
\label{sec:harmonic}

\subsection{The harmonic space: generic highest weights}

In this section we assume that $\vec\Lambda \in \C^n[2\m]$ is
generic. For each $x \in \C$ we denote $\Harm_{x,\vec\Lambda}$ the
image in $\Fun \o M_{\vec\Lambda}[0]$ of the operator
$\dH(\cdot,x;\vec\Lambda)$, defined as in \eqref{eq:def
hypergeometric qKZB operator}; in other words,
\begin{equation}\label{eq:def harmonic generic}
    \Harm_{x,\vec\Lambda} =
    \Span \left\{ \dH^{\vec m'}(\cdot,x;\vec\Lambda) \ \biggr| \
    \vec m' \in \ZZ^n[\m] \right\}.
\end{equation}

It is clear that $\Harm_{x,\vec\Lambda}$ is a subspace of
$\FunBA_{x,\m} \o M_{\vec\Lambda}$. Summarizing the results of the
previous sections, we obtain

\begin{thm}\label{thm:five definitions generic}
Let $x \in \C$  be generic. Then $\Harm_{x,\vec\Lambda}$ has
dimension
\begin{equation*}
    \dim \Harm_{x,\vec\Lambda} = \dim M_{\vec\Lambda}[0] = \qbinom
    {\m+n-1}\m,
\end{equation*}
and admits the following descriptions:
\begin{enumerate}
\item \label{item:def generic trace}
$\Harm_{x,\vec\Lambda} =  \Span \left\{ \F^{\vec
m'}(\cdot,x;\vec\Lambda) \ \biggr| \
    \vec m' \in \ZZ^n[\m] \right\}$.
\item \label{item:def generic resonance}
$\Harm_{x,\vec\Lambda} =  \Span \left\{ \Psi^{\vec
m'}(\cdot,x;\vec\Lambda) \ \biggr| \
    \vec m' \in \ZZ^n[\m] \right\}$.
\item \label{item:def generic MR}
$\Harm_{x,\vec\Lambda} =  \Span \left\{ \psi \in \FunBA_{x,\m} \o
M_{\vec\Lambda}[0] \ \biggr| \
    \psi \text { is a solution of the MR equations} \right\}$.
\end{enumerate}
Moreover, if $n>1$, then we also have
\begin{enumerate}
\setcounter {enumi} 3
\item \label{item:def generic qKZB} $\Harm_{x,\vec\Lambda} =
\Span \left\{ \psi \in \FunBA_{x,\m} \o M_{\vec\Lambda}[0] \ \biggr|
\
    \psi \text { is a solution of the qKZB equations} \right\}$.
\end{enumerate}
\end{thm}

\begin{proof}
Theorem \ref{thm:qKZ to qKZB} implies \eqref{item:def generic
trace}. Theorem \ref{thm:trace vs resonance} implies \eqref{item:def
generic resonance}. Theorem \ref{thm:hypergeometric Ruijsenaars} and
Theorem \ref{thm:FTV KZB solutions} show that
$\Harm_{x,\vec\Lambda}$ is contained in the right hand sides of
\eqref{item:def generic MR} and \eqref{item:def generic qKZB}. The
opposite inclusions in \eqref{item:def generic MR} and
\eqref{item:def generic qKZB} follow from Theorem
\ref{thm:completeness of hypergeometric MR} and Theorem
\ref{thm:completeness of hypergeometric qKZB}.
\end{proof}

\subsection{The harmonic space: integral highest weights}
In this subsection we assume that $\vec\Lambda \in \ZZ^n[2\m]$. Then
the natural projection from $M_{\vec\Lambda}$ to $L_{\vec\Lambda}$
gives rise to a map from $\Fun \o M_{\vec\Lambda}[0]$ to $\Fun \o
L_{\vec\Lambda}[0]$:
\begin{equation*}
    \varphi(\lambda) = \sum_{\vec m \in \ZZ^n[\m]}
    \varphi_{\vec m}(\lambda)\  v_{\vec\Lambda}^{\vec m},
    \qquad \mapsto \qquad
    \adm\varphi(\lambda) = \sum_{\vec m \in \Adm_{\vec\Lambda}[\m]}
    \varphi_{\vec m}(\lambda) \  v_{\vec\Lambda}^{\vec m}.
\end{equation*}
For generic $x \in \C$ and $\vec m' \in \ZZ^n[\m]$ this yields
$\adm\Psi^{\vec m'}(\lambda,x;\vec\Lambda), \, \adm\F^{\vec
m'}(\lambda,x;\vec\Lambda) \in \FunBA_{x,\m} \o L_{\vec\Lambda}[0]$.
The matrix elements of the hypergeometric matrices have poles for
integral values of $\Lambda_i$. However, $\dH_{\vec m}^{\vec
m'}(\lambda,x;\vec\Lambda)$ and $\qqH_{\vec m}^{\vec
m'}(\lambda;\vec\Lambda)$ are regular at $\vec\Lambda \in
\ZZ^n[2\m]$, provided that at least one of the indices $\vec m,\vec
m'$ is $\vec\Lambda$-admissible, see \cite{MV}. Therefore, vectors
$\adm\qqH^{\vec m'}(\lambda;\vec\Lambda) \in L_{\vec\Lambda}[0]$ and
functions $\adm\dH^{\vec m'}(\lambda,x;\vec\Lambda) \in \Fun \o
L_{\vec\Lambda}[0]$ are well-defined for any $\vec m' \in
\ZZ^n[\m]$.

Let $\adm\qqH(x;\vec\Lambda) = \{\qqH_{\vec m}^{\vec
m'}(x;\vec\Lambda) \}_{\vec m,\vec m' \in \Adm_{\vec\Lambda}[\m]}$
denote the submatrix of the hypergeometric qKZ matrix, corresponding
to the $\vec\Lambda$-admissible indices, and similarly for
$\adm\dH(\lambda,x;\vec\Lambda), \adm\Psi(\lambda,x;\vec\Lambda)$
and $\adm\F(\lambda,x;\vec\Lambda)$. For any fixed $x \in \C$ these
matrices can be regarded as operators
\begin{equation*}
    \adm\qqH (x;\vec\Lambda): L_{\vec\Lambda}[0] \to
    L_{\vec\Lambda}[0],
    \qquad\qquad
    \adm\dH(\lambda,x;\vec\Lambda),
    \adm\Psi(\lambda,x;\vec\Lambda),
    \adm\F(\lambda,x;\vec\Lambda):
    L_{\vec\Lambda}[0] \to \Fun \o
    L_{\vec\Lambda}[0].
\end{equation*}

\medskip

For each $x \in \C$ we denote $\adm\Harm_{x,\vec\Lambda}$ the image
in $\Fun \o L_{\vec\Lambda}[0]$ of the operator
$\adm\dH(\cdot,x;\vec\Lambda)$, i.e.
\begin{equation}\label{eq:def harmonic integral}
    \adm\Harm_{x,\vec\Lambda} =
    \Span \left\{ \adm\dH^{\vec m'}(\cdot,x;\vec\Lambda) \ \biggr| \
    \vec m' \in \ZZ^n[\m] \right\}.
\end{equation}

\begin{thm}\label{thm:five definitions integral}
Let $x \in \C$  be generic. Then $\adm\Harm_{x,\vec\Lambda}$ admits
the following descriptions:
\begin{enumerate}
\item \label{item:def integral trace}
$\adm\Harm_{x,\vec\Lambda} =  \Span \left\{ \adm\F^{\vec
m'}(\cdot,x;\vec\Lambda) \ \biggr| \
    \vec m' \in \Adm_{\vec\Lambda}[\m] \right\}$.
\item \label{item:def integral resonance}
$\adm\Harm_{x,\vec\Lambda} =  \Span \left\{ \adm\Psi^{\vec
m'}(\cdot,x;\vec\Lambda) \ \biggr| \
    \vec m' \in \Adm_{\vec\Lambda}[\m] \right\}$.
\item \label{item:def integral MR}
$\adm\Harm_{x,\vec\Lambda} =  \Span \left\{ \psi \in \FunBA_{x,\m}
\o L_{\vec\Lambda}[0] \ \biggr| \
    \psi \text { is a solution of the MR equations} \right\}$.
\end{enumerate}
Moreover, if $n>1$ and at least one of $\Lambda_i$ is odd, then we
also have
\begin{enumerate}
\setcounter {enumi} 3
\item \label{item:def integral qKZB}
$\adm\Harm_{x,\vec\Lambda} =  \Span \left\{ \psi \in \FunBA_{x,\m}
\o L_{\vec\Lambda}[0] \ \biggr| \
    \psi \text { is a solution of the qKZB equations} \right\}$.
\end{enumerate}
The dimension of the space $\adm\Harm_{x,\vec\Lambda}$ is equal to
the cardinality of $\Adm_{\vec\Lambda}[\m]$.
\end{thm}

\begin{proof}
The argument is an obvious "admissible" modification of the proof of
Theorem \ref{thm:five definitions generic}. The assumption that one
of $\Lambda_i$ is odd guarantees that $\vec\Lambda \ne 2\vec m$, see
Remark \ref{rem:simple spectrum}.
\end{proof}

\begin{rem}
It is plausible that even for special values of $x$ the space
$\Harm_{x,\vec\Lambda}$ contains all trigonometric quasi-polynomial
solutions of the corresponding qKZB and MR equations.
\end{rem}

\subsection{The Weyl reflection}

We keep the assumption $\vec\Lambda \in \ZZ^n[2\m]$. Define the
linear map
\begin{equation*}
    \mathrm S: L_{\vec\Lambda}[0]
    \to L_{\vec\Lambda}[0],
    \qquad\qquad
    v_{\vec\Lambda}^{(\vec m)} \mapsto
    v_{\vec\Lambda}^{(\vec\Lambda-\vec m)},
\end{equation*}
and let $\mathbb S$ denote the composition of $\mathrm S$ with the
involution $\varphi(\lambda) \mapsto \varphi(-\lambda)$ of the space
$\Fun$:
\begin{equation*}
    \mathbb S: \Fun \o L_{\vec\Lambda}[0] \to \Fun \o
    L_{\vec\Lambda}[0],
    \qquad\quad
    \xi(\lambda) \o v_{\vec\Lambda}^{(\vec m)} \mapsto
    \xi(-\lambda) \o v_{\vec\Lambda}^{(\vec\Lambda-\vec m)}.
\end{equation*}

The involution $\mathbb S$ is called the Weyl reflection, and
defines an action of the Weyl group in the space of
$L_{\vec\Lambda}[0]$-valued functions.  The following lemma is
similar to Theorem 42 in \cite{FV}.

\begin{lem}\label{thm:Weyl commutes with MR and qKZB}
The Weyl reflection $\mathbb S$ commutes with the qKZB operators
$\qKZB_j$ and the Macdonald-Ruijsenaars operators $\mathbb
M_\Theta$.
\end{lem}
\begin{proof}
 It is easy to prove by induction (see also
Theorem 41 in \cite{FV}) that for any $\Lambda_1,\Lambda_2$
\begin{equation*}
    \dR_{\Lambda_1,\Lambda_2}(-\lambda) \, (s_{\Lambda_1} \o s_{\Lambda_2}) =
    (s_{\Lambda_1} \o s_{\Lambda_2}) \, \dR_{\Lambda_1,\Lambda_2}(\lambda),
\end{equation*}
where $s_{\Lambda} \in \End(L_\Lambda)$ is defined by $s_\Lambda
v_\Lambda^{(m)} = v_\Lambda^{(\Lambda-m)}$. Since $\mathrm S =
s_{\Lambda_1}  \o \dots \o s_{\Lambda_n}$, the desired statement
easily follows from the definitions and the relation $\Gamma_j \,
\mathbb S = \mathbb S \, \Gamma_j$.
\end{proof}

\medskip

\begin{thm}\label{thm:Weyl on hypergeometric}
For any $x \in \C$ we have
\begin{equation*}\label{eq:commutation Weyl with hypergeometric}
    \mathbb S \, \adm\dH(\lambda,x;\vec\Lambda) =
    \adm\dH(\lambda,-x;\vec\Lambda) \, \mathrm S.
\end{equation*}
\end{thm}

\begin{proof} We need to show that for any $\vec m \in
\Adm_{\vec\Lambda}$ the following holds:
\begin{equation}\label{eq:Weyl on hypergeometric}
    \mathbb S \, \adm\dH^{\vec m}(\lambda,x;\vec\Lambda) =
    \adm\dH^{\vec\Lambda-\vec m}(\lambda,-x;\vec\Lambda).
\end{equation}
Assume that $x \in \C$ is generic. It is clear that $\mathbb S \(
\FunBA_{x,\m} \o L_{\vec\Lambda}[0] \) \subset \FunBA_{-x,\m} \o
L_{\vec\Lambda}[0]$. Since the subspace $\adm\Harm_{x,\vec\Lambda}
\subset \FunBA_{x,\m}$ is characterized as the joint eigenspace for
the MR operators, Lemma \ref{thm:Weyl commutes with MR and qKZB}
implies that $\mathbb S \, \adm\Harm_{x,\vec\Lambda} \subset
\adm\Harm_{-x,\vec\Lambda}$. Therefore, $\mathbb S \, \adm\dH^{\vec
m}(\lambda,x;\vec\Lambda)$ can be written as a linear combination of
$\adm\dH^{\vec m'}(\lambda,-x;\vec\Lambda)$ with $\vec m' \in
\Adm_{\vec\Lambda}$. Applying Lemma \ref{thm:Weyl commutes with MR
and qKZB} again, we conclude that $\mathbb S \, \dH^{\vec
m}(\lambda,x;\vec\Lambda)$ must satisfy the qKZB equation with
eigenvalues
\begin{equation*}
    \eps_j = \eps_j^{\vec m}(x;\vec\Lambda) =
    \eps_j^{\vec\Lambda - \vec m}(-x;\vec\Lambda).
\end{equation*}
It now follows that $\mathbb S \, \adm\dH^{\vec
m}(\lambda,x;\vec\Lambda) = C \ \adm\dH^{\vec\Lambda - \vec
m}(\lambda,-x;\vec\Lambda)$ for some scalar $C$, and comparing the
principal terms of both sides, we get $C=1$. We conclude that
\eqref{eq:Weyl on hypergeometric} is valid for generic $x \in \C$,
and since both sides are holomorphic in $x$, the desired equality
holds identically.
\end{proof}

\subsection{The generalized Weyl character formula}

\begin{lem}
For $\vec m\in \Adm_{\vec\Lambda}$ we have $\Phi_{\vec
m}^{\vec\Lambda}(\lambda) \bigr(M_{\lambda-1}\bigr) \subset
M_{\lambda-1-\sum_{i=1}^n(\Lambda_i-2m_i)} \o
L_{\vec\Lambda}^\dual$.
\end{lem}
\begin{proof} The explicit definition \eqref{eq:n=1 operator}
shows that the image of $\Phi_m^\Lambda(\lambda)$ lies in the
submodule $M_{\lambda-\Lambda+2m-1} \o L_\Lambda^\dual$, settling
the case $n=1$. Obvious induction extens it for $n>1$.
\end{proof}

\begin{lem}\label{thm:admissibility check}
Let $\lambda \in \Z_{>0}$. Let $\vec\Lambda \in \ZZ^n[2\m]$ and
$\vec m\in \Adm_{\vec\Lambda}[\m]$. Then $\Phi_{\vec
m}^{\vec\Lambda}(\lambda)$ induces an operator
\begin{equation}\label{eq:admissible intertwiner}
    {}' \Phi_{\vec m}^{\vec\Lambda}(\lambda) : L_{\lambda-1} \to
    L_{\lambda-1} \o L_{\Lambda_1}^\dual \o \dots \o
    L_{\Lambda_n}^\dual.
\end{equation}
\end{lem}
\begin{proof}
Lemma \ref{thm:diagram1} and Lemma \ref{thm:diagram2} imply that for
$\Lambda \in \Z_{\ge0}$ and $m \in \{0,1,\dots,\Lambda\}$ we have
\begin{equation}\label{eq:tuda suda}
\begin{split}
    \Phi_m^\Lambda(\lambda) \circ \iota(\lambda) =
    (\iota(\lambda-\Lambda+2m )\o 1) \circ \Phi_{\Lambda-m}^\Lambda(-\lambda),
    \qquad &\text{ if } \lambda-\Lambda+2m\ge0,
    \\
    (\iota(-\lambda+\Lambda-2m) \o 1) \circ \Phi_m^\Lambda(\lambda) \circ \iota(\lambda) =
    \Phi_{\Lambda-m}^\Lambda(-\lambda),
    \qquad &\text{ if } \lambda-\Lambda+2m\le0.
\end{split}
\end{equation}
Therefore, we can form a commutative diagram
\begin{equation*}
\begin{gathered}\xymatrix@C=100pt{
    M_{-\lambda-1}
    \ar[r]^{\iota(\lambda)}
    \ar[d]_{\Phi_{\Lambda_n-m_n}^{\Lambda_n}(-\lambda)}
    &
    M_{\lambda-1}
    \ar[d]^{\Phi_{m_n}^{\Lambda_n}(\lambda)}
    \\
    M_{-\lambda+\Lambda_n-2m_n-1} \o L_{\Lambda_n}^\dual
    \ar@{<->}[r]^{\iota(|\lambda-\Lambda_n+2m_n|)\o 1}
    \ar@{.>}[d]
    &
    M_{\lambda-\Lambda_n+2m_n-1} \o L_{\Lambda_n}^\dual
    \ar@{.>}[d]
    \\
    M_{-\lambda-\Lambda_1+2m_1-1} \o L_{\Lambda_2}^\dual \o\dots\o L_{\Lambda_n}^\dual
    \ar@{<->}[r]_{\iota(|\lambda+\Lambda_1-2m_1|)\o 1^{n-1}}
    \ar[d]_{\Phi_{\Lambda_1-m_1}^{\Lambda_1}(-\lambda-\Lambda_1+2m_1) \o 1^{n-1}}
    &
    M_{\lambda+\Lambda_1-2m_1-1} \o L_{\Lambda_2}^\dual \o\dots\o L_{\Lambda_n}^\dual
    \ar[d]^{\Phi_{m_1}^{\Lambda_1}(\lambda+\Lambda_1-2m_1) \o 1^{n-1}}
    \\
    M_{-\lambda-1}\o L_{\Lambda_1}^\dual \o\dots\o L_{\Lambda_n}^\dual
    \ar[r]_{\iota(\lambda)\o 1^n}
    &
    M_{\lambda-1}\o L_{\Lambda_1}^\dual \o\dots\o L_{\Lambda_n}^\dual
}\end{gathered},
\end{equation*}
where the directions of the horizontal arrows are determined by the
positivity/negativity of the corresponding highest weights, and each
square is commutative due to \eqref{eq:tuda suda}. Thus we obtain
$\Phi_{\vec m}^{\vec\Lambda}(\lambda) M_{-\lambda-1} \subset
M_{-\lambda-1} \o L_{\vec\Lambda}^\dual[0]$, and the desired
statement follows.
\end{proof}

Let $\lambda \in \Z_{>0}$. Define the trace functions ${}'\F_{\vec
m}(\lambda,x;\vec \Lambda)$, associated with ${}' \Phi_{\vec
m}^{\vec\Lambda}(\lambda)$, by analogy with \eqref{eq:definition
trace function}:
\begin{equation*}
    {}'\F_{\vec m}(\lambda,x;\vec \Lambda) =
    q^{-\frac{\m (\m+1)} 2} (q-q^{-1})^{\m} (q^x-q^{-x}) \ \,
    \Tr\biggr|_{L_{\lambda-1}} \({}'\Phi_{\vec m}^{\vec \Lambda} (\lambda) \ q^{x \mathbf h}  \).
\end{equation*}
The coordinates $\left\{{}'\F_{\vec m}^{\vec m'}(\lambda,x;\vec
\Lambda)\right\}_{\vec m, \vec m' \in \Adm_{\vec\Lambda}[\m]}$ of
the $L_{\vec\Lambda}^\dual[0]$-valued functions ${}'\F_{\vec
m}(\lambda,x;\vec \Lambda)$ form a matrix, which can be regarded as
a linear map ${}'\F(\lambda,x;\vec \Lambda): L_{\vec\Lambda}[0] \to
\Fun \o L_{\vec\Lambda}[0]$.

\begin{thm} \label{thm:Weyl formula}
Let $\lambda \in \Z_{>0}, \m \in \Z_{\ge0}$ and $\vec\Lambda \in
\ZZ^n[2\m]$. Then for any $x \in \C$ we have
\begin{equation}\label{eq:Weyl formula}
    {}'\F(\lambda,x;\vec \Lambda)
    = \F(\lambda,x;\vec\Lambda) -
    \mathbb S \F(\lambda,x;\vec\Lambda).
\end{equation}
\end{thm}
\begin{proof}
The desired statement is equivalent to a family of relations for
$\vec m \in \Adm_{\vec\Lambda}$:
\begin{equation}\label{eq:coordinate Weyl formula}
    {}'\F_{\vec m}(\lambda,x;\vec \Lambda) =
    \F_{\vec m}(\lambda,x;\vec\Lambda) - \F_{\vec\Lambda-\vec m}(-\lambda,x;\vec\Lambda).
\end{equation}
We see from the commutative diagram in the proof of Lemma
\ref{thm:admissibility check} that
\begin{equation*}\label{eq:dynamical restriction}
    \Phi_{\vec m}^{\vec\Lambda}(\lambda) \circ \iota(\lambda) =
    (\iota(\lambda)\o1^n) \circ \Phi_{\vec\Lambda-\vec m}^{\vec\Lambda}(-\lambda),
\end{equation*}
or equivalently that the restriction of $\Phi_{\vec
m}^{\vec\Lambda}(\lambda)$ to the submodule
$\iota(\lambda)(M_{-\lambda-1}) \subset M_{\lambda-1}$ is equal to
$(\iota(\lambda)\o1^n) \circ \Phi_{\vec\Lambda-\vec
m}^{\vec\Lambda}(-\lambda)$. Therefore,
\begin{equation*}
\begin{split}
    \Tr\biggr|_{L_{\lambda-1}} \( {}'\Phi_{\vec m}^{\vec\Lambda}(\lambda) q^{x\bh} \)
    &=
    \Tr \biggr|_{M_{\lambda-1}} \( \Phi_{\vec m}^{\vec\Lambda}(\lambda) q^{x\bh} \) -
    \Tr \biggr|_{\iota(\lambda)(M_{-\lambda-1})} \( \Phi_{\vec m}^{\vec\Lambda}(\lambda) q^{x\bh} \)
    \\
    &=
    \Tr \biggr|_{M_{\lambda-1}} \( \Phi_{\vec m}^{\vec\Lambda}(\lambda) q^{x\bh} \) -
    \ \Tr \biggr|_{M_{-\lambda-1}} \( \Phi_{\vec\Lambda-\vec m}^{\vec\Lambda}(-\lambda) q^{x\bh}
    \),
\end{split}
\end{equation*}
which implies \eqref{eq:coordinate Weyl formula}, and thus proves
the theorem.
\end{proof}

Theorem \ref{thm:Weyl formula} can be regarded as a generalization
of the Weyl character formula. It also allows to extend the
definition of ${}'\F(\lambda,x;\vec\Lambda)$ from the positive
integral cone $\lambda \in \Z_{>0}$ to arbitrary values of
$\lambda$.

\section{Towards quantum conformal blocks}
\label{sec:conformal blocks}

\subsection{Fusion rules for $\mathfrak{sl}_2$ and the Grothendieck ring}
Let symbols $\binom{\lambda \ \mu}\nu$ be defined by
\begin{equation*}
    \binom{\lambda \ \ \mu}\nu =
    \begin{cases}
    1, &
    \lambda,\mu,\nu,\frac{\lambda+\mu+\nu}2 \in \Z_{\ge0} \text{ and }
    \lambda+\mu \ge \nu, \
    \lambda+\nu \ge \mu, \
    \mu+\nu \ge \lambda,
    \\
    0, & {\rm otherwise}.
    \end{cases}
\end{equation*}
These numbers arise in representation theory as the fusion rules for
finite-dimensional complex representations of $\mathfrak{sl}_2$, or
the associated quantum group $\U$ with generic $q$. Namely, for
$\lambda,\mu,\nu \in \Z_{\ge0}$ we have the equivalent definition
\begin{equation*}
    \binom{\lambda \ \ \mu}\nu =
    \text{multiplicity of $L_\nu$ in the composition series of $L_\lambda\o L_\mu$}.
\end{equation*}
We say that a triple $(\lambda,\mu,\nu)$ violates the $\sl_2$ fusion
rules, if $\binom {\lambda \ \mu}\nu = 0$.

The Grothendieck ring $\Groth$, associated with the tensor category
of finite-dimensional representations, is the ring with generators
$\{L_\lambda\}_{\lambda \in \Z_{\ge0}}$ and multiplication
\begin{equation*}
    [L_\lambda] \cdot [L_\mu] = \sum_\nu \binom {\lambda\ \ \mu}\nu
    \, [L_\nu].
\end{equation*}
The ring $\Groth$ is commutative, associative, and has the unit
element $[L_0]$.

We will need a combinatorial description of multiple products in
$\Groth$. For $\mu,\nu \in \Z$ define
$\Path_{\vec\Lambda}[\mu\rightsquigarrow\nu]$ to be the subset of
$\Adm_{\vec\Lambda}$, consisting of $\vec m$ satisfying
\begin{equation}\label{eq:def path general}
    \mu - \sum_{i=1}^n (\Lambda_i-2m_i) = \nu,
    \qquad\qquad
    \mu - \sum_{i=j}^n (\Lambda_i-2m_i) \ge m_j
    \quad \text{ for all } j=1,\dots, n,
\end{equation}
A representation-theoretic interpretation of
$\Path_{\vec\Lambda}[\mu \rightsquigarrow\nu]$ is as follows. For
$\vec m \in \ZZ^n$ denote
\begin{equation*}
    \mu_j = \mu - \sum_{i=j+1}^n (\Lambda_i-2m_i),
    \qquad j=0,\dots,n,
\end{equation*}
so that, in particular, $\mu_n = \mu$ and $\mu_0 = \nu$. Then $\vec
m \in \Path_{\vec\Lambda}[\mu \rightsquigarrow \nu]$ if and only if
the triples $(\mu_j,\Lambda_j,\mu_{j-1})$ do not violate the $\sl_2$
fusion rules for all $j=1,\dots,n$.

\begin{lem}
Let $\vec\Lambda \in \ZZ^n$. Then in the Grothendieck ring $\Groth$
we have
\begin{equation*}
    [L_{\Lambda_1}] \cdot \ldots \cdot [L_{\Lambda_n}] = \sum_{\nu}
    \# \Path_{\vec\Lambda}[\nu \rightsquigarrow 0] \  [L_\nu].
\end{equation*}
\end{lem}
\begin{proof} By induction on $n$.
\end{proof}

In particular, $\# \Path_{\vec\Lambda}[0 \rightsquigarrow 0]$ equals
the dimension of the subspace $(L_{\Lambda_1} \o \dots \o
L_{\Lambda_n})^\U$ of $\U$-invariants in the tensor product
$L_{\Lambda_1} \o \dots \o L_{\Lambda_n}$.

\medskip

\ex Let $n=4$ and $\Lambda_1 = \dots \Lambda_4 = 1$. Then
\begin{align*}
    \Path_{\vec\Lambda}[0\rightsquigarrow 0] &= \bigl\{ (0,0,1,1),(0,1,0,1) \bigr\},
    \qquad\\
    \Path_{\vec\Lambda}[2\rightsquigarrow 0] &= \bigl\{ (0,0,0,1),(0,0,1,0),(0,1,0,0) \bigr\},
    \qquad\\
    \Path_{\vec\Lambda}[4\rightsquigarrow 0]
    &= \bigl\{ (0,0,0,0) \bigr\},
\end{align*}
and in the Grothendieck ring we have
\begin{equation*}
    [L_1] \cdot [L_1] \cdot [L_1] \cdot [L_1] = 2[L_0] + 3[L_2] +
    [L_4].
\end{equation*}

\subsection{Weyl anti-symmetric functions and vanishing conditions}
In this section we assume $\vec\Lambda \in \ZZ^n[2\m]$.

Weyl anti-symmetric elements of the harmonic space satisfy the
so-called vanishing conditions, which can be elegantly formulated in
terms of the fusion rules for finite-dimensional representations.

\begin{lem} \label{thm:vanishing lemma}
Let $\delta\in\Z_{\ge0}, \vec m\in\Adm_{\vec\Lambda}[\m]$ be such
that $\vec m \notin
\Path_{\vec\Lambda}[\delta-1\rightsquigarrow\delta-1]$. Then for any
$x \in \C$ and $\vec m' \in \Adm_{\vec\Lambda}[\m]$ we have
\begin{equation}\label{eq:hypergeometric vanishing}
    \adm\dH_{\vec m}^{\vec m'}(\delta,x;\vec\Lambda)  =
    \adm\dH_{\vec\Lambda-\vec m}^{\vec m'}(-\delta,x;\vec\Lambda).
\end{equation}
\end{lem}
\begin{proof} The operator ${}'\Phi_{\vec
m}^{\vec\Lambda}(\delta)$, defined as in \eqref{eq:admissible
intertwiner}, vanishes under the given assumptions on $\delta$ and
$\vec m$. Therefore, ${}' \F_{\vec m}(\delta,x;\vec\Lambda) = 0$,
and the generalized Weyl formula \eqref{eq:coordinate Weyl formula}
yields
\begin{equation*}
    \adm\F_{\vec m}(\delta,x;\vec\Lambda) =
    \adm\F_{\vec\Lambda-\vec m}(-\delta,x;\vec\Lambda),
\end{equation*}
where the equality holds for all $x \in \C$, such that both sides
are well-defined. Theorem \ref{thm:five definitions integral} now
implies that the equation \eqref{eq:hypergeometric vanishing} is
valid for generic $x$. Since $\dH(\lambda,x;\vec\Lambda)$ is
holomorphic in $x$, we see that \eqref{eq:hypergeometric vanishing}
holds for all $x \in \C$.
\end{proof}

\begin{thm}\label{thm:vanishing conditions}
Suppose that $\psi(\lambda) \in \adm\Harm_{\vec\Lambda}$ is Weyl
anti-symmetric. Then for any $\delta \in \Z_{\ge0}$ and $\vec m \in
\Adm_{\vec\Lambda}$ such that $\vec m \notin
\Path_{\vec\Lambda}[\delta-1\rightsquigarrow\delta-1]$, we have
\begin{equation}\label{eq:vanishing conditions}
    \psi_{\vec m}(\delta) =
    \psi_{\vec\Lambda - \vec m}(-\delta) = 0.
\end{equation}
In particular, we have $\psi_{\vec m}(0) = 0$ for every $\vec m \in
\Adm_{\vec\Lambda}$.
\end{thm}
\begin{proof}
It follows from the Lemma \ref{thm:vanishing lemma} that the desired
equations \eqref{eq:vanishing conditions} are valid for functions of
the form $\psi(\lambda) = \adm\dH^{\vec m'}(\delta,x;\vec\Lambda) -
\mathbb S \adm\dH^{\vec m'}(\delta,x;\vec\Lambda)$ with $\vec m' \in
\Adm_{\vec\Lambda}$. On the other hand, such $\psi(\lambda)$ span
the subspace of Weyl anti-symmetric functions, which establishes
\eqref{eq:vanishing conditions} in the general case. Finally, for
$\delta = 0$ the inequality \eqref{eq:def path general},
corresponding to $j=n$, is violated for all $\vec m$, and we get the
last statement.
\end{proof}

Equations \eqref{eq:vanishing conditions} appeared in \cite{FV} for
the elliptic Weyl anti-symmetric hypergeometric qKZB solutions, and
are called the \textbf{\emph{vanishing conditions}}. The proof in
\cite{FV} was based on the detailed analysis of the combinatorics of
resonance relations. Our representation-theoretic argument might be
useful in the study of higher rank analogues of this phenomenon.

\medskip

\ex Let $n=4, \  \m=2$, and $\Lambda_1 = \dots  = \Lambda_4 = 1$.
Then for any Weyl anti-symmetric function $\psi(\lambda)$ the
vanishing conditions $\psi_{\vec m}(\delta) = 0$ hold  for the
following $\vec m$ and $\delta$:


$$
\begin{tabular}{c|c|c|c|c|c|c}
    $\vec m$ & (0,0,1,1) & (0,1,0,1) & (0,1,1,0) & (1,0,0,1) & (1,0,1,0) & (1,1,0,0)\\
    \hline
    $\delta$ & -2,-1,0 & -1,0 & -1,0,1 & -1,0,1 & 0,1 & 0,1,2
\end{tabular}
$$
\qed

As a consequence, we obtain the Weyl symmetry of some of the
hypergeometric qKZB solutions.

\begin{thm}\label{thm:some Weyl symmetric functions}
Let $\delta\in\Z_{\ge0}, \vec m\in\Adm_{\vec\Lambda}$ be such that
$\vec m \notin \Path_{\vec\Lambda}[\delta-1 \rightsquigarrow
\delta-1]$. Then
\begin{equation*}
    \mathbb S \dH^{\vec m}(\lambda,-\delta;\vec\Lambda) =
    \dH^{\vec m}(\lambda,-\delta;\vec\Lambda).
\end{equation*}
In particular, for any $\vec m \in \Adm_{\vec\Lambda}$ the function
$\dH^{\vec m}(\lambda,0;\vec\Lambda)$ is Weyl symmetric.
\end{thm}
\begin{proof} Observe that $\vec m \in\Path[\mu\rightsquigarrow\nu]$ if and
only if $\Opp(\vec\Lambda - \vec m) \in
\Path_{\Opp(\vec\Lambda)}[\nu \rightsquigarrow \mu]$. Using the
symmetry \eqref{eq:symmetry of hypergeometric qKZB} and Lemma
\ref{thm:vanishing lemma}, for any $\lambda \in \C$ and $\vec m' \in
\Adm_{\vec\Lambda}$ we get
\begin{equation*}
    \adm\dH^{\vec m}_{\vec m'}(\lambda,-\delta;\vec\Lambda)  =
    \adm\dH_{\Opp(\vec m)}^{\Opp(\vec m')}(-\delta,\lambda;\Opp(\vec\Lambda))  =
    \adm\dH_{\Opp(\vec\Lambda)-\Opp(\vec m)}^{\Opp(\vec m')}(\delta,\lambda;\Opp(\vec\Lambda))
    = \adm\dH^{\vec\Lambda-\vec m}_{\vec m'}(\lambda,\delta;\vec\Lambda).
\end{equation*}
Therefore, $\adm\dH^{\vec m}(\lambda,\delta;\vec\Lambda) =
\adm\dH^{\vec\Lambda-\vec m}(\lambda,-\delta;\vec\Lambda)$, and the
desired statement follows from \eqref{eq:Weyl on hypergeometric}.
\end{proof}

\subsection{The special value identity}
For any $\vec m \in \ZZ^n[\m]$ denote
\begin{equation}\label{eq:conformal spanning set}
    \cb^{\vec m}(\lambda) = \frac {\dH^{\vec m}(\lambda,-1;\vec\Lambda) -
    \dH^{\vec\Lambda - \vec
    m}(\lambda,1;\vec\Lambda)}{(q-q^{-1})^{\m+1}}.
\end{equation}

According to \eqref{eq:Weyl on hypergeometric}, the functions
$\cb^{\vec m}(\lambda)$ are Weyl anti-symmetric. The following
important result, expressing the special value $\vartheta^{\vec
m}(1)$ as a product, is somewhat similar to the Macdonald special
value identity, though the precise connection between them is
unclear.
\begin{thm} \label{thm:theta(1)}
Let $\vec\Lambda \in \ZZ^n[2\m]$ and $\vec m \in \ZZ^n[\m]$. Then
\begin{equation*}
    \cb^{\vec m}(1)  =
    -\mathbb Q_{\vec m}^{\vec \Lambda}(1;\vec\Lambda) \ v_{\vec\Lambda}^{(\vec m)}.
\end{equation*}
\end{thm}
\begin{proof}
Let the matrix elements $\{\beta_{\vec m}^{\vec m'} \}$ be defined
by
\begin{equation*}
    \Phi_{\vec m}^{\vec \Lambda}(1) v_0 = v_0 \o
    \sum_{\vec m'} \beta_{\vec m}^{\vec m'} u_{\vec m'}^{\vec\Lambda}
    + \dots,
\end{equation*}
where the omitted terms involving $v_0^{(k)}$ with $k>0$. Then we
have
\begin{equation*}
    {}'\F_{\vec m}(1,x;\vec\Lambda)
    = q^{-\frac{\m(\m+1)}2} (q-q^{-1})^\m (q^x-q^{-x}) \sum_{\vec m'} \beta_{\vec m}^{\vec m'} u_{\vec m'}^{\vec\Lambda}.
\end{equation*}
Our next goal is to compute the matrix elements $\qqH_{\vec m}^{\vec
m'}(-1;\vec\Lambda)$. We have
\begin{equation}\label{eq:qKZ(-1)}
\begin{split}
    (\Xi_{\vec\Lambda}^{\vec m} )^{-1} & \, \qqH_{\vec m}^{\vec m'}(-1;\vec\Lambda) =
    \lim_{q^{\lambda} \to 0} q^{\lambda(\m+1)} \,
    \dH_{\vec m}^{\vec m'}(\lambda,-1;\vec\Lambda) =
    \lim_{q^{\lambda} \to 0} q^{\lambda(\m+1)}
    \dH_{\vec m'}^{\vec m}(1,-\lambda;\vec\Lambda) \\
    &=
    \lim_{q^{x} \to \infty} q^{-x(\m+1)} \dH_{\vec m'}^{\vec m}(1,x;\vec\Lambda) =
    \lim_{q^{x} \to \infty} \( q^{-x(\m+1)} \sum_{\vec k} \F_{\vec m'}^{\vec k}(1,x;\vec\Lambda)
    \qqH_{\vec k}^{\vec m}(x;\vec\Lambda) \)
    \\
    &=
    \sum_{\vec k} \( \lim_{q^{x} \to \infty} q^{-x} \, \F_{\vec m'}^{\vec k}(1,x;\vec\Lambda) \)
    \(\lim_{q^{x} \to \infty} q^{-\m x} \qqH_{\vec k}^{\vec m}(x;\vec\Lambda)\).
\end{split}
\end{equation}
Observe now that $\lim_{q^{x} \to \infty} q^{-\m x}
\qqH(x;\vec\Lambda)$ is a diagonal matrix. Indeed, $\qqH_{\vec
k}^{\vec m}(x;\vec\Lambda)\equiv 0$ unless $\vec k \preccurlyeq \vec
m$ due to Proposition \ref{thm:qKZ is triangular}. If $\vec k
\preccurlyeq \vec m$, but $\vec k \ne \vec m$, then $\qqH^{\vec
k}_{\vec m}(x;\vec\Lambda)\equiv 0$, and we compute
\begin{equation*}
\begin{split}
    \lim_{q^{x} \to \infty} q^{-\m x} \qqH_{\vec k}^{\vec m}(x;\vec\Lambda) &=
    \lim_{q^{x} \to \infty} q^{-\m x} \lim_{q^\lambda \to 0}
    q^{\m\lambda} \, \dH_{\vec k}^{\vec m}(\lambda,x;\vec\Lambda) =
    \lim_{q^\lambda \to 0} q^{\m\lambda} \lim_{q^{x} \to \infty} q^{-\m x}  \,
    \dH_{\vec k}^{\vec m}(\lambda,x;\vec\Lambda) \\
    &=
    \lim_{q^{\lambda} \to 0} q^{\m\lambda} \lim_{q^{x} \to \infty} q^{-\m x} \,
    \dH^{\vec k}_{\vec m}(-x,-\lambda;\vec\Lambda)
    = \lim_{q^{\lambda} \to 0} q^{\m\lambda}  \,
    \qqH^{\vec k}_{\vec m}(-\lambda;\vec\Lambda) = 0.
\end{split}
\end{equation*}
The only remaining entries $\lim_{q^{x} \to \infty} q^{-\m x}
\qqH(x;\vec\Lambda)_{\vec m}^{\vec m}$ are easily computed using
\eqref{eq:qKZ matrix diagonal}:
\begin{equation*}
    \lim_{q^{x} \to \infty} q^{-\m x} \qqH_{\vec m}^{\vec m}(x;\vec\Lambda) =
    \frac{q^{\m(\m-1)/2}}{(q-q^{-1})^\m} \ \prod_{i=1}^n
    \frac{q^{m_i \sum_{k=i+1}^n(\Lambda_k-2m_k) - \frac {m_i(m_i-1)}2} \, [m_i]!}
    {[\Lambda_i]\dots[\Lambda_i-m_i+1]}.
\end{equation*}
It is also clear that $\lim_{q^{x} \to \infty} q^{-x} \, \F_{\vec
m'}^{\vec m}(1,x;\vec\Lambda) = q^{-\frac{\m(\m+1)}2} (q-q^{-1})^\m
\beta_{\vec m'}^{\vec m}$, and thus \eqref{eq:qKZ(-1)} yields
\begin{equation*}
    \qqH_{\vec m}^{\vec m'}(-1;\vec\Lambda) = \Xi_{\vec\Lambda}^{\vec m} \,  \prod_{i=1}^n
    \frac{q^{m_i \sum_{k=i+1}^n(\Lambda_k-2m_k) - \frac {m_i(m_i+1)}2} \, [m_i]!}
    {[\Lambda_i]\dots[\Lambda_i-m_i+1]}
    \beta_{\vec m'}^{\vec m} =
    q^{\frac{\m(\m+1)}2} \, \mathcal Q_{\vec\Lambda}^{\vec m} \, \beta_{\vec m'}^{\vec m}.
\end{equation*}
Finally, we compute using Lemma \ref{thm:orthogonality lemma}:
\begin{equation*}
\begin{split}
    \cb^{\vec m}(1) &= \frac1{(q-q^{-1})^{\m+1}}\sum_{\vec k} {}'\F^{\vec k}_{\vec m}(1,-1;\vec\Lambda) \,
    \qqH^{\vec m'}_{\vec k}(-1;\vec\Lambda) =
    -\sum_{\vec k} \beta^{\vec k}_{\vec m} \,
    \mathcal Q_{\vec k}^{\vec \Lambda} \, \beta_{\vec m'}^{\vec k}
    \\
    &=
    -\<\sum_{\vec k} \beta_{\vec m}^{\vec k} u_{\vec k}^{\vec\Lambda},
    \sum_{\vec k} \beta_{\vec m'}^{\vec k} u_{\vec k}^{\vec\Lambda}\> =
    -\<\Phi_{\vec m}^{\vec\Lambda}(1)v_0,\Phi_{\vec m'}^{\vec\Lambda}(1) v_0\> =
    -\mathbb Q_{\vec m}^{\vec \Lambda}(1;\vec\Lambda) \ v_{\vec\Lambda}^{(\vec m)}.
\end{split}
\end{equation*}

\end{proof}

\subsection{The space $\Conf_{\vec\Lambda}$}

In this section we propose a quantum analogue of the space of
conformal blocks on the torus in conformal field theory.

Let $\vec\Lambda \in \ZZ^n[2\m]$. Define the space
$\Conf_{\vec\Lambda}$ by
\begin{equation*}
    \Conf_{\vec\Lambda} = \Span \left\{  \vartheta(\lambda)\in \adm\Harm_{-1,\Vec\Lambda} +
    \adm\Harm_{1,\Vec\Lambda} \ \biggr| \ \mathbb S_{\vec\Lambda} \cb(\lambda) = - \cb(\lambda)
    \right\}.
\end{equation*}

\begin{thm}\label{thm:conformal blocks generic}
Let $\vec\Lambda \in \ZZ^n$. Then $\{\cb^{\vec m}(\lambda)\}_{\vec m
\in \Path_{\vec\Lambda}[0 \rightsquigarrow 0]}$ is a basis of
$\Conf_{\vec\Lambda}$. In particular,
\begin{equation*}
    \dim \Conf_{\vec\Lambda} = \dim (L_{\Lambda_1} \o \dots \o
    L_{\Lambda_n})^\U.
\end{equation*}
\end{thm}

\begin{proof}
It is clear from the definitions that $\{\cb^{\vec
m}(\lambda)\}_{\vec m \in \Adm_{\vec\Lambda}}$ is a spanning set for
$\Conf_{\vec\Lambda}$. Moreover, Theorem \ref{thm:some Weyl
symmetric functions} implies that $\cb^{\vec m}(\lambda) \equiv 0$
if $\vec m \notin \Path_{\vec\Lambda}[0 \rightsquigarrow 0]$.
Therefore, we only need to check that the functions $\{\cb^{\vec
m}(\lambda)\}_{\vec m \in \Path_{\vec\Lambda}[0 \rightsquigarrow
0]}$ are linearly independent.

Using the explicit formulae, one checks that $\mathbb Q_{\vec
m}^{\vec \Lambda}(1;\vec\Lambda) \ne 0$ when $\vec m \in
\Path_{\vec\Lambda}[0 \rightsquigarrow 0]$, and Theorem
\ref{thm:theta(1)} implies that $\{\cb^{\vec m}(1)\}_{\vec m \in
\Path_{\vec\Lambda}[0 \rightsquigarrow 0]}$ are linearly independent
vectors. In particular, $\{\cb^{\vec m}(\lambda)\}_{\vec m \in
\Path_{\vec\Lambda}[0 \rightsquigarrow 0]}$ are linearly independent
functions, which completes the proof.
\end{proof}
\medskip

\ex Let $n=4, \  \m=2$, and $\Lambda_1 = \dots  = \Lambda_4 = 1$.
Then the space $\Conf(\vec\Lambda)$ is spanned by two functions
$\vartheta^{(0,0,1,1)}(\lambda)$ and
$\vartheta^{(0,1,0,1)}(\lambda)$, described in the table below:
$$
\begin{tabular}{c|c|c}
    $\vec m$ & $\vartheta^{(0,0,1,1)}_{\vec m}(\lambda)$ & $\vartheta^{(0,1,0,1)}_{\vec m}(\lambda)$
    \\
    \hline
    (0,0,1,1) & $-[2][\lambda][\lambda+1][\lambda+2]$ & 0
    \\
    (0,1,0,1) & $[\lambda-1][\lambda][\lambda+1]$ & $-[\lambda][\lambda+1]^2$
    \\
    (0,1,1,0) & $[\lambda-1][\lambda][\lambda+1]$ & $[\lambda-1][\lambda][\lambda+1]$
    \\
    (1,0,0,1) & $[\lambda-1][\lambda][\lambda+1]$ & $[\lambda-1][\lambda][\lambda+1]$
    \\
    (1,0,1,0) & $[\lambda-1][\lambda][\lambda+1]$ & $-[\lambda-1]^2[\lambda]$
    \\
    (1,1,0,0) & $-[2][\lambda][\lambda-1][\lambda-2]$ & 0
    \\
\end{tabular}
$$


\medskip

In conclusion, we note that according to Theorem
\ref{thm:hypergeometric Ruijsenaars}, the conformal blocks
$\cb^{\vec m}(\lambda)$ satisfy the MR equations \eqref{eq:general
MR equations} with eigenvalues
\begin{equation*}
    \mathcal X_\Theta = \dim_q L_\Theta = \sum_\mu \dim L_\Theta[\mu] \
    q^\mu.
\end{equation*}
We conjecture that there are no other trigonometric polynomial
solutions.
\begin{conj}\label{conj:conformal generic}
Let $\vec\Lambda \in \ZZ^n[2\m]$, and let $\cb(\lambda)$ be a Weyl
anti-symmetric $L_{\vec\Lambda}[0]$-valued trigonometric polynomial,
satisfying $\mathbb M_\Theta \cb(\lambda) = ( \dim_q L_\Theta ) \
\cb(\lambda)$. Then $\cb(\lambda)\in \Conf_{\vec\Lambda}$.
\end{conj}

In other words, the space $\Conf_{\vec\Lambda}$ should be
characterized as the distinguished eigenspace of the MR operators.
Similarly, we expect that $\Conf_{\vec\Lambda}$ can be described as
the span of trigonometric polynomial solutions of the qKZB equations
of trigonometric degree at most $\m+1$.
\medskip

\subsection{Two-point quantum conformal blocks}

In this subsection we assume that $n=2$. In this case the
Macdonald-Ruijsenaars operator $\mathbb M_1$ can be written out
explicitly.

\begin{lem}\label{thm:explicit MR for n=2}
The Macdonald-Ruijsenaars operator $\mathbb M_1$ is given by
\begin{equation}
    \mathbb M_1 \psi(\lambda) = A_+(\lambda) \psi(\lambda+1) +
    A_-(\lambda) \psi(\lambda-1),
\end{equation}
where $A_\pm(\lambda) \in \End(L_{\vec\Lambda}[0])$ are represented
by matrices with matrix elements
\begin{align*}
    \( A_+(\lambda) \)_{m_1,m_2}^{m_1,m_2} &=
    \frac{[\lambda-m_1][\lambda-\Lambda_2+m_2]}{[\lambda][\lambda-\Lambda_2+2m_2]},
    &
    \( A_+(\lambda) \)_{m_1,m_2}^{m_1+1,m_2-1} &=
    \frac{[m_2][m_1-\Lambda_1]}{[\lambda][\lambda-\Lambda_2+2m_2]},
    \\
    \( A_-(\lambda) \)_{m_1,m_2}^{m_1,m_2} &=
    \frac{[\lambda+m_2][\lambda+\Lambda_1-m_1]}{[\lambda][\lambda-\Lambda_2+2m_2]},
    &
    \( A_-(\lambda) \)_{m_1,m_2}^{m_1-1,m_2+1} &=
    \frac{[m_1][m_2-\Lambda_2]}{[\lambda][\lambda-\Lambda_2+2m_2]},
\end{align*}
and all other matrix elements of $A_\pm(\lambda)$ are zero.
\end{lem}
\begin{proof}
It is clear from the definition that the only nonzero matrix
elements of $A_+(\lambda)$ are
\begin{equation*}
\begin{split}
    \( A_+(\lambda) \)_{m_1,m_2}^{m_1,m_2} &= q^{m_1+m_2}
    \( \dR_{1,\Lambda_1}(\lambda-\Lambda+2m_2)\)_{1,m_2}^{1,m_2}
    \ \( \dR_{1,\Lambda_2}(\lambda)\)_{1,m_2}^{1,m_2} ,
    \\
    \( A_+(\lambda) \)_{m_1,m_2}^{m_1+1,m_2-1}  &= q^{m_1+m_2}
    \( \dR_{1,\Lambda_1}(\lambda-\Lambda+2m_2)\)_{1,m_1}^{0,m_1+1}
    \ \( \dR_{1,\Lambda_2}(\lambda)\)_{0,m_2}^{1,m_2-1}.
\end{split}
\end{equation*}
From the computations of the dynamical $R$-matrix
$\dR_{1,\Lambda}(\lambda)$ in \cite{FTV Bethe}, we obtain
\begin{align*}
    \( \dR_{1,\Lambda}(\lambda) \)_{0,m}^{0,m} &= \frac{q^{-m}
    [\lambda+m]}{[\lambda]},
    \qquad
    &\( \dR_{1,\Lambda}(\lambda) \)_{0,m}^{1,m-1} &= -\frac{q^{-\lambda-m} [m]}{[\lambda]},
    \\
    \( \dR_{1,\Lambda}(\lambda) \)_{1,m}^{0,m+1} &= \frac{q^{\lambda-\Lambda+m} [\Lambda-m]}{[\lambda]},
    &\( \dR_{1,\Lambda}(\lambda) \)_{1,m}^{1,m} &= \frac{q^{-\Lambda+m}
    [\lambda-\Lambda+m]}{[\lambda]},
\end{align*}
which implies the desired statement for $A_+(\lambda)$. The case of
$A_-(\lambda)$ is completely analogous.
\end{proof}

\ex Let $\m = 1$. Then the MR operator $\mathbb M_1$ is given by
\begin{equation*}
    \mathbb M_1  =
    \begin{pmatrix}
    \frac{[\lambda-\Lambda_2+1]}{[\lambda-\Lambda_2+2]} & -\frac{[\Lambda_1]}{[\lambda][\lambda-\Lambda_2+2]} \\
    0 & \frac{[\lambda-1]}{[\lambda]}
    \end{pmatrix}
    \mathbb T_1 +
    \begin{pmatrix}
    \frac{[\lambda+1]}{[\lambda]} & 0\\
    -\frac{[\Lambda_2]}{[\lambda][\lambda-\Lambda_2]} & \frac{[\lambda+\Lambda_1-1]}{[\lambda-\Lambda_2]}
    \end{pmatrix}
    \mathbb T_{-1}.
\end{equation*}

\ex Let $\m = 2$. Then the MR operator $\mathbb M_1$ is given by
\begin{equation*}
    \mathbb M_1  =
    \begin{pmatrix}
    \frac{[\lambda-\Lambda_2+2]}{[\lambda-\Lambda_2+4]} & -\frac{[2][\Lambda_1]}{[\lambda][\lambda-\Lambda_2+4]} & 0 \\
    0 & \frac{[\lambda-1][\lambda-\Lambda_2+1]}{[\lambda][\lambda-\Lambda_2+2]} & \frac{[1-\Lambda_1]}{[\lambda][\lambda-\Lambda_2+2]} \\
    0 & 0 & \frac{[\lambda-2]}{[\lambda]}
    \end{pmatrix}
    \mathbb T_1 +
    \begin{pmatrix}
    \frac{[\lambda+2]}{[\lambda]} & 0 & 0 \\
    \frac{[1-\Lambda_2]}{[\lambda][\lambda-\Lambda_2+2]} & \frac{[\lambda+1][\lambda+\Lambda_1-1]}{[\lambda][\lambda-\Lambda_2+2]} & 0 \\
    0 & -\frac{[2][\Lambda_2]}{[\lambda][\lambda-\Lambda_2]} & \frac{[\lambda+\Lambda_1-2]}{[\lambda-\Lambda_2]}
    \end{pmatrix}
    \mathbb T_{-1}.
\end{equation*}
\qed
\bigskip

The following gives an explicit description of the two-point
conformal block spaces.

\begin{thm} Let $\Lambda_1,\Lambda_2 \in \Z_{\ge0}$ be such that $\Lambda_1+\Lambda_2 = 2\m$.
\begin{enumerate}
\item
If $\Lambda_1 \ne \Lambda_2$, then $\Conf_{\vec\Lambda} = 0$.
\item
If $\Lambda_1 = \Lambda_2 = \m$, then $\Conf_{\vec\Lambda} = \C \,
\cb^{(0,\m)}$, where
\begin{equation*}
    \cb^{(0,\m)}(\lambda) = \sum_{j=0}^\m (-1)^{j+1} \, \( \prod_{k=j-\m}^{j} [\lambda - k] \) \,
    \mathbf F^j v_{\m} \o \mathbf F^{\m-j} v_{\m}.
\end{equation*}
\end{enumerate}
\end{thm}

\begin{proof}
The dimensions of the spaces $\Conf_{\vec\Lambda}$ are computed by
Theorem \ref{thm:conformal blocks generic}, and it remains to
describe it explicitly for the case $\Lambda_1 = \Lambda_2 = \m$.
The vanishing conditions of Theorem \ref{thm:vanishing conditions}
for the function $\psi(\lambda) = \cb^{(0,\m)}(\lambda)$ are given
by
\begin{equation*}
    \psi_{m_1,m_2}(\delta) = 0, \qquad\qquad
    \delta = -m_2,\dots,m_1.
\end{equation*}
Therefore, $\psi_{m_1,m_2}(\lambda)$ must be divisible by
$\prod_{k=-m_2}^{m_1} [\lambda-k]$ as Laurent polynomials in
$q^\lambda$, and comparing the degrees we see that the ratios are in
fact constants, i.e. we can write
\begin{equation*}
    \psi(\lambda) = \sum_{j = 0}^\m c_j \,
    \( \prod_{k=j-\m}^{j} [\lambda - k] \) \,
    \mathbf F^{j} v_{\m} \o  \mathbf F^{\m-j} v_{\m}
\end{equation*}
for some $c_j \in \C$. From Theorem \ref{thm:theta(1)} we get
$\psi(1) = -[\m]![\m+1]! \ v_\m^{(0)} \o v_\m^{(\m)} = -[\m+1]! \
v_\m \o \mathbf F^\m v_\m$, which gives $c_0=-1$. The
Macdonald-Ruijsenaars equation $\mathbb M_1 \psi(\lambda) =
(q+q^{-1}) \, \psi(\lambda)$ yields recurrent relations on $c_j$.
Using the explicit description of the operator $\mathbb M_1$ from
Lemma \ref{thm:explicit MR for n=2}, they are simplified to $c_j = -
c_{j-1}$ for $j=1,\dots,n$.
\end{proof}

\section{Conformal blocks at roots of unity and the Verlinde algebra}

\subsection{Representations of $\U$ at roots of unity and the Verlinde algebra}
In this section we drop the assumptions $\Im \eta > 0$ and consider
the case $\eta = \frac 1\ell$ for some positive integer $\ell\ge3$.
In other words, $q = \exp\(\frac{\pi i}{\ell} \)$ becomes a
primitive $2\ell$-th root of unity.

We define the quantum group $\U$ in the same way as in
\eqref{eq:quantum commutation} and \eqref{eq:def quantum
comultiplication}. We have the infinite-dimensional $\U$-modules
$M_{\lambda}$ and $M_{\lambda}^\dual$ defined as in $\eqref{eq:Verma
action}, \eqref{eq:dual Verma action}$, and for $\lambda
\in\Z_{\ge0}$ the finite-dimensional modules $L_{\lambda},
L_{\lambda}^\dual$ are defined as before. The modules
$L_\lambda,L_\lambda^\dual$ are irreducible if and only if $\lambda
\in \{0,1,\dots,\ell-1\}$, in which case we have $L_\lambda^{} \cong
L_\lambda^\dual$.

When $q$ is a root of unity, the tensor products $L_\lambda \o
L_\mu$ need not be completely reducible. Nevertheless, there exists
a semisimple truncated tensor category, associated with $\U$-modules
$\{L_0,\dots,L_{\ell-2}\}$, with the fusion rules given by
\begin{equation*}
    \binom{\lambda \ \ \mu}{\nu}_\ell =
    \begin{cases}
    1, &
    \lambda,\mu,\nu,\frac{\lambda+\mu+\nu}2 \in \{0,1,\dots,\ell-2\} \text{ and }
    \lambda+\mu \ge \nu, \
    \lambda+\nu \ge \mu, \
    \mu+\nu \ge \lambda,
    \\
    0, & {\rm otherwise}.
    \end{cases}
\end{equation*}
For $\lambda,\mu,\nu \in \{0,1,\dots,\ell-2\}$ we have the
equivalent definition
\begin{equation*}
    \binom{\lambda \ \ \mu}{\nu}_\ell =
    \text{multiplicity of $L_{\nu}$ as a direct summand of $L_{\lambda}\o L_{\mu}$
    }.
\end{equation*}
These numbers differ from the composition series multiplicities,
since $L_\nu$ may occur as Jordan-H\"older constituents in reducible
but indecomposable tilting modules. We say that a triple
$(\lambda,\mu,\nu)$ violates the $\ell$-truncated $\sl_2$ fusion
rules, if $\binom {\lambda \ \mu}{\nu}_\ell = 0$.

The associated Verlinde algebra $\Ver_\ell(\sl_2)$ is the algebra
with generators $\{L_\lambda\}_{\lambda \in \{0,1,\dots,\ell-2\}}$
and multiplication
\begin{equation*}
    [L_\lambda] \cdot [L_\mu] = \sum_\nu \binom {\lambda\ \
    \mu}{\nu}_\ell
    \, [L_\nu].
\end{equation*}
The algebra $\Ver_\ell(\sl_2)$ is commutative, associative, and has
the unit element $[L_0]$.

\medskip \ex Let $\ell=3$. Then $\Ver_\ell(\mathfrak{sl}_2)$ is
generated by $[L_0],[L_1]$ with multiplication
\begin{equation*}
    [L_0] \cdot [L_0] = [L_0], \qquad [L_0] \cdot [L_1] = [L_1],
    \qquad [L_1] \cdot [L_1] = [L_0].
\end{equation*}
In other words, the Verlinde algebra for $\ell=3$ is isomorphic to
the group algebra $\C[\Z_2]$.\qed

\bigskip
For $\mu,\nu \in \Z$ denote
$\Path_{\vec\Lambda}^{(\ell)}[\mu\rightsquigarrow\nu]$ the subset of
$\Path_{\vec\Lambda}[\mu\rightsquigarrow\nu]$, consisting of $\vec
m$, satisfying in addition to \eqref{eq:def path general} the extra
conditions
\begin{equation}\label{eq:def path root}
    \mu - \sum_{i=j+1}^n (\Lambda_i-2m_i) \le \ell-m_j-2
    \quad \text{ for all } j=1,\dots, n.
\end{equation}

\medskip

\begin{lem}
Let $\Lambda_1,\dots,\Lambda_n \in \{0,1,\dots,\ell-2\}$. Then in
the Verlinde algebra we have
\begin{equation*}
    [L_{\Lambda_1}] \cdot \ldots \cdot [L_{\Lambda_n}] = \sum_\nu
    \# \Path^{(\ell)}_{\vec\Lambda}[\nu \rightsquigarrow 0] \ [L_\nu].
\end{equation*}
\end{lem}
\begin{proof} By induction on $n$.
\end{proof}

\medskip

Denote by $(L_{\Lambda_1} \o \dots \o
L_{\Lambda_n})^{\Ver_\ell(\sl_2)}$ the sum of all direct summands of
$L_{\Lambda_1} \o \dots \o L_{\Lambda_n}$, isomorphic to $L_0 = \C$.
Then, in particular, $\dim (L_{\Lambda_1} \o \dots \o
L_{\Lambda_n})^{\Ver_\ell(\sl_2)} = \#
\Path_{\vec\Lambda}^{(\ell)}[0 \rightsquigarrow 0]$. The remaining
part of this section is devoted to proving an analogue of Theorem
\ref{thm:conformal blocks generic}: for any $\vec\Lambda$ we
construct a distinguished subspace of $\adm\Harm_{\vec\Lambda}$,
whose dimension is precisely equal to $\dim (L_{\Lambda_1} \o \dots
\o L_{\Lambda_n})^{\Ver_\ell(\sl_2)}$, see Theorem
\ref{thm:conformal blocks roots}.

\subsection{The regularity lemmas}
The purpose of this subsection is to show that the hypergeometric
matrices and holomorphic intertwining operators remain well-defined
when $q$ becomes a root ot unity.

\begin{lem}
Let $\vec m,\vec m' \in \ZZ^n[\m]$. For $\lambda,x \in \C$ and
$\vec\Lambda \in \C^n$ consider the expression
\begin{equation*}
    H = q^{-\lambda x + \m x + \m\lambda} \( \prod_{j=1}^n
\prod_{k=1}^{\min(m_j,m'_j)} (1 - q^{2(\Lambda_j-k)}) \) \,
\dH_{\vec m}^{\vec m'}(\lambda,x;\vec\Lambda).
\end{equation*}
Then $H$ is a Laurent polynomial in
$q^2,q^{2x},q^{2\lambda},q^{\Lambda_1},\dots,q^{\Lambda_n}$.
\end{lem}
\begin{proof}
It is easy to see that each residue term in the definition of the
hypergeometric pairing has the desired representation. We leave
details to the reader.
\end{proof}

\begin{lem}\label{thm:matrix elements regularity}
Let $\Lambda \in \Z$. Then the matrix elements of the operator
$(q-q^{-1})^m \Phi_m^\Lambda(\lambda)$ with respect to bases
$\{v_{\lambda-1}^{(k)}\}$ and $\{v_{\lambda-\Lambda+2m-1}^{(k)} \o
u_l^{\Lambda}\}$ are Laurent polynomials in $q,q^{2\lambda}$.
\end{lem}

\begin{proof}
Denote for convenience $\mu = \lambda-\Lambda+2m$. Using
\eqref{eq:n=1 operator}, we obtain
\begin{equation}\label{eq:regularity of Phi}
\begin{split}
    \Phi_m^\Lambda(\lambda) v_{\lambda-1}^{(k)} &=
    \( \sum_{j=0}^k \frac{q^{j(k-j)}}{[j]![k-j]!}
    \mathbf F^j \o \mathbf F^{k-j} q^{-j\bh} \)
    \( \sum_{i=0}^m \frac{q^{i(\mu-i)}}{[i]!} [-\mu+i+1]_{m-i}
    \mathbf F^i v_{\mu-1} \o \mathbf E^i u_m^{\Lambda} \)
    \\
    &= \sum_{l=0}^{m+k} C^k_l(\mu,\Lambda,m) \  v_{\mu-1}^{(l)} \o
    u_{m-l+k}^{\Lambda},
\end{split}
\end{equation}
where we denoted
\begin{equation}\label{eq:Phi matrix elements}
    C^k_l(\mu,\Lambda,m) = \sum_{i+j=l}
    \frac{q^{j(k-j-\Lambda+2m-2i)+i(\mu-i)} \, [l]!}{[j]!} \  \qbinom{m}{i}
    \qbinom{\Lambda-m+i}{k-j}  \, [-\mu+i+1]_{m-i} .
\end{equation}
Treating $q$ as a formal variable, we see that $(q-q^{-1})^m
[-\mu+i+1]_{m-i}$ are Laurent polynomials in $q, q^\mu$. Also, for
any $a \in \Z, b \in \Z_{>0}$ the binomial coefficient $\qbinom ab$
is a Laurent polynomial in $q$. The desired matrix elements are
linear combinations of $(q-q^{-1})^m C^k_l(\mu,\Lambda,m)$, and
therefore are well-defined for any nonzero $q\in\C$.
\end{proof}

\medskip

\subsection{Weyl symmetry and vanishing conditions at roots of unity}

For $\mu \in \{0,\dots,\ell-1\}$ there exists a proper inclusion of
$\U$-modules $\jmath(\mu) : L_{\ell-\mu-1} \hookrightarrow
L_{\ell+\mu-1}$, uniquely determined by $\jmath(\mu)
(v_{\ell-\mu-1}) = \mathbf F^{\mu} \, v_{\ell+\mu-1}$.

\begin{lem}\label{thm:diagrams roots}
Let $\mu,\Lambda \in \{0,1,\dots,\ell-1\}$, and let $m \in
\Z_{\ge0}$ be such that $m+\mu \le \Lambda$. Then we have
commutative diagrams
\begin{equation*}
\begin{gathered}\xymatrix{
    &
    L_{\ell-\mu+\Lambda-2m-1}
    \ar[rd]^(0.6){\quad\Phi_{m+\mu}^\Lambda(\ell-\mu+\Lambda-2m)}
    \ar[ld]_(0.6){\Phi_{m}^\Lambda(\ell-\mu+\Lambda-2m)\quad}
    &
    \\
    L_{\ell-\mu-1} \o L_{\Lambda}^\dual
    \ar[rr]_{\jmath(\mu)\o 1}
    &&
    L_{\ell+\mu-1} \o L_{\Lambda}^\dual
}\end{gathered}.
\end{equation*}
\begin{equation*}
\begin{gathered}\xymatrix{
    L_{\ell-\mu-1}
    \ar[rd]_(0.4){\Phi_{m+\mu}^\Lambda(\ell-\mu)\ }
    \ar[rr]^{\jmath(\mu)}
    &&
    L_{\ell+\mu-1}
    \ar[ld]^(0.4){\Phi_m^\Lambda(\ell+\mu)}
    \\
    &
    L_{\ell-\Lambda+2m+\mu-1} \o L_{\Lambda}^\dual
    &
}\end{gathered}.
\end{equation*}
\end{lem}
\begin{proof}
Similar to proofs of Lemma \ref{thm:diagram1} and Lemma
\ref{thm:diagram2}; we skip further details.
\end{proof}

\begin{lem} \label{thm:vanishing lemma roots}
Let $\delta \in \{0,\dots,\ell-1\}$, $\vec m \in
\Path_{\vec\Lambda}[\delta-1 \rightsquigarrow \delta -1]$ be such
that at least one of the inequalities \eqref{eq:def path root} is
violated. Then for any $x \in \C$ and $\vec m' \in
\Adm_{\vec\Lambda}[\m]$ we have
\begin{equation}\label{eq:hypergeometric vanishing roots}
    \adm\dH_{\vec m}^{\vec m'}(\delta,x;\vec\Lambda)  =
    \adm\dH_{\vec\Lambda-\vec m}^{\vec m'}(2\ell-\delta,x;\vec\Lambda).
\end{equation}
\end{lem}

\begin{proof} Using the operators ${}'\Phi_{\vec
m}^{\vec\Lambda}(\delta)$, defined as in \eqref{eq:admissible
intertwiner}, we form the commutative diagram
\begin{equation*}
\begin{gathered}\xymatrix@C=100pt{
    L_{\delta-1}
    \ar[r]^{\jmath(\lambda)}
    \ar[d]_{'\Phi_{m_n}^{\Lambda_n}(\lambda)}
    &
    L_{2\ell-\delta-1}
    \ar[d]^{'\Phi_{\Lambda_n-m_n}^{\Lambda_n}(2\ell-\lambda)}
    \\
    L_{\delta-\Lambda_n+2m_n-1} \o L_{\Lambda_n}^\dual
    \ar@{<->}[r]^{\jmath(|\lambda-\Lambda_n+2m_n|)\o 1}
    \ar@{.>}[d]
    &
    L_{2\ell-\delta+\Lambda_n-2m_n-1} \o L_{\Lambda_n}^\dual
    \ar@{.>}[d]
    \\
    L_{\delta+\Lambda_1-2m_1-1} \o L_{\Lambda_2}^\dual \o\dots\o L_{\Lambda_n}^\dual
    \ar@{<->}[r]_{\jmath(|\lambda+\Lambda_1-2m_1|)\o 1^{n-1}}
    \ar[d]_{'\Phi_{m_1}^{\Lambda_1}(\lambda+\Lambda_1-2m_1) \o 1^{n-1}}
    &
    L_{2\ell-\delta-\Lambda_1+2m_1-1} \o L_{\Lambda_2}^\dual \o\dots\o L_{\Lambda_n}^\dual
    \ar[d]^{'\Phi_{\Lambda_1-m_1}^{\Lambda_1}(2\ell-\lambda-\Lambda_1+2m_1) \o 1^{n-1}}
    \\
    L_{\delta-1}\o L_{\Lambda_1}^\dual \o\dots\o L_{\Lambda_n}^\dual
    \ar[r]_{\jmath(\lambda)\o 1^n}
    &
    L_{2\ell-\delta-1}\o L_{\Lambda_1}^\dual \o\dots\o L_{\Lambda_n}^\dual}
\end{gathered}.
\end{equation*}
Multiplying by $q^{x\bh}$ and taking the trace, we obtain
${}'\adm\F_{\vec m}(\delta,x;\vec\Lambda) = {}'\adm\F_{\vec\Lambda-
\vec m}(2\ell-\delta,x;\vec\Lambda)$, which can be rewritten using
the Weyl formula \eqref{eq:Weyl formula} as
\begin{equation*}
    \adm\F_{\vec m}(\delta,x;\vec\Lambda) - \adm\F_{\vec\Lambda - \vec m}(-\delta,x;\vec\Lambda)
    =
    \adm\F_{\vec\Lambda - \vec m}(2\ell-\delta,x;\vec\Lambda) - \adm\F_{\vec
    m}(\delta-2\ell,x;\vec\Lambda).
\end{equation*}
Using the quasi-periodicity $\adm\F(\lambda+2\ell,x;\vec\Lambda) =
e^{2\pi i x} \adm\F(\lambda,x;\vec\Lambda)$, we obtain
\begin{equation*}
    \adm\F_{\vec m}(\delta,x;\vec\Lambda) =
    \adm\F_{\vec\Lambda - \vec m}(2\ell-\delta,x;\vec\Lambda).
\end{equation*}
The desired equation \eqref{eq:hypergeometric vanishing roots} for
the hypergeometric qKZB matrix now follows from the admissible
version of Theorem \ref{thm:five definitions integral}.
\end{proof}

\subsection{Conformal blocks at roots of unity}

In this subsection we assume that $\vec\Lambda \in \ZZ^n[2\m]$ is
such that $\Lambda_j \in \{0,1,\dots,\ell-1\}$ for all
$j=1,\dots,n$.

To each function $\varphi: \C \to L_{\vec\Lambda}[0]$ we associate
its restriction $\tilde\varphi: \Z \to L_{\vec\Lambda}[0]$ to the
integral lattice. We refer to $\tilde \varphi$ as the sequence,
corresponding to the function $\varphi$.

Define the space $\Conf^{(\ell)}_{\vec\Lambda}$ of
\textbf{\emph{discrete conformal blocks}} to be the space of
sequences, corresponding to functions $\cb(\lambda) \in
\Conf_{\vec\Lambda}$. Due to the quasi-periodicity
$\cb(\lambda+\ell) = (-1)^{\m+1} \, \cb(\lambda)$, a discrete
conformal block is completely determined by its values at $\lambda =
0,1,\dots,\ell-1$, and thus $\Conf_{\vec\Lambda}^{(\ell)}$ can be
thought of as a subspace of $\(L_{\vec\Lambda}[0]\)^\ell$.

\begin{thm} \label{thm:conformal blocks roots}
The set $\{\tilde\cb^{\vec m}(\lambda)\}_{\vec m \in
\Path_{\vec\Lambda}^{(\ell)}[0 \rightsquigarrow 0]}$ is a basis of
$\Conf_{\vec\Lambda}^{(\ell)}$. In particular,
\begin{equation*}
    \dim \Conf^{(\ell)}_{\vec\Lambda} = \dim (L_{\Lambda_1} \o \dots \o
    L_{\Lambda_n})^{\Ver_\ell(\sl_2)}.
\end{equation*}
\end{thm}
\begin{proof}
The argument is similar to that of Theorem \ref{thm:conformal blocks
generic}. If $\vec m \notin \Path^{(\ell)}_{\vec\Lambda}[0
\rightsquigarrow 0]$, then we claim that $\tilde\cb^{\vec
m}(\lambda) \equiv 0$. Indeed, for $\vec m \notin
\Path_{\vec\Lambda}$ this follows from the argument in the proof of
Theorem \ref{thm:conformal blocks generic}. Assume now that $\vec m
\in \Path_{\vec\Lambda}$, or in other words that the inequalities
\eqref{eq:def path general} hold, but some of \eqref{eq:def path
root} fail. Using Lemma \ref{thm:vanishing lemma roots} and the
symmetry of $\dH(\lambda,x;\vec\Lambda)$, for any $\vec m'$ we
compute
\begin{equation*}
\begin{split}
    \dH^{\vec m}_{\vec m'}(\lambda,-1;\vec\Lambda) &=
    \dH_{\vec m}^{\vec m'}(1,-\lambda;\vec\Lambda) =
    \dH_{\Lambda-\vec m}^{\vec m'}(2\ell - 1,-\lambda;\vec\Lambda) =\\
    &=
    e^{-2\pi i \lambda} \, \dH_{\Lambda-\vec m}^{\vec m'}(- 1,-\lambda;\vec\Lambda) =
    e^{-2\pi i \lambda} \, \dH^{\Lambda-\vec m}_{\vec m'}(\lambda,1;\vec\Lambda),
\end{split}
\end{equation*}
and for $\lambda \in \Z$ we obtain $\dH^{\vec m}_{\vec
m'}(\lambda,-1;\vec\Lambda)= \dH^{\Lambda-\vec m}_{\vec
m'}(\lambda,1;\vec\Lambda)$, which means that $\tilde\cb^{\vec
m}(\lambda) = 0$.

It remains to verify the linear independence of $\{\tilde\cb^{\vec
m}(\lambda)\}_{\vec m \in \Path_{\vec\Lambda}^{(\ell)}[0
\rightsquigarrow 0]}$. As in the proof of Theorem \ref{thm:conformal
blocks generic}, one immediately checks from the definitions that
$\tilde\cb^{\vec m}(1) = -\mathbb Q_{\vec m}^{\vec
\Lambda}(1;\vec\Lambda) \ v_{\vec\Lambda}^{(\vec m)} \ne 0$ when
$\vec m \in \Path^{(\ell)}_{\vec\Lambda}[0 \rightsquigarrow 0]$.
This proves the desired linear independence, and concludes the proof
of the theorem.
\end{proof}

\ex Let $n=4, \  \m=2$, and $\Lambda_1 = \dots  = \Lambda_4 = 1$.
Then $\Conf_{\vec\Lambda}$ is two-dimensional, as described in the
example in Section \ref{sec:conformal blocks}. However, the values
$\tilde\vartheta^{(0,0,1,1)}(\lambda)$ for all $\lambda \in \Z$ are
easily seen to be zero when $q = e^{\pi i/3}$. Therefore,
$\Conf^{(\ell=3)}_{\vec\Lambda}$ in this case is one-dimensional,
and is spanned by the restriction of
$\tilde\cb^{(0,1,0,1)}(\lambda)$ to $\Z$.

On the other hand, in the Verlinde algebra, corresponding to
$\ell=3$, we compute
\begin{equation*}
    [L_1] \cdot [L_1] \cdot [L_1] \cdot [L_1] = [L_0],
\end{equation*}
and the space of Verlinde algebra invariants in this case also has
dimension one. \qed

\subsection{The qKZB and Macdonald-Ruijsenaars operators on the lattice}

The discrete conformal block space $\Conf_{\vec\Lambda}^{(\ell)}$
should admit a characterization as the space of solutions of the
qKZB and MR equations, restricted to the integral lattice. The
Macdonald-Ruijsenaars difference operators involve only integral
shifts in the $\lambda$ variable, and the same is true for the qKZB
operators because the highest weights $\vec\Lambda$ are integers.
Therefore, we can consider the restriction of these operators to
cosets $\C/\Z$, i.e. regard them as difference operators on the
lattice of the form $\eps + \Z$ with $\eps \in \C$.

However, in the most important case $\eps = 0$, the coefficients of
operators $\qKZB_j$ and $\mathbb M_\Theta$ have poles for small
integral values of $\lambda$, and in order to get well-defined
difference operators we introduce their modified versions
$\widetilde{\qKZB}_j$ and $\widetilde{\mathbb M}_\Theta$ by the
following procedure.

Write the Macdonald-Ruijsenaars operator as $\mathbb M_\Theta =
\sum_{\mu\in\Z} A_\mu(\lambda) \, \mathbb T_\mu$ for some
meromorphic in $\lambda$ coefficients $A_\mu(\lambda)\in
\End(L_{\vec\Lambda}[0])$ with matrix elements
$\(A_\mu(\lambda)\)_{\vec m}^{\vec m'}$, and for each $\delta \in
\Z$ define the reduced coefficients $\widetilde{A}_\mu(\lambda)\in
\End(L_{\vec\Lambda}[0])$ with matrix elements
\begin{equation*}
    (\widetilde{A}_\mu(\delta))_{\vec m}^{\vec m'} = \begin{cases}
    \(A(\delta)\)_{\vec m}^{\vec m'}, & \text {if $\(A(\lambda)\)_{\vec m}^{\vec m'}$ is well-defined at $\lambda = \delta$},\\
    0, & \text {if $\(A(\lambda)\)_{\vec m}^{\vec m'}$ has a pole at $\lambda = \delta$}.
    \end{cases}
\end{equation*}
This reduction procedure makes sense both when $q$ is generic, and
when $q$ is a root of unity; in the latter case the matrices
$\widetilde{A}_\mu(\delta)$ are $\ell$-periodic. For each
$\varphi:\Z \to L_{\vec\Lambda}[0]$ we now set
\begin{equation*}
    \widetilde{\mathbb M}_\Theta \varphi (\delta) = \sum_{\mu\in\Z}
    \widetilde{A}_\mu(\delta) \, \varphi(\delta+\mu).
\end{equation*}
The reduced coefficients of $\qKZB_j$ and the modified operator
$\widetilde{\qKZB}_j$ are defined in a similar fashion.

\medskip

The above reduction removes the singularities of the qKZB and MR
operators in a somewhat {\it ad hoc} fashion, and in general the
eigenfunctions of $\qKZB_j$ and $\mathbb M_\Theta$ need not restrict
on the lattice to eigensequences of $\widetilde{\qKZB}_j$ and
$\widetilde{\mathbb M}_\Theta$. However, the special values of
$\lambda$, affected by the reduction procedure, are precisely the
values participating in the resonance relations and the vanishing
conditions, and it seems plausible that the order of vanishing of
$\tilde\cb^{\vec m}(\lambda)$ at those points is higher than the
order of the removed poles of $\qKZB_j$ and $\mathbb M_\Theta$. The
latter property would imply that the restrictions $\tilde\cb^{\vec
m}(\lambda)$ remain the eigenfunctions of the modified operators:

\begin{conj}
The sequences $\tilde\cb^{\vec m}(\lambda)$ satisfy the equations
\begin{equation}\label{eq:reduced MR equations}
    \widetilde{\qKZB}_j \tilde\cb(\delta) = \eps_j^{\vec m}(-1;\vec\Lambda) \,
    \tilde\cb(\delta),
    \qquad\qquad
    \widetilde{\mathbb M}_\Theta \tilde\cb(\delta) = \( \dim_q
    L_\Theta \)\, \tilde\cb(\delta).
\end{equation}
\end{conj}

\medskip

The results of Section 12 in \cite{FV} imply that this conjecture
holds for the operators $\widetilde{\qKZB}_j$ in the case when
$\Lambda_1 = \dots = \Lambda_n = 1$.

\medskip

Finally, we formulate our last conjecture, which is similar to
Conjecture \ref{conj:conformal generic}.

\begin{conj}\label{conj:conformal root}
Let $q = \exp\(\frac{\pi i}\ell\)$ with $\ell \ge 3$. Let $\tilde
\cb: \Z \to
 L_{\vec\Lambda}[0]$ be a solution of the equations $\widetilde{\mathbb M}_\Theta
\tilde\cb(\delta) = ( \dim_q L_\Theta ) \, \tilde\cb(\delta)$ such
that $\tilde\cb(\delta+\ell) = (-1)^{\m+1}\tilde\cb(\delta)$. Then
$\tilde\cb(\delta) \in \Conf_{\vec\Lambda}^{(\ell)}$.
\end{conj}

We illustrate it on the simplest example.

\medskip
\ex Let $n=2, \m=1$ and $\Lambda_1 = \Lambda_2 = 1$. Suppose that a
vector-valued sequence $\varphi(\delta) = \binom
{\varphi_{0,1}(\delta)}{\varphi_{1,0}(\delta)}$, $\delta \in\Z$,
satisfies the reduced MR equations. Using the explicit formula for
the Macdonald operator in Lemma \ref{thm:explicit MR for n=2} we get
\begin{equation*}
    \mathbb M_1  =
    \begin{pmatrix}
    \frac{[\lambda]}{[\lambda+1]} & -\frac{1}{[\lambda][\lambda+1]} \\
    0 & \frac{[\lambda-1]}{[\lambda]}
    \end{pmatrix}
    \mathbb T_1 +
    \begin{pmatrix}
    \frac{[\lambda+1]}{[\lambda]} & 0\\
    -\frac{1}{[\lambda][\lambda-1]} & \frac{[\lambda]}{[\lambda-1]}
    \end{pmatrix}
    \mathbb T_{-1},
\end{equation*}
and from the definitions we see that
\begin{equation*}
\begin{split}
    \widetilde{\mathbb M}_1\varphi(0) &=
    \begin{pmatrix} 0 & 0 \\ 0 & 0 \end{pmatrix}
    \binom {\varphi_{0,1}(1)}{\varphi_{1,0}(1)} +
    \begin{pmatrix} 0 & 0 \\ 0 & 0 \end{pmatrix}
    \binom {\varphi_{0,1}(-1)}{\varphi_{1,0}(-1)}
    =
    \binom 0 0,
    \\
    \widetilde{\mathbb M}_1\varphi(1) &=
    \begin{pmatrix} \frac 1{[2]} & -\frac 1{[2]} \\ 0 & 0 \end{pmatrix}
    \binom {\varphi_{0,1}(2)}{\varphi_{1,0}(2)} +
    \begin{pmatrix} [2] & 0 \\ 0 & 0 \end{pmatrix}
    \binom {\varphi_{0,1}(0)}{\varphi_{1,0}(0)}
    =
    \binom {\frac {\varphi_{0,1}(2) - \varphi_{1,0}(2)}{[2]} +  [2] \varphi_{0,1}(0)} 0,
    \\
    \widetilde{\mathbb M}_1\varphi(-1) &=
    \begin{pmatrix} 0 & 0 \\ 0 & [2] \end{pmatrix}
    \binom {\varphi_{0,1}(0)}{\varphi_{1,0}(0)} +
    \begin{pmatrix} 0 & 0 \\ -\frac 1{[2]} & \frac 1{[2]} \end{pmatrix}
    \binom {\varphi_{0,1}(-2)}{\varphi_{1,0}(-2)}
    =
    \binom 0 {[2] \varphi_{1,0}(0) - \frac{\varphi_{0,1}(-2) -
    \varphi_{1,0}(-2)}{[2]}}.
\end{split}
\end{equation*}

The equation $\widetilde{\mathbb M}_1 \varphi(\delta) = [2] \,
\varphi(\delta)$ now implies the vanishing conditions
\begin{equation*}
    \varphi_{0,1}(0) = \varphi_{1,0}(0) = \varphi_{0,1}(-1) =
    \varphi_{1,0}(1) = 0.
\end{equation*}
Let $\ell = 4$. The quasi-periodicity shows that our sequence must
have the form
$$
\begin{tabular}{c|c|c|c|c|c|c|c|c|c|c|c|c|c}
    $\delta$ & \dots & -5 & -4 & -3 & -2 & -1 & 0 & 1 & 2 & 3 & 4 & 5 & \dots\\
    \hline
    $\varphi_{0,1}(\delta)$ & \dots & 0 & 0 & $\a_1$ & $\a_2$ & 0 & 0 & $\a_1$ & $\a_2$ &  0 & 0 & $\a_1$ & \dots \\
    \hline
    $\varphi_{1,0}(\delta)$ & \dots & $\beta_2$ & 0 & 0 & $\beta_1$ & $\beta_2$ & 0 & 0 & $\beta_1$ & $\beta_2$ & 0 & 0 & \dots
\end{tabular}
$$
for some $\a_1,\a_2,\beta_1,\beta_2 \in \C$. Since $q =
\exp\(\frac{\pi i}4\)$, we get $[2]=\sqrt 2$ and $[3]=1$, and it is
easy to see the Macdonald-Ruijsenaars equation $\widetilde{\mathbb
M}_1 \varphi(\delta) = [2] \, \varphi(\delta)$ simplifies to
$$
\begin{tabular}{c|c|c}
    $\delta=1$ & $\delta=2$ & $\delta=3$ \\
    \hline
    $2\a_1 = \a_2 - \beta_1$ & $2\a_2 = \a_1 - \beta_2, \ \ 2\beta_1 = \beta_2-\a_1$ & $2\beta_2 = \beta_1-\a_2$
\end{tabular}
$$

These equations have a unique up to proportionality solution, given
by $\a_1 = \a_2 = 1$, $\beta_1 = \beta_2 = -1$, which coincides with
the discrete conformal block $\tilde\cb^{(0,1)}(\delta)$. Therefore,
our Conjecture \ref{conj:conformal root} holds in this case. \qed

\end{document}